\documentclass[a4paper,12pt]{report}
\usepackage{amsmath,amssymb,bm}
\usepackage{amsthm}
\usepackage{amsfonts}
\usepackage{graphicx}
\usepackage[all]{xy}
\usepackage{a4wide}
\usepackage{enumerate}
\usepackage{calrsfs}
\usepackage[latin1]{inputenc}
\usepackage[T1]{fontenc}
\usepackage[francais]{babel}
\usepackage{ae,aeguill}

\DeclareMathOperator\pr{pr}
\DeclareMathOperator*\card{Card\,}
\DeclareMathOperator*\p{\mathsf{P}}
\newcommand{\N}{\mathbb{N}}

\newtheorem{theorem}{Théorème}
\newtheorem{proposition}[theorem]{Proposition}
\newtheorem{lemma}[theorem]{Lemme}
\theoremstyle{definition}
\newtheorem{definition}[theorem]{Définition}
\newtheorem{remark}[theorem]{Remarque}
\newtheorem{remarks}[theorem]{Remarques}
\newtheorem{example}[theorem]{Exemple}
\newtheorem{corollary}[theorem]{Corollaire}
\newtheorem{examples}[theorem]{Exemples}

\begin{document}

\begin{center}

{\LARGE N. Bourbaki, Théorie des ensembles, Hermann 1970}

\vspace{5cm}
{\LARGE \textbf{Notes et Solutions de Quelques Exercices}}

\vspace{2cm}
{\Large Mohssin Zarouali}

\vspace{0.5cm}
(email: mohssin.zarouali@hotmail.com)
\end{center}
\thispagestyle{empty}

%\tableofcontents

\chapter*{Introduction}

L'étude de la théorie des catégories m'a poussé à s'intéresser dernièrement à la théorie des ensembles et plus particulièrement à la notion d'ensemble. C'est pour cette raison que j'ai étudié une bonne partie du livre ``Théorie des Ensembles'' de N. Bourbaki. Il bien connu que G. Cantor est le fondateur de la théorie des ensembles. Je rappelle que D. Hilbert s'est servi de son opération $\varepsilon$ et de son axiome transfini \cite[S5 (I, p. 33)]{Bourbaki:1970} pour fonder les mathématiques \cite{Hilbert/Bernays:1939} (cf. aussi \cite{Froidevaux:1983}), et principalement dans la théorie de la démonstration dont il est l'iniciateur. Ensuite N. Bourbaki a adopté l'opération épsilon de Hilbert pour fonder sa théorie des ensembles, et il lui a préféré la lettre grecque $\tau$ parce qu'il a trouvé que la lettre de Hilbert est trop usuelle en mathématique. L'obligation d'utiliser cette opération en métamathématiques est très contestée. Bourbaki lui-même est conscient qu'on peut se passer d'elle.

Cette opération permet de donner des définitions convaincantes des quantificateurs existentiel et universel, d'ensemble, de relation, de théorème, de démonstration\dots. Elle permet aussi de caractériser les ensembles et les relations (cf. \cite[l'appendice du chapitre I]{Bourbaki:1970}).

La théorie des ensembles de Bourbaki a été adopté par R. Godement dans son cours d'algèbre \cite{Godement:1997}. Je reprends des phrases de ce dernier parce qu'elles décrivent très bien l'opération de Hilbert: - (\dots) l'intérêt de l'opération de Hilbert est de donner un procédé parfaitement artificiel mais purement mécanique pour construire effectivement un objet dont on sait \emph{seulement} qu'il satisfait à des conditions imposées d'avance (dans le cas où de tels objets existerait). - Comme le Dieu des philosophes, l'opération de Hilbert est incompréhensible et ne se voit pas; mais elle gouverne tout, et ses manifestations sensibles éclatent partout \cite[p. 40]{Godement:1997}.

A. Grothendieck lui aussi a fait référence à cette opération dans \cite{Grothendieck:1957} et l'a qualifié du symbole à tout faire $\tau$ de Hilbert. D'ailleurs il s'est situé dans la théorie des ensembles de Bourbaki afin de définir ses univers (cf. \cite{Grothendieck/Verdier:1963/64} et \cite{Grothendieck:1965}).

J'ai rédigé les solutions de quelques exercices qui ont suscités mon intérêt, et j'ai donné dans les annexes A, B et C, des notes concernant les trois premiers chapitres dans le but de les rendre plus accessibles. Par contre, je ne suis pas en accord avec la totalité du livre, par exemple ses définitions des nombres cardinaux, des nombres ordinaux et des nombres entiers naturels. À mon avis, N. Bourbaki a utilisé excessivement l'opération $\bm{\tau}$ de Hilbert pour définir le cardinal d'un ensemble (et aussi pour définir un ordinal, cf. III, p. 77, exerc. 14), ce qui l'a conduit à définir le terme 1 comme étant le cardinal de l'ensemble $\{\varnothing\}$. D'après N. Bourbaki (\cite[III, p. 24]{Bourbaki:1970}):

Une estimation grossière montre que le terme ainsi \emph{désigné} est un assemblage de plusieurs dizaines de milliers de signes (chacun de ces signes étant l'un des signes $\bm{\tau}$, $\square$, $\vee$, $\neg$, =, $\in$).

Vu le nombre estimé de signes figurant dans l'assemblage définissant 1, Yu. I. Manin \cite[p. 17]{Manin:1977} parle d'imprudence de N. Bourbaki et dit que cela semble trop pour 1.

En outre A. R. D. Mathias montre dans \cite{Mathias:2002} que le nombre de signes figurant dans l'assemblage définissant 1 est 4 523 659 424 929 avec 1 179 618 517 981 liens. Ce qui est vraiment trop pour 1.

Je signale que R. Godement a adopté dans \cite{Godement:1997} la définition du cardinal d'un ensemble et celle des entiers naturels de N. Bourbaki.

Afin de contourner ce sérieux problème, j'ai adopté les ordinaux de von Neumann et la définition du cardinal d'un ensemble $X$ comme étant le plus petit ordinal équipotent à $X$, cf. Annexe C.

Je rappelle que le successeur d'un ensemble $X$ est $X^+=X\cup\{X\}$; et les entiers naturels de von Neumann sont: $0=\varnothing$, $1=\varnothing^+=\{0\}$, $2=1^+=\{0,1\}$, $3=2^+=\{0,1,2\}$,\dots.

Je me place maintenant dans la théorie des ensembles de Bourbaki. On sait que l'ensemble vide contient $12$ signes et $3$ liens (cf. \cite[II, p. 6]{Bourbaki:1970}). J'ai calculé les nombres de signes et de liens figurant dans 1 et 2. Notre 1 contient $513$ signes et $134$ liens. Ce qui est raisonnable pour 1. Notre 2 contient $7527$ signes et $1968$ liens.

Ces notes constituent une occasion pour découvrir la beauté de la théorie des ensembles de Bourbaki, et j'espère qu'elles faciliteront la tâche à ceux qui ont trouvé des difficultés en lisant ce livre qui est l'un des plus critiqué des éléments de mathématique de Bourbaki.

\chapter*{Solutions de quelques exercices du chapitre I (Description de la mathématique formelle)}

\begin{center}
\S 1
\end{center}

1) Découle du fait que les éléments de chaque construction formative de $\mathcal{T}$ sont des lettres car il n'y a pas d'éléments du type e) y appartenant (voir I, p. 18).

\vspace{0.5cm}
2) Découle du fait que l'opération faisant intervenir les signes $\bm{\tau}$, $\square$ et les liens dans chaque construction formative est d) (voir I, p. 18).

\vspace{0.5cm}
3) $\bm{A}$ figure dans une construction formative. Soit $\bm{s}$ un signe spécifique figurant dans $\bm{A}$. La première occurence de $\bm{s}$ est dans une relation de type e) (voir I, p. 18) et il est suivi donc d'une lettre ou de $\bm{\tau}$. Si dans la suite on a une relation $\bm{R}$ où figure $\bm{s}$ qui est suivi d'une lettre $\bm{x}$ et $\bm{\tau}_{\bm{x}}(\bm{R})$ figure dans la construction formative, $\bm{s}$ est donc suivi par $\square$. D'où le résultat.

\vspace{0.5cm}
4) On propose une méthode plus simple que celle suggérée. On a $(\bm{A}\bm{B})^*$ est $\bm{A}^*\bm{B}^*$ qui n'est pas équilibré car $\bm{A}$ est équilibré (Critère 1, I, p. 45). Donc $\bm{A}\bm{B}$ n'est ni un terme ni une relation.

\vspace{0.5cm}
5) On a $\bm{\tau}_{\bm{x}}(\bm{A})$ est identique à $\bm{\tau}_{\bm{y}}(\bm{R})$ où $\bm{R}$ une relation et  $\bm{y}$ une lettre.  Si $\bm{x}$ figure dans $\bm{A}$, donc $\bm{y}$ figure dans $\bm{R}$. Donc $\bm{A}$ est identique à $(\bm{x}\mid\bm{y})\bm{R}$. D'après CF7, $\bm{A}$ est une relation. Si $\bm{x}$ ne figure pas dans $\bm{A}$, donc $\bm{y}$ ne figure pas dans $\bm{R}$. D'où $\bm{A}$ est identique à $\bm{R}$.

\vspace{0.5cm}
6) On propose de montrer un résultat plus général et de donner une méthode plus simple que celle suggérée. Soient $\bm{A}$ et $\bm{B}$ des assemblages d'une théorie $\mathcal{T}$ tels que $\bm{A}$ et $\vee\bm{A}\bm{B}$ sont des relations de $\mathcal{T}$, $\bm{B}$ est une relation de $\mathcal{T}$. D'abord on écrit la relation $\vee\bm{A}\bm{B}$ sous la forme $\vee\bm{A}'\bm{B}'$ où $\bm{A}$ et $\bm{B}$ sont deux relations. Puisque les assemblages $\vee\bm{A}\bm{B}$, $\bm{A}$ et $\bm{A}'$ sont équilibrés (Critère 1, I, p. 45) et le mot $(\vee\bm{A}\bm{B})^*$ est identique à $\vee\bm{A}^*\bm{B}^*$ et à $\vee\bm{A}'^*\bm{B}'^*$, $\bm{A}^*$ et $\bm{B}^*$ sont identiques à $\bm{A}'^*$ et $\bm{B}'^*$ respectivement. D'où $\bm{B}$ est identique à $\bm{B}'$, donc une relation.

\vspace{2cm}
\begin{center}
\S 2
\end{center}

1) a) On a $\mathcal{T}$ est plus forte que $\mathcal{T}'$. D'après C4, pour que $\bm{A}_n$ ne soit pas indépendant des autres axiomes de $\mathcal{T}$ (i.e. $\mathcal{T}'$ soit plus forte que $\mathcal{T}$) il faut et il suffit que $\bm{A}_n$ soit un théorème $\mathcal{T}'$.

b) Supposons que $\bm{A}_n$ est un théorème $\mathcal{T}'$. D'après C5 ($\mathcal{T}'$ au lieu de $\mathcal{T}$, $\mathcal{T}''$ au lieu de $\mathcal{T}'$ et noter que les constantes de $\mathcal{T}'$ sont parmis $\bm{a}_1,\bm{a}_2,\dots,\bm{a}_h$), $$(\bm{T}_1\mid\bm{a}_1)(\bm{T}_2\mid\bm{a}_2)\dots(\bm{T}_h\mid\bm{a}_h)\bm{A}_n$$ est un théorème de $\mathcal{T}''$. D'où $\mathcal{T}''$ est contradictoire.

\vspace{2cm}
\begin{center}
\S 3
\end{center}

1) D'abord la 6 ème relation de l'ex. 1) est un théorème même si on remplace ``$\Rightarrow$'' par ``$\Leftrightarrow$''. On va montrer ses théorèmes mais pas dans l'ordre suggéré.

\vspace{0.5cm}
\textbf{Dans la suite de ses solutions d'exercices on va utiliser C1, C6 et C22 sans les mentioner, et on se servira des théorèmes 1) à 18) une fois démontrés.}
\vspace{0.5cm}

Soient $\bm{A}$, $\bm{B}$, $\bm{C}$ et $\bm{D}$ des relations d'une théorie logique $\mathcal{T}$.

\begin{enumerate}[(1)]
\item \label{impl-ou} D'après C14, C21, S4 et S2,
$$(\bm{A}\Rightarrow\bm{B}\textrm{ et }\bm{C}\Rightarrow\bm{D})\Rightarrow\big((\bm{A}\textrm{ ou }\bm{C})\Rightarrow(\bm{B}\textrm{ ou }\bm{D})\big)$$
est un théorème.
\item \label{equiv-ou} D'après C14, C21, \eqref{impl-ou} et C20,
$$(\bm{A}\Leftrightarrow\bm{B}\textrm{ et }\bm{C}\Leftrightarrow\bm{D})\Rightarrow\big((\bm{A}\textrm{ ou }\bm{C})\Leftrightarrow(\bm{B}\textrm{ ou }\bm{D})\big)$$
est un théorème.
\item \label{impl-et} D'après C14, C21, C12 et \eqref{impl-ou},
$$(\bm{A}\Rightarrow\bm{B}\textrm{ et }\bm{C}\Rightarrow\bm{D})\Rightarrow\big((\bm{A}\textrm{ et }\bm{C})\Rightarrow(\bm{B}\textrm{ et }\bm{D})\big)$$
est un théorème.
\item \label{equiv-et} D'après C14, C21, \eqref{impl-et} et C20,
$$(\bm{A}\Leftrightarrow\bm{B}\textrm{ et }\bm{C}\Leftrightarrow\bm{D})\Rightarrow\big((\bm{A}\textrm{ et }\bm{C})\Leftrightarrow(\bm{B}\textrm{ et }\bm{D})\big)$$
est un théorème.

\item \label{C25} D'après C8 et C20, $\bm{C}\Leftrightarrow\bm{C}$ est un théorème. Si $\bm{B}$ est un théorème, d'après C21 et C20 et C14, $\bm{A}\Leftrightarrow(\bm{A}\textrm{ et }\bm{B})$ est aussi un théorème. Donc d'après \eqref{equiv-non}, \eqref{nonnon} (cf. ci-dessous) et \eqref{equiv-ou}, si $\textrm{non }\bm{B}$ est un théorème, $\bm{A}\Leftrightarrow(\bm{A}\textrm{ ou }\bm{B})$ est un théorème.

    Donc on a les théorèmes suivants:
\item \label{impl-ou'} $(\bm{A}\Rightarrow\bm{B})\Rightarrow\big((\bm{A}\textrm{ ou }\bm{C})\Rightarrow(\bm{B}\textrm{ ou }\bm{C})\big),$
\item \label{equiv-ou'} $(\bm{A}\Leftrightarrow\bm{B})\Rightarrow\big((\bm{A}\textrm{ ou }\bm{C})\Leftrightarrow(\bm{B}\textrm{ ou }\bm{C})\big),$
\item \label{impl-et'} $(\bm{A}\Rightarrow\bm{B})\Rightarrow\big((\bm{A}\textrm{ et }\bm{C})\Rightarrow(\bm{B}\textrm{ et }\bm{C})\big),$
\item \label{equiv-et'} $(\bm{A}\Leftrightarrow\bm{B})\Rightarrow\big((\bm{A}\textrm{ et }\bm{C})\Leftrightarrow(\bm{B}\textrm{ et }\bm{C})\big).$

\item \label{AouA} D'après S1, S2 et C20, $\bm{A}\Leftrightarrow(\bm{A}\textrm{ ou }\bm{A})$ est un théorème. Donc

\item d'après C14 et \eqref{impl-ou},
 \label{impl-ou''} $$(\bm{A}\Rightarrow\bm{C})\Rightarrow\big((\bm{B}\Rightarrow\bm{C})\Rightarrow((\bm{A}\textrm{ ou
    }\bm{B})\Rightarrow\bm{C})\big)$$ est un théorème.

\item \label{nonnon} D'après  C11, C16 et C20, $\bm{A}\Leftrightarrow(\textrm{non }\textrm{non }\bm{A})$ est un théorème. Donc

\item d'après C14, C12 et \eqref{impl-ou''},
 \label{impl-et''} $$(\bm{A}\Rightarrow\bm{B})\Rightarrow\big((\bm{A}\Rightarrow\bm{C})\Rightarrow(\bm{A}\Rightarrow(\bm{B}\textrm{ et }\bm{C})\big)$$ est un théorème.
\item \label{comm-ou} D'après S3 et C20, $(\bm{A}\textrm{ ou }\bm{B})\Leftrightarrow(\bm{B}\textrm{ ou }\bm{A})$ est un théorème.
\item D'après S2, C7, \eqref{impl-ou''} et C20, $\bm{A}\textrm{ ou }(\bm{B}\textrm{ ou }\bm{C})\Leftrightarrow(\bm{A}\textrm{ ou }\bm{B})\textrm{ ou }\bm{C}$ est un théorème. Donc si on interchange $\bm{A}$, $\bm{B}$, $\bm{C}$ et les parenthèses dans ``$\bm{A}\textrm{ ou }\bm{B}\textrm{ ou }\bm{C}$'', qui est par définition ``$\bm{A}\textrm{ ou }(\bm{B}\textrm{ ou }\bm{C})$'', on obtient une relation équivalente (on se sert de \eqref{comm-ou}).

    D'où la 7 ème et la 8 ème relations de l'ex. 1) sont des théorèmes.
\item \label{dist-et-ou} On propose de montrer que $$(\bm{A}\textrm{ et }(\bm{B}\textrm{ ou }\bm{C}))\Leftrightarrow((\bm{A}\textrm{ et }\bm{B})\textrm{ ou }(\bm{A}\textrm{ et }\bm{C}))$$ est un théorème. ``$\Rightarrow$'': Soit $\mathcal{T}'$ la théorie obtenue en adjoignant $\bm{A}\textrm{ et }(\bm{B}\textrm{ ou }\bm{C})$ aux axiomes explicites. $\bm{A}$, $\bm{B}\textrm{ ou }\bm{C}$ sont des théorèmes de $\mathcal{T}'$. D'après C14, S2 et C7, $\bm{B}\Rightarrow((\bm{A}\textrm{ et }\bm{B})\textrm{ ou }(\bm{A}\textrm{ et }\bm{C}))$, $\bm{C}\Rightarrow((\bm{A}\textrm{ et }\bm{B})\textrm{ ou }(\bm{A}\textrm{ et }\bm{C}))$ sont des théorèmes de $\mathcal{T}'$. D'après \eqref{impl-ou'}, $(\bm{A}\textrm{ et }\bm{B})\textrm{ ou }(\bm{A}\textrm{ et }\bm{C})$ est un théorème de $\mathcal{T}'$. C14 montre donc l'implication. ``$\Leftarrow$'': D'après S2 et C7, $\bm{B}\Rightarrow(\bm{B}\textrm{ ou }\bm{C})$, $\bm{C}\Rightarrow(\bm{B}\textrm{ ou }\bm{C})$ sont des théorèmes. Donc d'après \eqref{impl-et}, $(\bm{A}\textrm{ et }\bm{B})\Rightarrow(\bm{A}\textrm{ et }(\bm{B}\textrm{ ou }\bm{C}))$, $(\bm{A}\textrm{ et }\bm{C})\Rightarrow(\bm{A}\textrm{ et }(\bm{B}\textrm{ ou }\bm{C}))$ sont des théorèmes. \eqref{impl-ou''} montre donc cette implication.

    \item \label{equiv-non} D'après C21, C12 et C20, si $\bm{A}\Leftrightarrow\bm{B}$ est un théorème, il en est de même de $(\textrm{non }\bm{A})\Leftrightarrow(\textrm{non }\bm{B})$.
    \item \label{dist-ou-et} En rempla\c{c}ant dans \eqref{dist-et-ou} $\bm{A}$, $\bm{B}$ et $\bm{C}$ par $\textrm{non }\bm{A}$, $\textrm{non }\bm{B}$ et $\textrm{non }\bm{C}$ respectivement, et en utilisant \eqref{equiv-non}, \eqref{nonnon}, \eqref{equiv-ou} et \eqref{equiv-et}, on montre le théorème
    $$(\bm{A}\textrm{ ou }(\bm{B}\textrm{ et }\bm{C}))\Leftrightarrow((\bm{A}\textrm{ ou }\bm{B})\textrm{ et }(\bm{A}\textrm{ ou }\bm{C})).$$
    \item D'après \eqref{dist-ou-et} et \eqref{equiv-et}, $(\bm{A}\textrm{ et }\bm{B})\textrm{ ou }((\textrm{non }\bm{A})\textrm{ et }(\textrm{non }\bm{B}))$ est équivalente \'{a}

    $\big((\bm{A}\textrm{ ou }\textrm{non }\bm{A})\textrm{ et }(\bm{B}\textrm{ ou }\textrm{non }\bm{A})\big) \textrm{ et } \big((\bm{A}\textrm{ ou }\textrm{non }\bm{B})\textrm{ et }(\bm{B}\textrm{ ou }\textrm{non }\bm{B})\big)$. D'après C10, C25 et \eqref{equiv-et}, cette dernière relation est équivalente \'{a} $(\bm{B}\textrm{ ou }\textrm{non }\bm{A})\textrm{ et }(\bm{A}\textrm{ ou }\textrm{non }\bm{B})$ qui à son tour équivalente \'{a} $\bm{A}\Leftrightarrow\bm{B}$. D'où la 5 ème relation de l'ex. 1) est un théorème.

D'après \eqref{equiv-ou'}, $(\textrm{non }\bm{A}\Rightarrow\bm{B})\Leftrightarrow(\bm{A}\textrm{ ou }\bm{B})$ est un théorème. Donc, d'après \eqref{equiv-et}, $$((\textrm{non }\bm{A})\Leftrightarrow\bm{B})\Leftrightarrow\big((\bm{A}\textrm{ ou }\bm{B})\textrm{ et }((\textrm{non }\bm{A})\textrm{ ou }(\textrm{non }\bm{B}))\big).$$ Donc $$(\textrm{non }((\textrm{non }\bm{A})\Leftrightarrow\bm{B}))\Leftrightarrow (\textrm{non }(\bm{A}\textrm{ ou }\bm{B})\textrm{ ou }(\bm{A}\textrm{ et }\bm{B})).$$
D'après \eqref{nonnon}, \eqref{equiv-ou} et \eqref{equiv-non}, $((\textrm{non }\bm{A})\textrm{ et }(\textrm{non }\bm{B}))\Leftrightarrow\textrm{non }(\bm{A}\textrm{ ou }\bm{B})$. En se servant de \eqref{equiv-ou'} et du 5 ème théorème de l'ex. 1) on montre que la 6 ème relation de l'ex. 1) est un théorème.

Maintenant on montre que les quatre premières relations de l'ex. 1) sont des théorèmes.
    \item $\bm{B}\Rightarrow\bm{A}$ est identique \'{a} $(\textrm{non }\bm{B})\textrm{ ou }\bm{A}$. On utilise C7.
\item On utilise C14 deux fois.
\item  $(\textrm{non }\bm{A})\Rightarrow\bm{B}$ est identique \'{a} $(\textrm{non }\textrm{non }\bm{A})\textrm{ ou }\bm{B}$. On utilise C11 et S2.
    \item D'après \eqref{nonnon}, \eqref{AouA}, \eqref{equiv-ou} et \eqref{equiv-non}, $\textrm{non }(\bm{A}\Rightarrow\bm{B})\Leftrightarrow (\bm{A}\textrm{ et }(\textrm{non }\bm{B}))$ est un théorème. On utilise ensuite \eqref{equiv-ou'}, \eqref{comm-ou}, \eqref{dist-ou-et}, C10 et \eqref{C25}.
\end{enumerate}

\vspace{0.5cm}
2) D'abord la relation $\bm{A}\Rightarrow(\textrm{non }\bm{A})$ n'est autre que ``$(\textrm{non }\bm{A})\textrm{ ou }(\textrm{non }\bm{A})$'', qui est équivalente d'après \eqref{AouA} \'{a} $\textrm{non }\bm{A}$. D'après \eqref{nonnon} et \eqref{equiv-ou'}, la relation $(\textrm{non }\bm{A})\Rightarrow\bm{A}$ est équivalente \'{a} ``$\bm{A}\textrm{ ou }\bm{A}$''. Or celle-ci est équivalente d'après \eqref{AouA} \'{a} $\bm{A}$. On conclut en utilisant C21.

\vspace{0.5cm}
3) a) On le montre par récurrence sur $n$. Cas $n=2$: D'après C14, $(\textrm{non }\bm{A}_1)\Rightarrow\bm{A}_2$ est un théorème. D'après \eqref{impl-ou}, $(\bm{A}_1\textrm{ ou }\textrm{non }\bm{A}_1)\Rightarrow(\bm{A}_1\textrm{ ou }\bm{A}_2)$ est aussi un théorème. D'après C10, il en est de même de $\bm{A}_1\textrm{ ou }\bm{A}_2$. Soient $\mathcal{T}'$ (resp. $\mathcal{T}''$) la théorie obtenue en adjoignant \'{a} $\mathcal{T}$ les axiomes $\textrm{non }\bm{A}_1$, $\textrm{non }\bm{A}_2$, \dots, $\textrm{non }\bm{A}_{n-1}$ (resp. $\textrm{non }\bm{A}_1$). Supposons que $\bm{A}_n$ est un théorème dans $\mathcal{T}'$. Par hypothèse de récurrence $\bm{A}_2\textrm{ ou }\dots\textrm{ ou }\bm{A}_n$ est un théorème dans $\mathcal{T}''$. D'après le cas où $n=2$, $\bm{A}_1\textrm{ ou }\dots\textrm{ ou }\bm{A}_n$ est un théorème dans $\mathcal{T}$.

b) est immédiat par application de \eqref{impl-ou''}.

\vspace{0.5cm}
4) évident.

\vspace{0.5cm}
5) D'abord soit $\mathcal{T}$ une théorie logique contradictoire. Soit $\bm{A}$ une relation telle que $\bm{A}$, $\textrm{non }\bm{A}$ sont des théorèmes. Donc toute relation $\bm{B}$ de $\mathcal{T}$ est un théorème. En effet, d'après S2, $(\textrm{non }\bm{A})\Rightarrow(\bm{A}\Rightarrow\bm{B})$ est un théorème.

Maintenant on montre l'ex. 5). Soit $\mathcal{T}'$ (resp. $\mathcal{T}''$) la théorie dont les signes et les schémas sont ceux de $\mathcal{T}$, et les axiomes explicites sont $\bm{A}_1$, $\bm{A}_2$, \dots, $\bm{A}_{n-1}$ (resp. $\bm{A}_1$, $\bm{A}_2$, \dots, $\bm{A}_{n-1}$, $(\textrm{non }\bm{A}_n)$). D'après l'ex. 1), \S2, pour que $\bm{A}_n$ ne soit pas indépendant des autres axiomes de $\mathcal{T}$ il faut et il suffit que $\bm{A}_n$ soit un théorème de $\mathcal{T}'$. Si $\mathcal{T}''$ est contradictoire, $\bm{A}_n$ est un théorème de $\mathcal{T}''$, donc, d'après, C14, \eqref{impl-ou''} et C10, $\bm{A}_n$ est un théorème $\mathcal{T}'$. Finalement, si $\bm{A}_n$ est un théorème $\mathcal{T}'$, il est aussi un théorème $\mathcal{T}''$ car $\mathcal{T}''$ est plus forte que $\mathcal{T}'$. Or $\textrm{non }\bm{A}_n$ est un théorème $\mathcal{T}''$. D'où $\mathcal{T}''$ est contradictoire.

\vspace{2cm}
\begin{center}
\S 4
\end{center}

1) Immédiat d'après C33.

\vspace{0.5cm}
2) Immédiat d'après C31.

\vspace{0.5cm}
Comme dans le paragraphe 3 du \S4, soit $\mathcal{T}_0$ la théorie sans axiomes explicites qui possède les mêmes signes que $\mathcal{T}$ et les seuls schémas S1 à S5. $\mathcal{T}$ est plus forte que $\mathcal{T}_0$. Il suffit donc de montrer les théorèmes dans les exercices 3) à 6) du \S4 dans $\mathcal{T}_0$.

\vspace{0.5cm}
3) Le premier théorème est immédiat d'après C30 et C31, le second d'après S5 et C31.

\vspace{0.5cm}
4) D'après S5 et S4, $(\bm{A}\textrm{ ou }\bm{B})\Rightarrow(\bm{A}\textrm{ ou }(\exists\bm{x})\bm{B})$ est un théorème. Le premier théorème est donc immédiat d'après C31 et C33. On montre le second théorème soit directement comme pour le premier, soit en utilisant ce dernier et C12, C29, C28, \eqref{equiv-ou}, \eqref{equiv-non} respectivement.

\vspace{0.5cm}
5) Immédiat d'après C32 et C31.

\vspace{0.5cm}
6) Le premier théorème est immédiat d'après C21 et C31, le second d'après le premier théorème et C12.

\vspace{0.5cm}
7) D'après l'ex. 6), pour montrer le premier théorème il suffit de montrer la première implication, or celle-ci  est immédiate d'après \eqref{impl-et''} et C31. Le second théorème est immédiat en utilisant le premier et \eqref{equiv-non}. Le cas particulier découle de C7.

\vspace{0.5cm}
8) On montre le premier théorème par C14. Soit $\mathcal{T}'$ la théorie obtenue en adjoignant $(\bm{T}\mid\bm{x})\bm{R}$ aux axiomes. D'après CS6 et C20, $(\bm{T}\mid\bm{x})(\bm{A}\textrm{ et }\bm{R})$ est un théorème dans $\mathcal{T}'$. D'après S5, $(\exists\bm{x})(\bm{A}\textrm{ et }\bm{R})$ est aussi un théorème dans $\mathcal{T}'$. Le deuxième découle du premier en utilisant C12, CS5 et \eqref{nonnon}.

\vspace{2cm}
\begin{center}
\S 5
\end{center}

1) Plus généralement, soient $\bm{T}$ un terme et $\bm{x}$ une lettre ne figurant pas dans $\bm{T}$. La relation ``$\bm{x}=\bm{T}$'' est fonctionnelle en $\bm{x}$ dans $\mathcal{T}$. D'abord ``$\bm{T}=\bm{T}$'' (voir la généralisation du Th. 1, I, p. 39), donc d'après S5, $(\exists \bm{x})(\bm{x}=\bm{T})$ est vraie. Pour voir que la relation ``$\bm{x}=\bm{T}$'' est univoque en $\bm{x}$ il suffit d'utiliser la généralisation du Th. 3, I, p. 40.

\vspace{0.5cm}
2) Plus généralement, soient $\bm{R}$ une relation, $\bm{T}$ un terme et $\bm{x}$ une lettre ne figurant pas dans $\bm{T}$. On a $(\exists \bm{x})(\bm{x}=\bm{T}\textrm{ et } \bm{R})$, $(\bm{T}\mid\bm{x})\bm{R}$ sont équivalentes dans $\mathcal{T}$. D'après l'ex. 1) et C47, $(\exists \bm{x})(\bm{x}=\bm{T}\textrm{ et }\bm{R})$ est équivalente à $(\bm{\tau}_{\bm{x}}(\bm{x}=\bm{T})\mid\bm{x})\bm{R}$. D'après C46, $(\bm{x}=\bm{T})\Leftrightarrow (\bm{x}=\bm{\tau}_{\bm{x}}(\bm{x}=\bm{T}))$ est un théorème (dans la première partie de C46 et de C47 la condition ``$\bm{x}$ n'est pas une constante'' n'est pas nécessaire). Comme dans la démonstration de l'ex. 1), $(\exists \bm{x})(\bm{x}=\bm{T})$ est vraie, or ce n'est que $\bm{\tau}_{\bm{x}}(\bm{x}=\bm{T})=\bm{T}$. D'après S6, $(\bm{\tau}_{\bm{x}}(\bm{x}=\bm{T})\mid\bm{x})\bm{R}$, $(\bm{T}\mid\bm{x})\bm{R}$ sont équivalentes.

\vspace{0.5cm}
3) D'après C14, $\bm{S}\Rightarrow (\exists\bm{x})\bm{S}$ est un théorème de $\mathcal{T}$. D'après C3 et CS5, $(\bm{T}\mid\bm{y})\bm{S}\Rightarrow (\bm{T}\mid\bm{y})(\exists\bm{x})\bm{R}$. D'après CS9, $(\bm{T}\mid\bm{y})(\exists\bm{x})\bm{R}$ est identique \'{a} $(\exists\bm{x})(\bm{T}\mid\bm{y})\bm{R}$. Donc $(\exists\bm{x})(\bm{T}\mid\bm{y})\bm{R}$ est un théorème de $\mathcal{T}$.

Soient $\bm{z}$ et $\bm{z}'$ deux lettres distinctes entre elles, distinctes de $\bm{x}$, de $\bm{y}$ et ne figurant pas dans $\bm{R}$ ni dans $\bm{T}$ (donc non plus dans $(\bm{T}\mid\bm{y})\bm{R}$). D'après C14,
$$\bm{S}\Rightarrow(\forall\bm{z})(\forall\bm{z}')((\bm{z}\mid\bm{x})\bm{R}\textrm{ et }(\bm{z}'\mid\bm{x})\bm{R}\Rightarrow (\bm{z}=\bm{z}'))$$ est un théorème de $\mathcal{T}$. De même, d'après C3, CS5, CS9 et CS2, $$(\bm{T}\mid\bm{y})\bm{S}\Rightarrow(\forall\bm{z})(\forall\bm{z}')((\bm{z}\mid\bm{x})(\bm{T}\mid\bm{y})\bm{R}
\textrm{ et }(\bm{z}'\mid\bm{x})(\bm{T}\mid\bm{y})\bm{R}\Rightarrow (\bm{z}=\bm{z}'))$$ est un théorème de $\mathcal{T}$.

\textbf{Remarque.} Soient $\bm{T}$ un terme et $\bm{x}$ une lettre ne figurant pas dans $\bm{T}$. On prend $\bm{y}$ une lettre distincte de $\bm{x}$, n'est pas une constante et ne figurant pas dans $\bm{T}$. Soient $\bm{R}$ la relation ``$\bm{x}=\bm{y}$'' et $\bm{S}$ la relation ``$\bm{y}=\bm{T}$''. L'application de l'ex. 3) donne notre généralisation de l'ex. 1).

\vspace{0.5cm}
4) D'après C31, $(\exists\bm{x})\bm{S}$ est vraie. Soient $\bm{y}$ et $\bm{z}$ deux lettres distinctes entre elles, distinctes de $\bm{x}$, distinctes des constantes et ne figurant pas dans $\bm{R}$ ni dans $\bm{S}$. D'après C3 et CS7, $(\bm{y}\mid\bm{x})\bm{R}\Leftrightarrow (\bm{y}\mid\bm{x})\bm{S}$, $(\bm{z}\mid\bm{x})\bm{R}\Leftrightarrow (\bm{z}\mid\bm{x})\bm{S}$ sont des théorèmes.

D'après \eqref{equiv-et}, $(\bm{y}\mid\bm{x})\bm{R}\textrm{ et }(\bm{z}\mid\bm{x})\bm{R}$, $(\bm{y}\mid\bm{x})\bm{S}\textrm{ et }(\bm{z}\mid\bm{x})\bm{S}$ sont équivalentes. D'après C21, $$\big((\bm{y}\mid\bm{x})\bm{S}\textrm{ et }(\bm{z}\mid\bm{x})\bm{S}\big)\Rightarrow\big((\bm{y}\mid\bm{x})\bm{R}\textrm{ et }(\bm{z}\mid\bm{x})\bm{R}\big)$$ est un théorème. D'après C30 et C6, $$\big((\bm{y}\mid\bm{x})\bm{S}\textrm{ et }(\bm{z}\mid\bm{x})\bm{S}\big)\Rightarrow (\bm{y}=\bm{z}).$$ On finit en utilisant C27.

\vspace{0.5cm}
5) D'abord soient $\bm{x}$ une lettre, $\bm{R}$ une relation fonctionnelle en $\bm{x}$ et $\bm{x}'$ une lettre ne figurant pas dans $\bm{R}$. D'après CS1 et CS8, $(\bm{x}'\mid\bm{x})\bm{R}$ est fonctionnelle en $\bm{x}'$.

Donc on peut supposer que $\bm{x}$ n'est pas une constante quitte à considérer une lettre $\bm{x}'$ qui n'est pas une constante et ne figurant pas dans $\bm{R}$ et dans $\bm{S}$ (et dans $\bm{T}$ en ce qui concerne les deux derniers théorèmes) (on utilise CS8, CS6 et CS5). Par exemple la première relation est identique \'{a}
$$\big(\textrm{non }(\exists\bm{x}')((\bm{x}'\mid\bm{x})\bm{R}\textrm{ et }(\bm{x}'\mid\bm{x})\bm{S})\big)\Leftrightarrow (\exists\bm{x}')\big((\bm{x}'\mid\bm{x})\bm{R}\textrm{ et }(\textrm{non }(\bm{x}'\mid\bm{x})\bm{S})\big).$$

Le premier théorème est immédiat en utilisant C47, \eqref{nonnon} et CS5.

Le second théorème est immédiat en utilisant C47, CS6 et \eqref{impl-et'}.

Le dernier théorème est immédiat en utilisant C47, CS5 et \eqref{equiv-ou}.

\vspace{0.5cm}
6) Soient $\bm{x}$ et $\bm{y}$ deux lettres distinctes, $\bm{R}$ une relation où figurent $\bm{x}$ et $\bm{y}$, $\bm{A}$ la relation $(\exists \bm{x})\bm{R}\Rightarrow\bm{R}$.  La relation $(\bm{y}\mid\bm{x})\bm{A}$ est identique \'{a} $(\exists \bm{x})\bm{R}\Rightarrow(\bm{y}\mid\bm{x})\bm{R}$ (CS5). Supposons que $(\bm{y}\mid\bm{x})\bm{A}$ est de la forme $(\exists \bm{z})\bm{R}'\Rightarrow\bm{R}'$. Si $\bm{B}$ et $\bm{C}$ deux relations, les assemblages antécédents à $\vee\bm{B}\bm{C}$ sont $\bm{B}$ et $\bm{C}$ (voir I, p. 46). Donc $\bm{R}'$ est identique \'{a} $(\bm{y}\mid\bm{x})\bm{R}$ et $(\exists \bm{z})\bm{R}'$ est identique \'{a} $(\exists \bm{x})\bm{R}$. Si $\bm{z}$ est identique \'{a} $\bm{x}$, on aurait $(\exists \bm{x})\bm{R}$ est identique \'{a} $(\bm{y}\mid\bm{x})\bm{R}$. Donc $((\exists \bm{x})\bm{R})^*$ et $((\bm{y}\mid\bm{x})\bm{R})^*$ aurait la même longeur que $\bm{R}^*$. Par conséquent $\bm{\tau}_{\bm{x}}(\bm{R})$ est identique \'{a} $\bm{y}$. D'où $\bm{z}$ est distincte de $\bm{x}$. Puisque $\bm{z}$ ne figure pas dans $(\exists \bm{z})(\bm{y}\mid\bm{x})\bm{R}$ qui est identique à $(\exists \bm{x})\bm{R}$, elle ne figure non plus dans $\bm{R}$. Donc $(\exists \bm{x})\bm{R}$ est identique \'{a} $(\bm{y}\mid\bm{x})\bm{R}$ et on a vu que cela n'est pas possible. On conclut que $(\bm{y}\mid\bm{x})\bm{A}$ ne peut être de la forme $(\exists \bm{z})\bm{R}'\Rightarrow\bm{R}'$.

\vspace{0.5cm}
7) La méthode est la même que celle de 6). Soient $\bm{x}$ et $\bm{y}$ deux lettres distinctes, $\bm{R}$ une relation où figurent $\bm{x}$ et $\bm{y}$, $\bm{S}$ une relation quelconque, $\bm{A}$ la relation $(\bm{R}\Leftrightarrow\bm{S})\Rightarrow(\bm{\tau}_{\bm{x}}(\bm{R})=\bm{\tau}_{\bm{x}}(\bm{S}))$. La relation $(\bm{y}\mid\bm{x})\bm{A}$ est identique \'{a} $((\bm{y}\mid\bm{x})\bm{R}\Leftrightarrow(\bm{y}\mid\bm{x})\bm{S})\Rightarrow(\bm{\tau}_{\bm{x}}(\bm{R})=
\bm{\tau}_{\bm{x}}(\bm{S}))$ (CS5 et CS7). Supposons que $(\bm{y}\mid\bm{x})\bm{A}$ est de la forme $(\bm{R}'\Leftrightarrow\bm{S}')\Rightarrow(\bm{\tau}_{\bm{z}}(\bm{R}')=\bm{\tau}_{\bm{z}}(\bm{S}'))$. Donc $\bm{R}'$, $\bm{S}'$, $\bm{\tau}_{\bm{z}}(\bm{R}')$, $\bm{\tau}_{\bm{z}}(\bm{S}')$ sont identiques respectivement \'{a} $(\bm{y}\mid\bm{x})\bm{R}$, $(\bm{y}\mid\bm{x})\bm{S}$, $\bm{\tau}_{\bm{x}}(\bm{R})$, $\bm{\tau}_{\bm{x}}(\bm{S})$. Puisque $\bm{\tau}_{\bm{z}}((\bm{y}\mid\bm{x})\bm{R})$ est identique à $\bm{\tau}_{\bm{x}}(\bm{R})$ et $\bm{x}$ figure dans $\bm{R}$ et non dans $(\bm{y}\mid\bm{x})\bm{R}$, $\bm{x}$ est distincte de $\bm{z}$. Puisque $\bm{z}$ ne figure pas dans $\bm{\tau}_{\bm{z}}((\bm{y}\mid\bm{x})\bm{R})$, elle ne figure non plus dans $\bm{R}$ ni donc dans $(\bm{y}\mid\bm{x})\bm{R}$. D'où $(\bm{y}\mid\bm{x})\bm{A}$ ne peut être de la forme $(\bm{R}'\Leftrightarrow\bm{S}')\Rightarrow(\bm{\tau}_{\bm{z}}(\bm{R}')=\bm{\tau}_{\bm{z}}(\bm{S}'))$.

\chapter*{Solutions de quelques exercices du chapitre II (Théorie des ensembles)}

\begin{center}
\S 1
\end{center}

1) L'implication ``$\Rightarrow$'' est imm\'ediate par application de S6 \'a la relation ``$z\in X$''.

L'implication ``$\Leftarrow$''. On utilise la méthode de l'hypothèse auxiliaire (C14, I, p. 27). D'après C31, $(\{x\}\mid X)(x\in X \Rightarrow y\in X)$ est un théorème, or cette relation est identique à ``$x\in \{x\} \Rightarrow y\in \{x\}$''. Or ``$z\in \{x\}\Leftrightarrow z=x$'' est aussi un théorème. Puisque ``$x=x$'' est vraie, il en est de même de ``$x\in \{x\}$'' et ``$y=x$''.

\vspace{0.5cm}
2) On a vu que $(\forall x)(x\notin X)\Leftrightarrow X=\varnothing$ est vraie (II, p. 6). Donc si $\bm{x}$ une lettre ne figurant pas dans un $\bm{T}$, alors $(\forall \bm{x})(\bm{x}\notin \bm{T})\Leftrightarrow \bm{T}=\varnothing$ est aussi vraie, or cette dernière relation est équivalente à équivalente à $(\exists \bm{x})(\bm{x}\in \bm{T})\Leftrightarrow \bm{T}\neq\varnothing$. Puisque $x\in\{x\}$ est vraie (voir 1)), vraie.

D'après C28 et C31, $((\exists x)(\exists y)(x\neq y))\Leftrightarrow (\textrm{non }(\forall x)(\forall y)(x=y))$ est un théorème. Puisque les deux implications $((\forall x)(\forall y)(x=y))\Rightarrow (\forall y)(\{x\}=y)$ et $(\forall y)(\{x\}=y)\Rightarrow \{x\}=\varnothing$ sont vraies et d'après C12, $\{x\}\neq\varnothing \Rightarrow ((\exists x)(\exists y)(x\neq y))$ est aussi vraies.

\vspace{0.5cm}
3) Immédiate en utilisant que si $\bm{A}$ et $\bm{B}$ deux parties  d'un ensemble $\bm{X}$, les deux relations suivantes sont vraies: $\bm{A}=\complement_{\bm{X}}(\complement_{\bm{X}}\bm{A})$, $\bm{A}\subset\bm{B}\Leftrightarrow \complement_{\bm{X}}\bm{B}\subset\complement_{\bm{X}}\bm{A}$ (II, p. 6).

\vspace{0.5cm}
4) L'implication ``$\Leftarrow$'' est imm\'ediate par application de C14, C18 et S6.

L'implication ``$\Rightarrow$''. Supposons que $X\subset\{x\}$ et $X\neq\varnothing$ sont vraies. On a $(\forall z)(z\in X\Rightarrow z=x)$ est un théorème. Puisque $X\neq\varnothing$ il existe un ensemble $\bm{T}$ tel que $\bm{T}\in X$. Donc $\bm{T}=x$ est vraie et par S6, $x\in X$ est de même. Or $x\in X$ est équivalente \'{a} $\{x\}\subset X$ (II, p. 6). D'après l'axiome d'extensionalité $\{x\}= X$ est vraie. On termine en utilisant C14 et le fait que si $\bm{A}$ et $\bm{B}$ deux relations alors $(\textrm{non }\bm{A})\Rightarrow \bm{B}$ et $(\bm{A} \textrm{ ou } \bm{B})$ sont équivalentes.

\vspace{0.5cm}
5) On a $(\exists x)(x\in X)$ et identique \'{a} $\bm{\tau}_x(x\in X)\in X$, donc $(\forall x)(x\notin X)$ est équivalente à $\bm{\tau}_x(x\in X)\notin X$. On conclut par application de C27 et S7.

\vspace{0.5cm}
6) D'abord soit $\mathcal{T}$ une théorie égalitaire dans laquelle figure le signe $\in$. On commence par quelque considérations utiles.

\begin{enumerate}[(a)]
\item D'après CS5 la relation $x\subset y$ est l'assemblage
$$\,\neg\,\neg\,\vee\,\neg\,\in\,\underbrace{\bm{\tau}\,\neg\,\vee\,\neg\,\in \,\square \,x \,\in\,\square \,y}\,x\in \,\underbrace{\bm{\tau}\,\neg\,\vee\,\neg\,\in \,\square \,x \,\in\,\square \,y}\,y $$

où dans la relation entre acolade il ya deux liens entre $\bm{\tau}$ et les deux carrés.

\item Soient $\bm{x}$ et $\bm{y}$ deux lettres distinctes entre elles.  La relation A1 est identique à $$(\forall\bm{x})(\forall\bm{y})((\bm{x}\subset\bm{y}\textrm{ et }\bm{y}\subset\bm{x})\Rightarrow \bm{x}=\bm{y}).$$  Pour le voir on distingue 6 cas: \begin{itemize}
\item l'un des $\bm{x}$ et $\bm{y}$ coincide avec $x$ et l'autre avec $y$. Soit $\bm{z}$ une lettre distincte de $x$ et de $y$. D'après CS8, A1 est identique à $(\forall x)(\forall\bm{z})((x\subset\bm{z}\textrm{ et }\bm{z}\subset x)\Rightarrow x=\bm{z}).$ Donc d'après CS8 et CS9, A1 est identique à $(\forall y)(\forall\bm{z})((y\subset\bm{z}\textrm{ et }\bm{z}\subset y)\Rightarrow y=\bm{z}).$ Enfin d'après CS8, A1 est identique à $(\forall y)(\forall x)((y\subset x\textrm{ et }x\subset y)\Rightarrow y=x).$ Les justifications dans les cas suivants se font de même.
\item $\bm{x}$ et $\bm{y}$ sont distinctes de $x$ et $y$.
\item $\bm{y}$ coincide avec $y$ et $\bm{x}$ est distincte de $x$ et de $y$.
\item $\bm{x}$ coincide avec $x$ et $\bm{y}$ est distincte de $x$ et de $y$.
\item $\bm{x}$ coincide avec $y$ et $\bm{y}$ est distincte de $x$ et de $y$.
\item $\bm{y}$ coincide avec $x$ et $\bm{x}$ est distincte de $x$ et de $y$.
\end{itemize}

\item \textbf{Relations collectivisantes dans $\mathcal{T}$.} Soient $\bm{R}$ une relation, $\bm{x}$ une lettre. Soit $\bm{y}$ une lettre distincte de $\bm{x}$ et ne figurant pas dans $\bm{R}$. D'après CS8 et CS9, la relation $(\exists \bm{y})(\forall\bm{x})((\bm{x}\in\bm{y})\Leftrightarrow \bm{R})$ est indépendente du choix de $\bm{y}$ (distincte de $\bm{x}$ et ne figurant pas dans $\bm{R}$). On la désigne par $\textrm{Coll}_{\bm{x}}\bm{R}$. D'après CS9, la relation $\textrm{Coll}_{\bm{x}}\bm{R}$ est identique à $(\forall \bm{x})((\bm{x}\in\bm{\tau}_{\bm{y}}(\bm{S}))\Leftrightarrow \bm{R})$ où $\bm{S}$ désigne la relation $(\forall \bm{x})((\bm{x}\in\bm{y})\Leftrightarrow \bm{R})$. D'après CS3 et CS9, $\bm{\tau}_{\bm{y}}(\bm{S})$ est indépendente du choix de $\bm{y}$ (distincte de $\bm{x}$ et ne figurant pas dans $\bm{R}$). On la désigne par $\{\bm{x}\mid\bm{R}\}$ (indépendente de $\bm{x}$). On dit que $\bm{R}$ est \emph{collectivisante en} $\bm{x}$ lorsque $\textrm{Coll}_{\bm{x}}\bm{R}$ est un théorème. Inversement, supposons que $(\forall \bm{x})((\bm{x}\in\bm{T})\Leftrightarrow \bm{R})$ est un théorème pour un certain ensemble $\bm{T}$ ne contenant pas $\bm{x}$.  D'après CS9 et S5, $\bm{R}$ est collectivisante en $\bm{x}$.

    Soient $\bm{x}$ une lettre, $\bm{R}$ une relation collectivisante en $\bm{x}$, et $\bm{x}'$ une lettre distincte de $\bm{x}$ et ne figurant pas dans $\bm{R}$. Soit $\bm{y}$ une lettre distincte de $\bm{x}$, de $\bm{x}'$ et ne figurant pas dans $\bm{R}$. D'après CS8, $(\forall\bm{x})((\bm{x}\in\bm{y})\Leftrightarrow \bm{R})$ est identique \'{a} $(\forall\bm{x}')((\bm{x}'\in\bm{y})\Leftrightarrow (\bm{x}'\mid\bm{x})\bm{R})$. Donc $(\exists \bm{y})(\forall\bm{x})((\bm{x}\in\bm{y})\Leftrightarrow \bm{R})$ est identique à $(\exists \bm{y})(\forall\bm{x}')((\bm{x}'\in\bm{y})\Leftrightarrow (\bm{x}'\mid\bm{x})\bm{R})$. $(\bm{x}'\mid\bm{x})\bm{R}$ est enfin collectivisante en $\bm{x}'$ et $\{\bm{x}'\mid(\bm{x}'\mid\bm{x})\bm{R}\}$ est $\{\bm{x}\mid\bm{R}\}$.

    Soit $\bm{T}$ un terme. Soit $\bm{R}$ la relation ``$\bm{z}\in\bm{T}$'' où $\bm{z}$ est une lettre ne figurant pas dans $\bm{T}$. D'après C8 et C20, $\bm{z}\in\bm{T}\Leftrightarrow\bm{z}\in\bm{T}$ est vraie. Donc par C27 et CS8, $(\forall \bm{z})(\bm{z}\in\bm{T}\Leftrightarrow\bm{z}\in\bm{T})$ est vraie (on peut supposer que $\bm{z}$ n'est pas une constante). Soit $\bm{x}$ une lettre distincte de $\bm{z}$ et ne figurant pas dans $\bm{T}$. D'après S5, $(\exists\bm{x})(\forall \bm{z})(\bm{z}\in\bm{x}\Leftrightarrow\bm{z}\in\bm{T})$ est vraie. Donc la relation $\bm{R}$ est collectivisante en $\bm{z}$ et $\{\bm{z}\mid\bm{R}\}$ est $\bm{\tau}_{\bm{x}}((\forall\bm{z})(\bm{z}\in\bm{x}\Leftrightarrow\bm{z}\in\bm{T}))$.

\item Soit $\bm{T}$ un terme. Soient $\bm{x}$ et $\bm{z}$ deux lettres distinctes entre elles et ne figurant pas dans $\bm{T}$ (*). D'après CS3, CS8 et CS9, $\bm{\tau}_x((\forall \bm{z})(z\in \bm{x}\Leftrightarrow \bm{z}\in \bm{T}))$ est indépendante du choix de $\bm{x}$ et de $\bm{z}$ (en effet, soient $\bm{x}'$ et $\bm{z}'$ deux lettres vérifiant les même conditions que $\bm{x}$ et $\bm{z}$, on le voit en distinguant deux cas: $\bm{x}$ et $\bm{x}'$ sont identiques, et le cas où elles sont distinctes). Donc, d'après C26 et C30, dire que la relation A1$'$ est vraie revient au même de dire que pour tout terme $\bm{T}$, $\bm{T}=\bm{\tau}_{\bm{x}}((\forall\bm{z})(\bm{z}\in\bm{x}\Leftrightarrow\bm{z}\in\bm{T}))$ est vraie où $\bm{x}$ et $\bm{z}$ sont deux lettres comme dans (*); ou bien de dire que tout terme est égal à ``l'ensemble'' de ses éléments.
\end{enumerate}

Maintenant soit $\mathcal{T}$ une théorie égalitaire dans laquelle figure le signe $\in$ et ayant A1$'$ comme axiome. Soient $\bm{x}$ et $\bm{y}$ deux lettres distinctes entre elles. on a $\bm{y}=\bm{\tau}_{\bm{x}}((\forall \bm{z})(\bm{z}\in\bm{x}\Leftrightarrow\bm{z}\in\bm{y}))$ où $\bm{z}$ est une lettre distincte de $\bm{x}$ et de $\bm{y}$. D'après C32, C27, CS8 et CS9, $$(\forall \bm{x})\Big((\forall \bm{z})(\bm{z}\in\bm{x}\Leftrightarrow\bm{z}\in\bm{y})\Leftrightarrow (\bm{x}\subset\bm{y}\textrm{ et }\bm{y}\subset\bm{x})\Big)$$ est vraie (on peut supposer que $\bm{x}$ n'est pas une constante).

D'après S7, $\bm{\tau}_{\bm{x}}((\forall \bm{z})(\bm{z}\in\bm{x}\Leftrightarrow\bm{z}\in\bm{y}))=\bm{\tau}_{\bm{x}}(\bm{x}\subset\bm{y}\textrm{ et }\bm{y}\subset\bm{x})$ est vraie. Supposons que $\bm{x}$ et $\bm{y}$ ne sont pas des constantes et ``$\bm{x}\subset\bm{y}\textrm{ et }\bm{y}\subset\bm{x}$'' est vraie. Montrons que celle-ci est univoque en $\bm{x}$ (i.e. fonctionelle en $\bm{x}$). Soient $\bm{x}'$ et $\bm{x}''$ deux lettres distinctes entre elles et distinctes de $\bm{x}$, de $\bm{y}$ et des constantes. Supposons que ($\bm{x}'\subset\bm{y}\textrm{ et }\bm{y}\subset\bm{x}'$) et ($\bm{x}''\subset\bm{y}\textrm{ et }\bm{y}\subset\bm{x}''$) sont vraies. D'après S2, S3, C27 et S6, $\bm{\tau}_{\bm{y}}(\bm{x}'\subset\bm{y}\textrm{ et }\bm{y}\subset\bm{x}')=\bm{\tau}_{\bm{y}}(\bm{x}''\subset\bm{y}\textrm{ et }\bm{y}\subset\bm{x}'')$ est vraie. D'après A1$'$, $\bm{x}'=\bm{x}''$ est vraie. On conclut par C45 que $\bm{x}=\bm{y}$ est vraie. Finalement on utilise  C14 et C27.

\vspace{2cm}
\begin{center}
\S 2
\end{center}

1) La première équivalence se démontre de la même fa\c{c}on que dans II, pp. 7, 8, en utilisant la généralisation de l'ex. 2), Chap.1, \S5 (I, p. 49).

Soient $\bm{A}\{\bm{x},\bm{y}\}$, $\bm{B}$ and $\bm{C}$ des assemblages où $\bm{x}$ et $\bm{y}$ sont deux lettres. D'après CS5, $(\neg\bm{A})\{\bm{B},\bm{C}\}$ est identique à $\neg\bm{A}\{\bm{B},\bm{C}\}$ (voir E I.16).

Par application de la première équivalence à la relation $\textrm{non }\bm{R}\{\bm{x},\bm{y}\}$ on a $$\Big(\textrm{non }(\forall\bm{x})(\forall\bm{y})\bm{R}\{\bm{x},\bm{y}\}\Big)\Leftrightarrow\Big(\textrm{non }(\forall\bm{z})\big((\bm{z}\textrm{ est un couple})\Rightarrow \bm{R}\{\textrm{pr}_1\,\bm{z},\textrm{pr}_2\,\bm{z}\}\big)\Big)$$ est un théorème. Donc la deuxième équivalence est vraie (on utilise C28, C31, C23, C24 et le critère de substitution précédent).

\vspace{2cm}
\begin{center}
\S 3
\end{center}

2) La première implication est déjà démontrée, cf. II, p. 12. La deuxième découle immédiatement du fait que $G^{-1}\langle A\rangle\subset\pr_2G^{-1}=\pr_1G$ pour tout ensemble $A$ (pour l'inclusion cf. II, p. 10).

\vspace{0.5cm}
3) $\Rightarrow$. Soit $(h,h')\in H$. On a $h\in \pr_1H\subset\pr_1G$. Donc $(h,g)\in G$ pour un certain ensemble $g$. Alors $(g,h)\in G^{-1}$ et par suite $(h,h')\in H\circ G^{-1}\circ G$.

$\Leftarrow$. Soit $(h,h')\in H$. Donc $(h,h')\in H\circ G^{-1}\circ G$. Alors $(h,g)\in G$ pour un certain ensemble $g$ et par suite $h\in\pr_1G$.

\vspace{0.5cm}
4) Soit $\bm{R}\{\bm{x},\bm{y}\}$ une relation, $\bm{x}$ et $\bm{y}$ étant des lettres distinctes. Soit $\bm{T}$ un ensemble, où ne figurent ni $\bm{x}$ ni $\bm{y}$, tel que $\bm{R}\Rightarrow ((\bm{x},\bm{y})\in\bm{T})$ soit vraie. Soient $\bm{G}$ le graphe de $\bm{R}$ et $\bm{z}$ une lettre distincte de $\bm{x}$, $\bm{y}$ et ne figurant pas ni dans $\bm{R}$ ni dans $\bm{T}$. On a $\bm{z}\not\in\bm{G}\Leftrightarrow (\bm{z}\not\in\bm{T}\text{ ou } \textrm{non }\bm{S})$ pour une certaine relation $\bm{S}$. D'après S2, l'implication $\bm{T}=\varnothing\Rightarrow\bm{G}=\varnothing$ est vraie.

Maintenant soit $G$ un graphe. D'après la dernière considération et II, p. 8, Prop. 3:
\begin{itemize}
  \item $\varnothing^{-1}=\varnothing$.
  \item $\varnothing\circ G=G\circ\varnothing=\varnothing$.
\end{itemize}

L'implication $G=\varnothing\Rightarrow G^{-1}\circ G=\varnothing$ découle de la première partie et du fait que si $G_1$, $G_2$, $G'$ trois graphes, l'implication $G_1=G_2\Rightarrow G'\circ G_1=G'\circ G_2$ est vraie.

L'implication $G\neq\varnothing\Rightarrow G^{-1}\circ G\neq\varnothing$ découle du fait que si $x\in G$, $(x,x)\in G^{-1}\circ G$.

\vspace{0.5cm}
7) $\Rightarrow$. Supposons que $G$ est fonctionnel et soient $X$ un ensemble et $x\in G\langle G^{-1}\langle X\rangle\rangle$. Donc $(y,x)\in G$ pour un certain $y\in G^{-1}\langle X\rangle$. On a $(y,z)\in G$ pour un certain $z\in X$. Puisque $G$ est fonctionnel, $x=z$. D'où $x\in X$.

$\Leftarrow$. Soient $x$, $y$, $z$ trois ensembles tels que $(y,x)\in G$, $(y,z)\in G$. On a $y\in G^{-1}\langle\{z\}\rangle$, donc $x\in G\langle G^{-1}\langle\{z\}\rangle\rangle\subset\{z\}$. D'où $x=z$.

\vspace{0.5cm}
8) Soient $G$ et $G'$ les graphes de $\Gamma$ et $\Gamma'$ respectivement. Montrons que $\Gamma$ et $\Gamma'$ sont des applications. Soit $(x,y)\in G$. Donc $y\in \Gamma(x)$ ($\Gamma(x)$ est la coupe de $\Gamma$ suivant $x$: $\Gamma(\{x\})$). Ce qui entraîne que $\Gamma'(y)\subset\Gamma'(\Gamma(x))=\{x\}$ (cf. II, p. 10, Prop. 2). Puisque $\Gamma(\Gamma'(y))=\{y\}$, on a $\Gamma'(y)\neq\varnothing$. Donc $\Gamma'(y)=\{x\}$ (cf. II, p. 49, \S1, exerc. 4). Si en plus on a $(x,y')\in G$, on aurra $y'\in \Gamma(\Gamma'(y))=\{y\}$. Alors $y=y'$ et $G$ est fonctionnel. Soit $x\in A$. Par hypothèse $\Gamma(x)\neq\varnothing$. Donc $\pr_1G=A$. Finalement, $\Gamma$ est une application et par symmétrie $\Gamma'$ est de même; et par hypothèse, $\Gamma'(\Gamma(x))=x$ pour tout $x\in A$ et $\Gamma(\Gamma'(y))=y$ pour tout $y\in B$ (cf. II, p. 18, Corollaire).

\vspace{0.5cm}
9) D'après, Prop. 8, II, p. 18, l'application $g$ est injective et surjective donc bijective. Alors $f$, $g$ et $h$ sont des bijections.

\vspace{0.5cm}
11) L'erreur est ``$f(\varnothing)=\varnothing$''. On a vu (II, p. 10) que ``$f\langle\varnothing\rangle=\varnothing$''. L'erreur provient de la confusion de $f(\varnothing)$ avec $f\langle\varnothing\rangle$ (cf. II, p. 14, Remarque).

\vspace{2cm}
\begin{center}
\S 4
\end{center}

1) Si $X$ un ensemble, on a $G\langle X\rangle=G\langle X\cap\pr_1G\rangle$. Donc l'implication $a)\Rightarrow b)$ est un cas particulier de Compléments au Chapitre 2, Prop. 26 (3). L'implication $b)\Rightarrow c)$ est évidente puisque $G^{-1}\langle\varnothing\rangle=\varnothing$ (cf. II, p. 10). L'implication $c)\Rightarrow a)$. Soient $(x,y)\in G$ et $(x,y')\in G$. Donc $x\in G^{-1}(\{y\})\cap G^{-1}(\{y'\})$. Alors $\{y\}\cap\{y'\}\neq\varnothing$, c'est-à-dire $y=y'$.

\vspace{0.5cm}
8) a) D'abord on a $E=\bigcup\limits_{i\in I}X_i$ et on doit supposer que chacun des $Y_k$ est non vide. Pour montrer la propriété énoncée dans a), soient $k\in K$ et $x\in Y_k$. Donc il existe $i\in  I$ tel que $x\in X_i$. Par hypothèse, $X_i\subset Y_{k'}$ pour un certain $k'\in K$. On a $k'=k$ car sinon on aurait $x\in Y_k\cap Y_{k'}=\varnothing$. Donc $X_i\subset Y_{k'}=Y_k$.

b) Il suffit de prendre $E=\{1,2\}$, $\mathfrak{R}=\{E\}$ et $\mathfrak{S}=\{\{1\},E\}$.

c) Il suffit de prendre $E=\{1,2,3,4\}$, $\mathfrak{R}=\{\{1\},\{3\},\{2,4\}\}$ et $\mathfrak{S}=\{\{1,2\},\{3,4\}\}$.

\vspace{2cm}
\begin{center}
\S 5
\end{center}

1) C'est évident puisque $X\in \mathfrak{P}(X)$.

\vspace{0.5cm}
2) On propose de démontrer le résultat plus générale suivant: Soient $(E,\leq)$ un ensemble ordonné et $f:E\to E$ une application croissante. Supposons que la borne inférieure de l'ensemble $A$ des $z\in E$ tels que $f(z)\leq z$ existe, on la note par $v$, et supposons que la borne supérieure de l'ensemble $B$ des $z\in E$ tels que $z\leq f(z)$ existe, on la note par $w$. On a $f(v)=v$ et $f(w)=w$, et que pour tout $z\in E$ tel que $f(z)=z$, on a $v\leq z\leq w$.

On a $f(v)\leq f(z)\leq z$ pour tout $z\in A$ (car $f$ est croissante). Donc $f(v)$ est un minorant de $A$ et par suite $f(v)\leq v$, c'est-à-dire, $f(z)\in A$. D'où $f(v)$ est le plus petit élément de $A$ et finalement $f(v)=v$ (cf. III, p. 10). La seconde égalité découle de la première en considérant l'ordre opposé sur $E$. La dernière assertion est triviale.

\textbf{Remarque.} Après la rédaction de cette solution, j'ai trouvé que la généralisation que j'ai proposée généralise aussi une partie du théorème de Knaster-Tarski \cite[Theorem 1 (Lattice Theoretical Fixpoint Theorem)]{Tarski:1955}. Dans cet article A. tarski considère le cas où $E$ est un ensemble réticulé achevé (cf. \cite{Birkhoff: 1984} ou \cite[l'exerc. 11 de III, p. 71]{Bourbaki:1970}), et montre que l'ensemble des points fixes $P=\{z\in E\mid f(z)=z\}$ est non vide en montrant que $u$ et $v$ définis comme ci-dessus appartiennent  à $P$, et $P$ est aussi un ensemble réticulé achevé. Pour démonter cette dernière proposition, soit $Y$ une partie non vide de $P$. Puisque $f$ est croissante, alors $\sup_E Y\leq f(\sup_E Y)$. Donc $[\sup_E Y,1]$ est stable par $f$. En outre, si $S$ une partie de $E$ est réticulée achevée munie de l'ordre induit et stable par $f$, alors $P\cap S$ et non vide et admet un plus petit élément et un plus grand élément (on considère l'application induite $f_S:S\to S$). En appliquant cela à l'intervalle $S=[\sup_E Y,1]$, on trouve que le plus petit élément de $P\cap S$ est bien la borne supérieure de $Y$ dans $P$.

\vspace{0.5cm}
3) C'est évident d'après II, p. 32, Définition 1.

\vspace{0.5cm}
4) Soit $f:A\to \mathfrak{P}(B)$ une application. Il est clair qu'il existe $G\in\mathfrak{P}(A\times B)$ unique tel que $f(x)=G\langle\{x\}\rangle$, pour tout $x\in A$. Pour l'existence on prend $G=\bigcup\limits_{x\in A}(\{x\}\times f(x))$. Donc on a une bijection de $\mathfrak{P}(A\times B)$ sur $\mathcal{F}(A;\mathfrak{P}(B))$ (l'ensemble des applications de $A$ dans $\mathfrak{P}(B)$) qui s'identifie à $(\mathfrak{P}(B))^A$.

\vspace{0.5cm}
5) Soit $(X_i)_{1\leq i\leq n}$ une famille finie d'ensembles. D'après la prop. 8 de II, p. 35, $$\bigcap\limits_{H\in \mathfrak{F}_k}P_H=\bigcap\limits_{H\in \mathfrak{F}_k}(\bigcup\limits_{i\in H}X_i)=\bigcup\limits_{f\in I}(\bigcap\limits_{H\in \mathfrak{F}_k}X_{f(H)})$$ où $I=\prod\limits_{H\in \mathfrak{F}_k}H$. On a $\bigcup\limits_{H\in \mathfrak{F}_k}H=[1,n]$. Les éléments de $I$ sont les graphes fonctionnels $f=(f(H))_{H\in\mathfrak{F}_k}$ contenus dans $\mathfrak{F}_k\times[1,n]$ tels que $f(H)\in H$.

Si $k\leq \frac12(n+1)$. Soit $x\in\bigcup\limits_{f\in I}(\bigcap\limits_{H\in\mathfrak{F}_k}X_{f(H)})$. Donc il existe $f\in I$ tel que $x\in\bigcap\limits_{H\in\mathfrak{F}_k}X_{f(H)}$. On a $\pr_2f$ a au moins $k$ éléments. En effet, si $\pr_2f$ avait $l$ éléments avec $l<k$, l'ensemble $[1,n]-\pr_2f$ aurait $n-l$ éléments avec $k\leq n-l$, car $k-1\leq n-k<n-l$. Donc il existerait $H\in\mathfrak{F}_k$ tel que $H\subset [1,n]-\pr_2f$. Alors $f(H)\in H\cap\pr_2f$, ce qui serait une contradiction. D'où $x\in\bigcup\limits_{H\in\mathfrak{F}_k}Q_H=\bigcup\limits_{H\in\mathfrak{F}_k}(\bigcap\limits_{i\in H}X_i)$.

Si $k\geq \frac12(n+1)$. Soit $x\in\bigcup\limits_{H\in\mathfrak{F}_k}Q_H$. Donc il existe $H'\in\mathfrak{F}_k$ tel que $x\in Q_{H'}$. Pour tout $H\in\mathfrak{F}_k$, $H\cap H'\neq\varnothing$, car $2k\geq n+1$. D'après l'axiome de choix on a $\varnothing\neq\prod\limits_{H\in\mathfrak{F}_k}(H\cap H')\subset I$. Soit $f\in\prod\limits_{H\in\mathfrak{F}_k}(H\cap H')$. Donc $x\in Q_{H'}=\bigcap\limits_{i\in H'}X_i\subset\bigcap\limits_{H\in \mathfrak{F}_k}X_{f(H)}$. D'où $x\in\bigcup\limits_{f\in I}(\bigcap\limits_{H\in \mathfrak{F}_k}X_{f(H)})$.

\vspace{2cm}
\begin{center}
\S 6
\end{center}

1) Soient $G$ un graphe et $E$ un ensemble. D'abord si $\pr_1G=\pr_2G=E$ (resp. $\pr_1G\subset E$ et $\pr_2G\subset E$) et $\Delta_E\subset G$, alors $G\subset G\circ G$ et les relations suivantes sont équivalentes:
\begin{enumerate}[i)]
  \item $G\circ G=G$ et $G^{-1}=G$;
  \item $G\circ G^{-1}\circ G=G$;
  \item $G\circ G^{-1}\circ G\subset G$.
\end{enumerate}

En effet, $G\subset G\circ G$ est vraie car si $(x,y)\in G$, alors $(y,y)\in\Delta_E\subset G$ et par suite $(x,y)\in G\circ G$. Les relations ii) et iii) sont équivalentes d'après II, p. 50, \S3, exerc. 3.

$i)\Rightarrow ii)$. On a $G\circ G^{-1}\circ G=G\circ G\circ G=G\circ G=G$.

$iii)\Rightarrow i)$. $G\circ G\subset G$: Soit $(x,y)\in G\circ G$. Donc il existe un ensemble $z$ tel que $(x,z)\in G$ et $(z,y)\in G$. Alors $(z,z)\in\Delta_E\subset G$ et par suite $(x,y)\in G\circ G^{-1}\circ G\subset G$.

$G^{-1}=G$: Soit $(x,y)\in G$. Donc $(y,y)\in\Delta_E\subset G$ et $(y,x)\in G^{-1}$ et $(x,x)\in\Delta_E\subset G$. Par suite $(y,x)\in G\circ G^{-1}\circ G\subset G$.

Donc, d'après II, p. 41, Prop. 1, pour que $G$ soit le graphe d'une équivalence dans $E$, il faut et il suffit qu'une des conditions suivantes soit vraie:
\begin{enumerate}[(a)]
  \item $\pr_1G=\pr_2G=E$ (ou bien $\pr_1G\subset E$ et $\pr_2G\subset E$), $\Delta_E\subset G$, $G\circ G=G$ (ou bien $G\circ G\subset G$) et $G^{-1}=G$;
  \item $\pr_1G=\pr_2G=E$ (ou bien $\pr_1G\subset E$ et $\pr_2G\subset E$), $\Delta_E\subset G$, $G\circ G^{-1}\circ G=G$
  (ou bien $G\circ G^{-1}\circ G\subset G$).
  \end{enumerate}
Noter que pour que $G$ soit le graphe d'une équivalence dans $E$, il faut et il suffit que $G$ soit le graphe d'une relation équivalence dans $E$.

\vspace{0.5cm}
2) Soit $G$ un graphe. On a $\pr_1(G^{-1}\circ G)=\pr_2(G^{-1}\circ G)=\pr_1G$ et $(G^{-1}\circ G)^{-1}=G^{-1}\circ G$. Maintenant supposons que $G\circ G^{-1}\circ G=G$. Donc $G^{-1}\circ G$ est le graphe d'une équivalence dans $\pr_1G$ car $G^{-1}\circ G\circ G^{-1}\circ G=G^{-1}\circ G$ (cf. la condition (a) énoncée dans la solution de l'exercice 1). De même on montre l'autre relation.

\vspace{0.5cm}
3) Désignons cette application par $f$. Il suffit de décomposer canoniquement $f$ (II, p. 44) en notant que $f\langle\mathfrak{P}(E)\rangle=\mathfrak{P}(A)$.

\vspace{0.5cm}
4) Soit $G$ le graphe d'une équivalence dans un ensemble $E$. Supposons que $A$ est un graphe tel que $A\subset G$ et $\pr_1A=E$. On a $G\circ A\subset G$ car $A\subset G$ et $G\circ G=G$. Inversement, soit $(x,y)\in G$. Donc $x\in \pr_1G=E=\pr_1A$ et par suite il existe ensemble $z$ tel que $(x,z)\in A$. Donc $(z,x)\in G$ et $(z,y)\in G\circ G=G$. Enfin $(x,y)\in G\circ A$. Soit $B$ un graphe quelconque. $(G\cap B)\circ A$ est contenu dans $G\circ A=G$ et dans $B\circ A$ donc dans $G\cap (B\circ A)$. Inversement, soit $(x,y)\in G\cap (B\circ A)$. Donc $(x,y)\in G$ et il existe un ensemble $z$ tel que $(x,z)\in A$ et $(z,y)\in B$. Puisque $G\circ G=G$ et $G^{-1}=G$, alors $(z,y)\in G$. Enfin $(x,y)\in (G\cap B)\circ A$. De même on montre le cas respective.

\vspace{0.5cm}
5) Soit $(G_i)_{i\in I}$ une famille de graphes d'équivalence dans un ensemble $E$. Il est facile de vérifier que $\pr_1(\bigcap\limits_{i\in I}G_i)\subset E$, $\pr_2(\bigcap\limits_{i\in I}G_i)\subset E$, $\Delta_E\subset \bigcap\limits_{i\in I}G_i$, $(\bigcap\limits_{i\in I}G_i)\circ (\bigcap\limits_{i\in I}G_i)\subset\bigcap\limits_{i\in I}G_i$ et $(\bigcap\limits_{i\in I}G_i)^{-1}=\bigcap\limits_{i\in I}G_i$. D'où la première partie (cf. la condition (a) énoncée dans la solution de l'exercice 1). La réunion de la famille $(G_i)_{i\in I}$ vérifie tous les conditions vérifiées pour l'intersection sauf $(\bigcup\limits_{i\in I}G_i)\circ (\bigcup\limits_{i\in I}G_i)\subset\bigcup\limits_{i\in I}G_i$. Pour le contre exemple, soit $E=\{1,2,3\}$, $R$ et $S$ les relations d'équivalence associées aux partitions $\{\{1,2\},\{3\}\}$ et $\{\{1\},\{2,3\}\}$ de $E$ respectivement, et $G$ et $H$ leur graphes respectives. On a $(1,2)\in G$ et $(2,3)\in H$. Donc $(1,3)\in (G\cup H)\circ(G\cup H)$. Or $G=\{(1,2),(2,1),(3,3)\}$ et $H=\{(1,1),(2,3),(3,2)\}$. Par suite $(1,3)\not\in G\cup H$.

\vspace{0.5cm}
6) Soient $G$ et $H$ les graphes de deux équivalences dans un ensemble $E$. On a $\pr_1(G\circ H)=H^{-1}\langle\pr_1G\rangle=H^{-1}\langle E\rangle=E$ et $\pr_2(G\circ H)=G\langle\pr_2H\rangle=G\langle E\rangle=E$ (cf. II, p. 12). Soit $(x,y)\in G\circ H$. Donc il existe un ensemble $z$ tel que $(x,z)\in H$ et $(z,y)\in G$. Alors $(y,z)\in G$ et $(z,x)\in H$ et par suite $(y,x)\in H\circ G$. D'où $(G\circ H)^{-1}=H\circ G$. Ce qui achève la démonstration de la première partie (cf. la condition (a) énoncée dans la solution de l'exercice 1). Soit $(G_i)_{i\in I}$ la famille de graphes d'équivalence dans $E$ contenant $G$ et $H$. On a $G\circ H\subset(\bigcap\limits_{i\in I}G_i)\circ (\bigcap\limits_{i\in I}G_i)\subset\bigcap\limits_{i\in I}G_i$. D'autre part, $G\subset G\circ H$ et $H\subset G\circ H$, donc $\bigcap\limits_{i\in I}G_i\subset G\circ H$.

\vspace{0.5cm}
7) Soient $p$ l'application canonique de $F$ sur $F/R$ et $j$ l'injection canonique de $A$ dans $F$. L'application $f$ est compatible avec les relations $S$ et $R$ (cf. II, p. 44). Soient $h$ l'application déduite de $f$ par passage aux quotients suivant $S$ et $R$ et $l$ l'application déduite de $j$ par passage aux quotients suivant $S$ et $R$. Il facile de voir que $h$ est injective. Donc $l$ est aussi injective. Alors on obtient le diagramme suivant:
$$\xymatrix{E \ar[r]^f \ar[d] & F\ar[d]^p & A\ar[l]^j\ar[d]
\\ E/S\ar[r]^h \ar[dr]^{h'} & F/R & A/R_A\ar[l]^{l}\ar[dl]^{l'}
\\ & p\langle A\rangle\ar[u],}$$
où $h'$ et $l'$ sont les bijections canoniques. D'où $l'^{-1}\circ h'$ est une bijection canonique de $E/S$ sur $A/R_A$.

\vspace{0.5cm}
9) a) On a évidemment $(R\{x,y\}\text{ et }R\{y,x\})\Rightarrow (R\{y,x\}\text{ et }R\{x,y\})$. Donc la relation $R\{x,y\}\text{ et }R\{y,x\}$ est symétrique. Pour que cette relation soit réflexive dans un ensemble $E$ il faut et il suffit que $R\{x,y\}$ le soit par rapport aux lettres $x$ et $y$.

b) Soient $G$ le graphe de $R\{x,y\}$ et $z$ une lettre ne figurant pas dans $R\{x,y\}$. En considérant la suite $x'_i=x_{n-i}$ pour tout $i$, on voit que $S\{x,y\}$ est symétrique. Supposons que $S\{x,y\}$ et $S\{y,z\}$. En juxtaposant les deux suites conséquantes on voit que $S\{x,z\}$ est vraie aussi. On a évidemment $S\{x,y\}$ est réflexive dans $E$. Donc $S\{x,y\}$ est une relation d'équivalence dans $E$. Soit $G'$ sont graphe. En outre $G'=\bigcup\limits_{n\in \N-\{0\}}G^n$ où $G^n=G\circ G^{n-1}$ avec $G^1$=G. On a $(G^n)_{n\in \N-\{0\}}$ est une famille de parties de $G'$ car $G'$ est le graphe d'une équivalence dans $E$ (cf. la condition (a) énoncée dans la solution de l'exercice 1). Si $H$ est le graphe d'une équivalence dans $E$ contenant $G$, il contient aussi $G^n$ pour tout $n\in \N-\{0\}$ et par suite $G'$. D'où la dernière assertion.

c) Deux éléments de $E$ sont équivalents modulo $S$ s'il existe un chemin (peut être vide) de l'un vers l'autre. Soit $x\in E$. La composante connexe de $x$ suivant $R$ est l'ensemble des éléments de $E$ tels qu'il existe un chemin de $x$ vers eux, c'est donc l'ensemble $A\in \mathfrak{F}$ tel que $x\in A$.

\chapter*{Solutions de quelques exercices du chapitre III (Ensembles ordonnés, cardinaux, nombres entiers)}

\begin{center}
\S 1
\end{center}

1) D'abord $R\{x,y\}$ vérifie les deux premières conditions de III, p. 1, car $x<y \Rightarrow x\leq y$. Soient $X$ et $Y$ deux éléments distincts et comparables de $E$. Supposons que $X\leq Y$. Donc $R\{X,Y\}$ est vraie. On a les équivalences:
$(\text{non }R\{x,x\})\Leftrightarrow (x\not\in E\text{ ou } x\nless x)$; $x\nless x \Leftrightarrow (x\nleqslant x\text{ ou } x=x)$. Donc $x<x$ est fausse et par suite  $R\{x,x\}$ est aussi fausse. Puisque la relation $\text{non }(A\Rightarrow B)$ est équivalente à $(A\text{ et }\text{non }B)$, la relation $R\{X,Y\}\Rightarrow R\{X,X\}$ est fausse. Or $R\{x,y\}\Rightarrow R\{x,x\}$ implique $R\{X,Y\}\Rightarrow R\{X,X\}$. D'où $R\{x,y\}\Rightarrow R\{x,x\}$ est fausse et $R\{x,y\}$ ne vérifie pas la troisième relation de III, p. 1. Même chose si $Y\leq X$. On conclut à l'aide de C18 (I, p. 28).

\vspace{0.5cm}
7) Soient $E$ et $F$ deux ensembles ordonnés tels que $F$ a au moins deux éléments. D'abord si $x$ et $y$ sont deux éléments comparables de $E$ et $f:E\to F$ une application à la fois croissante et décroissante, alors $f(x)=f(y)$ (cf. C18 (I, p. 28)). Donc, si $E$ est connexe, toute application de $E$ dans $F$  qui est à la fois croissante et décroissante est constante. Supposons que $E$ n'est pas connexe. Donc il existe plus de deux composantes connexes pour la relation ``$x$ et $y$ sont comparables''. Soient $x'$ et $y'$ deux éléments distincts de $F$ et soit $C$ une composante connexe. Soit $f$ l'application qui associe $x'$ aux éléments de $C$ et $y'$ aux éléments des autres composantes connexes. Clairement $f$ est croissante et décroissante mais pas constante.

\vspace{0.5cm}
8) D'abord, si $x\in A$, on a $f(g(f(x)))=f(x)$, donc $f(x)\in B$. De même si $y\in B$, $g(y)\in A$. Donc on peut considérer les applications $f':A\to B$ et $g':B\to A$ induites de $f$ et $g$. On a $g'\circ f'=\text{Id}_{A'}$ et $f'\circ g'=\text{Id}_{B'}$, donc $f'$ est une bijection et $g'$ sa bijection réciproque. En plus, puisque $f$ et $g$ sont croissantes, $f'$ et $g'$ sont aussi croissantes. D'où le résultat (cf. III, p. 7).

\vspace{0.5cm}
9) Soient $E$ un ensemble réticulé et $(x_{ij})$ une famille double finie. D'abord, d'après la prop. 7, III, p. 11, toute partie finie de $E$ admet une borne supérieure et une borne inférieure dans $E$. On a $\inf\limits_i x_{ij}\leq x_{ij}$ pour tout couple $(i,j)$. D'après la prop. 6, III, p. 11, $\sup\limits_j(\inf\limits_i x_{ij})\leq \sup\limits_jx_{ij}$ pour tout $i$. Donc $\sup\limits_j(\inf\limits_i x_{ij})\leq \inf\limits_i(\sup\limits_jx_{ij})$ par la définition de la borne inférieure.

\vspace{0.5cm}
10) Soient $E$ et $F$ deux ensembles réticulés et $f:E\to F$ une application. Supposons que $f$ est croissante et soient $x,y\in E$. On a $f(\inf(x,y))\leq f(x)$ et $f(\inf(x,y))\leq f(y)$. Donc $f(\inf(x,y))\leq \inf(f(x),f(y))$.

Maintenant, supposons que, pour tout $x,y\in E$, $f(\inf(x,y))\leq \inf(f(x),f(y))$ et soient $x',y'\in E$ tels que $x'\leq y'$. Donc $\inf(x',y')=x'$ et $f(x')\leq \inf(f(x'),f(y'))$. Enfin $f(x')\leq f(y')$.

Soit $f$ l'application croissante de $\N\times\N$ (réticulé) dans $\N$ (réticulé) telle que $f((m,n))=m+n$ pour tout $(m,n)\in\N\times\N$. Soient $x=(0,1)$ et $y=(1,0)$. On a $\inf(x,y)=(0,0)$ et $f(x)=f(y)=1$. Donc $f(\inf(x,y))=0<1=\inf(f(x),f(y))$.

\vspace{2cm}
\begin{center}
\S 2
\end{center}

3) Soit $E$ un ensemble ordonné. Comme c'est indiqué on prend $B$ la réunion des parties de $E$ n'ayant pas de plus petit élément (peut être vide). Soit $A$ le complémentaire de $B$ par rapport à $E$. Donc $A$ est bien ordonné car s'il existait une partie non vide $A'$ de $A$ n'admettant pas un plus petit élément, on aurait $A'\subset A\cap B=\varnothing$. Par définition $B$ n'admet pas de plus petit élément. Pour l'exemple il suffit de prendre $E=\mathbb{Z}$ muni de l'ordre usuel, $A$ une partie finie de $\mathbb{Z}$ et $B$ son complémentaire.

\vspace{0.5cm}
4) Plus généralement et avec la même méthode on peut même monter que dans tout ensemble ordonné $E$, il existe une partie $F$ bien ordonné et cofinale (cf. III, p. 9) à $E$. Soient $(E,\leq)$ un ensemble ordonné et $\mathfrak{F}$ l'ensemble des parties bien ordonnées de $E$ muni de la relation d'ordre ``$X\subset Y$ et aucun élément de $Y-X$ n'est majoré par un élément de $X$'' entre $X$ et $Y$, que nous notons par $X\subset' Y$. En effet, cette relation est transitive et antisymétrique (2ème relation de III, p. 1). Donc cette relation est une relation d'ordre entre éléments de $\mathfrak{F}$ (cf. III, p. 2, Exemple 2, et l'annexe \ref{app3}). Soient $X,Y\in\mathfrak{F}$ tels que $X\subset' Y$. On a $X$ est un segment de $Y$ pour l'ordre induit par celui de $E$. En effet, soient $x\in X$ et $y\in Y$ tels que $y\leq x$. Si $y\not\in X$, on aurait $y\in Y-X$, ce qui contredirait $y\leq x$. Montrons que $\mathfrak{F}$ est inductif. Soit $\mathfrak{S}$ une partie totalement ordonné de $\mathfrak{F}$. D'après Prop. 3, III, p. 16, $\bigcup\limits_{X\in \mathfrak{S}}X$ est un ensemble bien ordonné pour l'ordre induit. Soient $Y\in\mathfrak{S}$ et $x\in (\bigcup\limits_{X\in\mathfrak{S}}X)-Y$. Donc $x\in X$ pour un certain $X\in \mathfrak{S}$. D'après C25 (I, p. 31), $Y\subset' X$. Par suite $x$ n'est pas majoré par un élément de $Y$. D'où $Y\subset'(\bigcup\limits_{X\in\mathfrak{S}}X)$ et $\bigcup\limits_{X\in\mathfrak{S}}X$ est un majorant de $\mathfrak{S}$ dans $\mathfrak{F}$. D'après le théorème de Zorn (III, p. 20, Théorème 2) $\mathfrak{F}$ possède un élément maximal $M$. On a $M$ est cofinale dans $E$. En effet, soit $x\in E-M$ et supposons que pour tout $y\in M$, $y<x$. Donc $x$ est un majorant de $M$. Par suite $M\cup\{x\}$ est bien ordonnée (cf. la remarque suivante) et contenant strictement $M$. D'où $x$ est majoré par un élément de $M$, ce qui donne une contradiction.

\textbf{Remarque.} Soient $(E,\leq)$ un ensemble ordonné, $A$ une partie de $E$ bien ordonnée et $x$ un majorant de $A$ n'appartenant pas à $A$. Il existe un isomorphisme de l'ensemble ordonné obtenu en adjoignant à $A$ un plus grand élément $x$ (cf. III, p. 9, Proposition 3) sur la partie $A\cup\{x\}$ de $E$ munie de l'ordre induit. D'après III, p. 16, Exemple 5, $A\cup\{x\}$ est bien ordonnée.

\vspace{0.5cm}
11) Notons $E$ le produit lexicographique des $E_i$. D'abord le cas où $I$ est l'ensemble vide est trivial ($E=\{\varnothing\}$, cf. II, p. 32). Donc assumons $I\neq\varnothing$.

Supposons que $E$ est bien ordonné. Soient $i\in I$ et $A$ une partie non vide de $E_i$. Soit $x_j\in E_j$ pour tout $j\in I-\{i\}$ (à l'aide de l'axiome de choix). Considérons l'ensemble $A'$ des éléments $x\in E$ tels que $\pr_ix\in A$ et $\pr_jx=x_j$. Donc $A'$ admet un plus petit élément $x_0$. Notons $a_0=\pr_ix_0$. On a $a_0<a$ pour tout $a\in A-\{a_0\}$. D'où $a_0$ est le plus petit élément de $A$ et $E_i$ est bien ordonné. Maintenant, supposons que $I$ est infini. Soient $x_i,y_i\in E_i$ tels que $x_i<y_i$, pour tout $i\in I$ (cf. les hypothèses). Soient $i_0$ le plus petit élément de $I$; $i_1$ le plus petit élément de $I-\{i_0\}$; $i_2$ le plus petit élément de $I-\{i_0,i_1\}$, et ainsi de suite. Donc $(i_n)_{n\in\N}$ est une suite strictement croissante d'éléments de $I$: $i_0<i_1<i_2\dots$. Prenons $y_0=(y_i)_{i\in I}\in E$; $y_1\in E$ tel que $\pr_iy_1=y_i$ pour tout $i\in I-\{i_0\}$ et $\pr_{i_0}y_1=x_{i_0}$; $y_2\in E$ tel que $\pr_iy_2=y_i$ pour tout $i\in I-\{i_0,i_1\}$, $\pr_{i_0}y_2=x_{i_0}$ et $\pr_{i_1}y_2=x_{i_1}$, et ainsi de suite. La suite $(y_n)_{n\in \N}$ est une suite strictement décroissante d'éléments de $E$: $\dots y_2<y_1<y_0$. Or l'image de cette suite admet un plus petit élément, ce qui entraîne une contradiction. D'où $I$ est fini.

Finalement, supposons que les $E_i$ sont bien ordonnés et $I$ est fini et soit $A$ une partie non vide de $E$. En considérant la suite strictement croissante $(i_n)_{n\in \N}$ d'éléments de $I$ définie ci-dessus, il existe $k\in\N$ tel que $I=\{i_0,\dots,i_k\}$. Soient $a_{i_0}$ le plus petit élément de $\pr_{i_0}A\subset E_{i_0}$; $a_{i_1}$ le plus petit élément de $\pr_{i_1}A_{i_1}$ où $A_{i_1}$ est l'ensemble des éléments $x$ de $A$ tels que $\pr_{i_0}x=a_{i_0}$; $a_{i_2}$ le plus petit élément de $\pr_{i_2}A_{i_2}$ où $A_{i_2}$ est l'ensemble des éléments $x$ de $A$ tels que $\pr_{i_0}x=a_{i_0}$ et $\pr_{i_1}x=a_{i_1}$, et ainsi de suite. On a $(a_i)_{i\in I}$ est le plus petit élément de $A$. D'où $E$ est bien ordonné.

\vspace{0.5cm}
12) Soit $(E_i)_{i\in I}$ une famille d'ensembles ordonnés, dont l'ensemble d'indices $I$ soit totalement ordonné. On muni l'ensemble $E=\prod\limits_{i\in I}E_i$ de la relation d'ordre: ``$x=(x_i)_{i\in I}\leq y=(y_i)_{i\in I}$'' si $$x=y \text{ ou } (x\neq y \text{ et } \{i\in I\mid x_i\neq y_i\} \text{ est bien ordonné} \text{ et } x_{i_0}<y_{i_0} \text{ où } i_0=\inf\{i\in I\mid x_i\neq y_i\}).$$
Le fait que c'est bien une relation d'ordre se justifie de la même façon que pour montrer que l'ordre lexicographique est bien un ordre (cf. Annexe \ref{app3}) en utilisant l'observation suivante*: Soient $x$, $y$ et $z$ sont des éléments de $E$ tels que $x$ et $y$ sont distincts et comparables pour $\leq$ et, $y$ et $z$ sont distincts et comparables pour $\leq$. On a $\{i\in I\mid x_i= y_i\}\cap\{i\in I\mid y_i= z_i\}\subset\{i\in I\mid x_i= z_i\}$. Donc $\{i\in I\mid x_i\neq z_i\}\subset\{i\in I\mid x_i\neq y_i\}\cup\{i\in I\mid y_i\neq z_i\}$ (cf. II, p. 6 et II, p. 27). D'après la remarque suivante, $\{i\in I\mid x_i\neq z_i\}$ est bien ordonné.

\textbf{Remarque.} Soit $E$ un ensemble totalement ordonné et $A$ et $B$ deux parties bien ordonnées de $E$. Donc $A\cup B$ est bien ordonnée.  En effet, soit $X$ une partie non vide de $A\cup B$. Si $X$ est contenu dans $A$ ou dans $B$, elle admet un plus petit élément. Dans le cas contraire, on a $X=(X\cap A)\cup (X\cap B)$, $\varnothing\neq X\cap A\subset A$ et $\varnothing\neq X\cap B\subset B$. Soient $a$ le plus petit élément de $X\cap A$ et et $b$ le plus petit élément de $X\cap B$. Donc $\inf(a,b)$ est le plus petit élément de $X$.

La seconde assertion est immédiate en utilisant l'observation (*). Maintenant, supposons que chaque $E_i$ a au moins deux éléments et $E$ est totalement ordonné. Donc, soient $x_i$ et $y_i$ deux éléments distincts de $E_i$ pour tout $i\in I$. On a $(x_i)_{i\in I}$ et $(y_i)_{i\in I}$ deux éléments distincts de $E$. Donc $I=\{i\in I\mid x_i\neq y_i\}$ est bien ordonné. D'où l'ordre définit sur $E$ est bien l'ordre lexicographique des $E_i$. Soient $j\in I$ et $x_j$ et $y_j$ deux éléments distints de $E_j$. Choisissons $x_i\in E_i$ pour tout $i\in I-\{j\}$. Les deux éléments $x=(x_i)_{i\in I}$ et $y=(y_i)_{i\in I}$ de $E$ tels que $y_i=x_i$ pour tout $i\in I-\{j\}$, sont donc distincts et comparables. Par suite $x_j<y_j\text{ ou }y_j<x_j$. D'où $x_j$ et $y_j$ sont comparables. L'autre implication est un résultat de cours (cf. III, p. 23).

\textbf{Remarque.} Nous nous sommes pas servi de l'exerc. 3 de III, p. 75, comme c'est indiqué dans l'énoncé de cet exercice.

\vspace{2cm}
\begin{center}
\S 3
\end{center}

1) D'abord la condition ``$f$ est injective'' n'est pas nécessaire. Il suffit de montrer qu'il existe une partie $A$ de $E$ telle que $E-A=g(F-f(A))$. Or cela est évident en appliquant l'exercice 2, II, p. 51, \S 5, en considérant l'application croissante de $\mathfrak{P}(E)$ dans lui-même qui a $X$ associe $E-g(F-f(X))$.

Comme conséquence immédiate de cet exercice et du lemme \ref{II, p. 29, prop. 8'} de l'annexe B, le théorème de Cantor-Bernstein: S'il existe une injection de $E$ dans $F$ et une injection de $F$ dans $E$, alors il existe une bijection de $E$ sur $F$ (cf. III, p. 25, le cor. 2 du th. 1).

\vspace{0.5cm}
2) Soient $E$ et $F$ deux ensembles distincts. Si l'un des $E$ et $F$ est l'ensemble vide, par exemple $E=\varnothing$ et $F\neq\varnothing$, alors $F^E=\{\varnothing\}$ et $E^F=\varnothing$. Maintenant, supposons que $E$ et $F$ sont non vides et $F^E=E^F$. Soit $y\in F$. Le graphe fonctionnel formé des éléments $(x,y)$ où $x\in E$ appartient à $F^E$ donc à $E^F$. D'où $E\subset F$, et par symétrie on a $E=F$. Ce qui est absurde.

On a $2^4=16=4^2$ comme cardinaux. Donc l'un au moins des ensembles $2^4$ et $4^2$ n'est pas un cardinal (car sinon on aurait $2=4$).

Plus généralement, soient $\mathfrak{a}$ et $\mathfrak{b}$ deux cardinaux tels que $\mathfrak{a}\neq0$ ou $\mathfrak{b}=0$. Le cardinal $\mathfrak{a}^\mathfrak{b}$ est distinct de l'ensemble des applications de $\mathfrak{b}$ dans $\mathfrak{a}$, $\mathfrak{a}^\mathfrak{b}$, car $\varnothing$ appartient au premier mais pas au second. Donc l'ensemble des applications de $\mathfrak{b}$ dans $\mathfrak{a}$ n'est pas un cardinal.

\vspace{0.5cm}
3) Soient $(\mathfrak{a}_i)_{i\in I}$ et $(\mathfrak{b}_i)_{i\in I}$ deux familles de cardinaux, telles que $\mathfrak{b}_i\geq2$ pour tout $i\in I$.

a) D'après la prop. 14 de III, p. 30, il suffit de montrer que $\sum\limits_{i\in I}\mathfrak{b}_i\leq\underset{i\in I}{\p}\mathfrak{b}_{i}$. Pour cela on distingue deux cas. D'abord si $I$ et les $\mathfrak{b}_i$ sont finis. Supposons que $I=\{1,\dots, n\}$. On raisonne par récurrence sur $n$. Le cas où $n=1$ est trivial. On a besoin de montrer aussi le cas où $n=2$. Pour cela on va utiliser la différence de deux entiers, cf. III, p. 37.
\begin{equation*}
\begin{split}
\mathfrak{b}_1\mathfrak{b}_2
& =((\mathfrak{b}_1-1)+1)((\mathfrak{b}_2-1)+1) \\
& =(\mathfrak{b}_1-1)(\mathfrak{b}_2-1)+(\mathfrak{b}_1-1)+(\mathfrak{b}_2-1)+1\\
& =((\mathfrak{b}_1-1)(\mathfrak{b}_2-1)-1)+1+(\mathfrak{b}_1+\mathfrak{b}_2-2)+1\\
& =\mathfrak{b}_1+\mathfrak{b}_2+((\mathfrak{b}_1-1)(\mathfrak{b}_2-1)-1).
\end{split}
\end{equation*}
Supposons que le résultat est vraie pour $n$. On a donc $$\sum\limits_{i=1}^{n+1}\mathfrak{b}_i\leq(\sum\limits_{i=1}^n\mathfrak{b}_i)\mathfrak{b}_n\leq\prod\limits_{i=1}^{n+1}\mathfrak{b}_i.$$

Si $I$ est infini ou l'un des $\mathfrak{b}_i$ est infini. D'après le cor. 1 de la prop. 14 de III, p. 30, $\mathfrak{b}_{i}\leq\underset{i\in I}{\p}\mathfrak{b}_{i}$ pour tout $i$. D'après cette prop. 14, la prop. 10 de III, p. 28 et la prop. 12 de III, p. 29, $$\card(I)\leq2^{\card(I)}=\underset{i\in I}{\p}2\leq\underset{i\in I}{\p}\mathfrak{b}_{i}.$$ Donc $$\sup(\card(I),\sup_{i\in I}{b}_i)\leq\underset{i\in I}{\p}\mathfrak{b}_{i}.$$

Or, d'après la même prop. 14 et le cor. 2 de la prop. 6 de III, p. 27, $$\sum\limits_{i\in I}\mathfrak{b}_i\leq\sum\limits_{i\in I}\sup_{i\in I}{b}_i=\card(I)\sup_{i\in I}{b}_i.$$ D'après le cor. 4 du th. 2 de III, p. 47, $$\card(I)\sup_{i\in I}{b}_i=\sup(\card(I),\sup_{i\in I}{b}_i).$$ Ce qui achève la démonstration de a).

b) C'est un des théorèmes de J. König. Supposons que $\mathfrak{a}_i<\mathfrak{b}_i$ pour tout $i\in I$. D'après a) il suffit de montrer que $\sum\limits_{i\in I}\mathfrak{a}_i\neq\prod\limits_{i\in I}\mathfrak{b}_i$. Supposons le contraire. Donc $\prod\limits_{i\in I}\mathfrak{b}_i$ est réunion d'une famille d'ensembles (mutuellement disjoints) $(A_i)_{i\in I}$ telle que $\card (A_i)=\mathfrak{a}_i$ pour tout $i\in I$. Soit $i\in I$. L'application $\pr_i:\prod\limits_{i\in I}\mathfrak{b}_i\to\mathfrak{b}_i$ définie une surjection de $A_i$ sur $\pr_i(A_i)$. Donc $\card(\pr_i(A_i))\leq\card (A_i)=\mathfrak{a}_i<\mathfrak{b}_i$. Par suite $\pr_i(A_i)\neq\mathfrak{b}_i$. Soit $x_i\in\mathfrak{b}_i-\pr_i(A_i)$. On a $(x_i)_{i\in I}\in\prod\limits_{i\in I}\mathfrak{b}_i$ et $(x_i)_{i\in I}\not\in\bigcup\limits_{i\in I}A_i$ (car sinon on aurait $(x_i)_{i\in I}\in A_j$ pour un certain $j\in I$ et donc $x_j\in\pr_j(A_j)$). Ce qui est absurde.

\textbf{Remarque.} Le théorème de J. König a des applications:
\begin{enumerate}
  \item Le théorème de Cantor III, p. 30, Th. 2 (ou \cite[Théorème 1.7.7]{Schwartz:1991}) (on prend $\mathfrak{a}_i=1$ et $\mathfrak{b}_i=2$ et on utilise III, p. 27, Cor. 2 de la prop. 6 et III, p. 28, Prop. 10).
  \item Soit $(\mathfrak{b}_n)_{n\in\N}$ une suite de cardinaux telle que $0<\mathfrak{b}_0<\mathfrak{b}_1<\mathfrak{b}_2<\cdots$. On a $\sum\limits_{n\in\N}\mathfrak{b}_n<\prod\limits_{n\in \N}\mathfrak{b}_n$ (on prend $\mathfrak{a}_0=0$, $\mathfrak{a}_n=\mathfrak{b}_{n-1}$ pour $n\geq1$ et on utilise III, p. 27, Prop. 6 et III, p. 26, Prop. 5, a) en considérant l'application bijective $\N\to\N^*$ définie par $n\mapsto n+1$).
\end{enumerate}

\vspace{0.5cm}
6) On se servira pas de l'indication. Soit $E$ un ensemble. Si on avait $\mathfrak{P}(E)\subset E$, on aurait $\card(\mathfrak{P}(E))\leq\card(E)$. À l'aide du résultat de Cantor III, p. 30, Th. 2, on obtiendrait $\card(E)<\card(\mathfrak{P}(E))\leq\card(E)$, et donc $\card(\mathfrak{P}(E))=\card(E)$, ce qui conduirait à une contradiction. Donc $\mathfrak{P}(E)\not\subset E$, c'est-à-dire, il existe $X\subset E$ tel que $X\not\in E$.

\vspace{2cm}
\begin{center}
\S 4
\end{center}

2) La condition est nécessaire d'après III, p. 34, Cor. 2 de la prop. 3. Inversement, soit $\mathfrak{F}(E)$ l'ensemble des parties finies de $E$. Par hypothèse, cet ensemble admet un élément maximal pour la relation d'inclusion. Soit $M$ un tel élément. On a $E=M$, car sinon il existerait $x\in E-M$, et par suite $M\cup\{x\}$ serait finie et contenant strictement $M$.

\vspace{0.5cm}
3) On se servira pas de l'indication. Soit $E$ un ensemble un ensemble bien ordonné tel que, munissant $E$ de l'ordre opposé, il est aussi bien ordonné.  D'après le th. 3 de III, p. 21, l'une au moins des deux propositions suivantes est vraie:
\begin{enumerate}[a)]
  \item il existe un isomorphisme et un seul de $E$ sur un segment de $\N$;
  \item il existe un isomorphisme et un seul de $\N$ sur un segment de $E$.
\end{enumerate}

D'après la prop. 1 de III, p. 16, les segments de l'ensemble bien ordonné $\N$ des entiers, sont $\N$ et les entiers, car $S_n=]\leftarrow,n[=n$ pour tout $n\in\N$. En outre, il n'existe pas un isomorphisme de $\N$ sur un segment de $E$ car ce muni de l'ordre induit par l'ordre de $E$ vérifie la même condition que vérifie $E$, donc aussi $\N$, en particulier $\N$ admet un plus grand élément, ce qui conduirait à une contradiction. Finalement, $E$ est en bijection avec un entier, c'est-à-dire, $E$ est fini (cf. \cite[Définition 1.4.17]{Schwartz:1991}).

\vspace{2cm}
\begin{center}
\S 6
\end{center}

1) Soit $E$ un ensemble. Supposons que $E$ est infini. Soient $f$ une application de $E$ dans lui-même, $a\in E$ et $(x_n)_{n\in\N}$ la suite d'éléments de $E$ définie par $x_0=a$ et $x_{n+1}=f(x_n)$ pour tout $n\in\N$. L'ensemble $S=\{x_n\mid n\in\N\}$ est stable par $f$, c'est-à-dire, $f(S)\subset S$. Si $S\neq E$, $S$ répond à la question. Maintenant, supposons que $S=E$. La suite $(x_n)_{n\in\N}$ est injective car s'il existait $n,k\in\N$ tels que $k\geq 1$ et $x_n=x_{n+k}$, on aurait $E=S=\{x_0,\dots,x_n,x_{n+1},\dots, x_{n+k-1}\}$ qui serait fini. Donc $S-\{a\}$ répond à la question. Inversement, supposons que $E$ est fini de cardinal $n$. Soit la bijection $f:\{0,\dots,n-1\}\to\{0,\dots,n-1\}$ telle que $f(i)=i+1$ pour $0\leq i<n-1$ et $f(n-1)=0$. Soit $S$ une partie non vide et propre de $\{0,\dots,n-1\}$ et soit $k$ son plus grand élément. On a $k<n-1$ car si $k=n-1$, on aurait $S=\{0,\dots,n-1\}$. Donc $S$ n'est pas stable par $f$ car $f(k)=k+1\not\in S$ (cf. C12 (I, p. 26) et C17 (I, p. 28)).

\vspace{0.5cm}
2) D'abord $\mathfrak{c}$ et $\mathfrak{d}$ sont non nuls. On va distinguer les cas suivants:
\begin{itemize}
  \item Cas 1: $\mathfrak{c}$ et $\mathfrak{d}$ sont finis. Donc $\mathfrak{a}$ et $\mathfrak{b}$ le sont aussi. Le résultat est immédiat de la prop. 3 de III, p. 36.
  \item Cas 2: l'un des $\mathfrak{c}$ et $\mathfrak{d}$ est infini. D'après le cor. 4 du th. 2 de III, p. 47, $\mathfrak{c}+\mathfrak{d}=\mathfrak{c}\mathfrak{d}=\sup(\mathfrak{c},\mathfrak{d})$.
  \begin{enumerate}[i)]
    \item Si $\mathfrak{a}$ et $\mathfrak{b}$ sont fini. Donc $\mathfrak{a}+\mathfrak{b}$ et $\mathfrak{a}\mathfrak{b}$ le sont aussi. C'est évident.
    \item Si l'un des $\mathfrak{a}$ et $\mathfrak{b}$ est nul. C'est évident car $\mathfrak{a}<\sup(\mathfrak{c},\mathfrak{d})$ et $\mathfrak{b}<\sup(\mathfrak{c},\mathfrak{d})$.
    \item Si $\mathfrak{a}$ et $\mathfrak{b}$ sont non nuls et l'un des deux est infini. Donc, $\mathfrak{a}+\mathfrak{b}=\mathfrak{a}\mathfrak{b}=\sup(\mathfrak{a},\mathfrak{b})$ (cf. le cor. 4 du th. 2 de III, p. 47). On a $\mathfrak{a}<\mathfrak{c}\leq\sup(\mathfrak{c},\mathfrak{d})$ et $\mathfrak{b}<\mathfrak{d}\leq\sup(\mathfrak{c},\mathfrak{d})$. D'où $\sup(\mathfrak{a},\mathfrak{b})<\sup(\mathfrak{c},\mathfrak{d})$.
  \end{enumerate}
\end{itemize}

\vspace{0.5cm}
9) Soient $a$ le plus petit élément de $E$ et $f$ l'application de $E$ dans lui-même définie par: $f(a)=a$ et $f(x)$ est le plus grand élément de $S_x=]\leftarrow,x[$ pour tout $x\in E-\{a\}$. Soient $x\in E-\{a\}$ et $(x_n)_{n\in\N}$ la suite d'éléments de $E$ définie par $x_0=x$ et $x_{n+1}=f(x_n)$ pour tout $n\in\N$. Cette suite est décroissante, donc, d'après le cor. 1 de III, p. 51, elle stationnaire. On a $$S_x=]\leftarrow,x_1]=
\begin{cases}\{x_1\} & \text{si}\quad x_1=a\\
]\leftarrow,x_2]\cup\{x_1\} & \text{si}\quad x_1\neq a
\end{cases}=
\begin{cases}\{x_1\} & \text{si}\quad x_1=a\\
\{x_1,x_2\} & \text{si}\quad x_1\neq a \text{ et }x_2=a\\
]\leftarrow,x_3]\cup\{x_1,x_2\}& \text{si}\quad x_1\neq a \text{ et }x_2\neq a
\end{cases}=\cdots $$

D'où $S_x=\{x_1,\dots,x_k\}$ où $k$ est le plus petit entier tel que $x_k=a$. Finalement, d'après la prop. 1 de III, p. 16, tout segment $\neq E$ est fini. Le th. 3 de III, p. 21, achève donc la démonstration.

\appendix
\chapter{Notes concernant le chapitre I}\label{app1}

I, p. 1: Une théorie mathématique peut être sans signe spécifique, cf. I, p. 47, \S1, exerc. 1, et sa solution.

\vspace{0.5cm}
I, pp. 17, 18: La juxtaposition de deux constructions formatives est aussi une construction formative.

\vspace{0.5cm}
I, pp. 21, 22: Soit $\mathcal{T}$ une théorie mathématique. Il se peut que $\mathcal{T}$ soit sans axiome explicite, cf. I, p. 30. Si $\mathcal{T}$ est sans axiome explicite et sans axiome implicite, alors $\mathcal{T}$ est sans théorème! .

En outre, les axiomes explicites et implicites de $\mathcal{T}$ sont des théorèmes de $\mathcal{T}$.

\vspace{0.5cm}
I, p. 25: Soit $\mathcal{T}$ une théorie logique contadictoire. Soit $\bm{A}$ une relation telle que $\bm{A}$, $\text{non}\bm{A}$ sont des théorèmes. Donc toute relation $\bm{B}$ de $\mathcal{T}$ est un théorème. En effet, d'après S2, $(\text{non}\bm{A})\Rightarrow ((\text{non}\bm{A})\text{ ou }\bm{B})$ c'est-à-dire, $(\text{non}\bm{A})\Rightarrow (\bm{A}\Rightarrow\bm{B})$, est un théorème. On conclut, à l'aide de C1, que $\bm{B}$ est un théorème.

\vspace{0.5cm}
I, p. 31: Pour les démonstrations de C22, C23, C24 et C25 cf. la solution de l'exerc. 1 de I, p. 47, \S3.

\chapter{Notes concernant le chapitre II}

Je crois que N. Bourbaki a utilisé au début, dans ce livre, la terminologie ``fonction'' ou ``application de\dots dans\dots'' pour désigner ce qu'on appelle ``fonction'' ou ``application'' et ``application de\dots sur\dots'' pour désigner ce qu'on appelle ``application surjective''. Donc il faut remplacer dans plusieurs situations ``application de\dots sur\dots'' par ``application surjective'':

\begin{itemize}
\item II, p. 20, Prop. 9: a) Soient $E$, $F$, $G$ des ensembles, $g$ une application surjective de $E$ dans $F$\dots
\item II, p. 33, Prop. 5: \dots la projection $\pr_J$ est une application surjective\dots
\item II, p. 34, Corollaire 1: \dots la projection $\pr_\alpha$ est une application surjective\dots
\item III, p. 18, C60: \dots Il existe un ensemble $U$ et une application surjective $f$ de $E$ sur $U$\dots et dans sa démonstration: \dots il existe un ensemble $U_S$ et une application surjective $f_S$ de $S$ sur $U_S$\dots
\item III, p. 19, 7ème ligne:\dots pour toute application surjective $h$ d'un sement de $E$ sur une partie de $F$,\dots
\item III, p. 26, 11ème ligne: En effet, il existe une application surjective de la somme des $E_\iota$ sur leur réunion\dots
\item III, p. 50, 4ème ligne: Proposition 4.--- Soit $f$ une application surjective d'un ensemble $E$ sur un ensemble infini $F$, \dots

etc.
\end{itemize}

\section*{\S1. RELATIONS COLLECTIVISANTES}

\subsection*{4. Relations collectivisantes}

Pour un cadre plus général cf. p. 9.

\subsection*{5. L'axiome de l'ensemble à deux éléments}

Soient $\bm{T}$ et $\bm{U}$ deux termes. On dit que le terme obtenu par remplacement simultané de $x$ par $\bm{T}$ et de $y$ par $\bm{U}$ dans $\{x,y\}$, est \emph{l'ensemble à deux éléments $\bm{T}$ et $\bm{U}$}, et on le note par $\{\bm{T},\bm{U}\}$.

Soit $\bm{z}$ une lettre ne figurant ni dans $\bm{T}$ ni dans $\bm{U}$ (donc non plus dans $\{\bm{T},\bm{U}\}$) et $\bm{V}$ un terme. On sait que $z\in\{x,y\}\Leftrightarrow (z=x \text{ ou } z=y)$ est un théorème. D'après C3 ( I, p. 23), si on remplace simultanément $x$ par $\bm{T}$, $y$ par $\bm{U}$ et $z$ par $\bm{z}$ (resp. $\bm{V}$) dans cette équivalence on obtient un théorème
$$\bm{z}\in\{\bm{T},\bm{U}\}\Leftrightarrow (\bm{z}=\bm{T} \text{ ou } \bm{z}=\bm{U})$$
(resp.
$$\bm{V}\in\{\bm{T},\bm{U}\}\Leftrightarrow (\bm{V}=\bm{T} \text{ ou } \bm{V}=\bm{U})$$
) (voir I, p. 16 et II, p. 2). Donc la relation ``$(\bm{z}=\bm{T}) \text{ ou } (\bm{z}=\bm{U})$'' est collectivisante en $\bm{z}$ et par l'axiome d'extensionalité
$$\{\bm{T},\bm{U}\}=\{\bm{z}\mid (\bm{z}=\bm{T}) \text{ ou } (\bm{z}=\bm{U})\}.$$
Lorsque $\bm{T}$ et $\bm{U}$ sont identiques, $\{\bm{T}\}$ est par définition $\{\bm{T},\bm{T}\}$. Donc $\{\bm{T}\}$ est identique à $(\bm{T}\mid x)\{x\}$.

\subsection*{7. Complémentaire d'un ensemble. L'ensemble vide}

Soient $\bm{X}$ un terme, $\bm{A}$ une partie de $\bm{X}$ et $\bm{x}$ une lettre ne figurant pas dans $\bm{X}$ ni dans $\bm{A}$. La relation ``$\bm{x}\notin\bm{A}\text{ et } \bm{x}\in\bm{X}$'' est collectivisante en $\bm{x}$ d'après C51.

\begin{definition}
On appelle \emph{complémentaire} de $\bm{A}$ par rapport à $\bm{X}$, et on le note par $\complement_{\bm{X}}\bm{A}$, le terme $\{\bm{x}\mid\bm{x}\notin\bm{A}\text{ et } \bm{x}\in\bm{X}\}$ qui est une partie de $\bm{X}$ (les lettres figurant dans $\complement_{\bm{X}}\bm{A}$ sont celles figurant dans $\bm{X}$ ou $\bm{A}$).
\end{definition}

La relation $\bm{A}=\complement_{\bm{X}}(\complement_{\bm{X}}\bm{A})$ est un théorème d'après (16), p. 4, C10 et (5), p. 4. Si $\bm{B}$ une partie de $\bm{X}$, $(\bm{A}\subset\bm{B})\Leftrightarrow(\complement_{\bm{X}}\bm{B}\subset\complement_{\bm{X}}\bm{A})$ est un théorème. Il suffit de montrer la première implication, or celle-ci est immédiate d'après C12 et (3), p. 3 en prenant une lettre $\bm{x}$ ne figurant pas dans $\bm{A}$, ni dans $\bm{B}$, ni dans $\bm{X}$.

\begin{theorem}
La relation $(\forall x) (x\notin X)$ est fonctionnelle en $X$ (D'après CS8 on peut prendre au lieu de $x$ n'importe qu'elle lettre qui soit distincte de $X$).
\end{theorem}

On désigne par $\varnothing$ le terme correspondant à cette relation fonctionnelle: $\bm{\tau}_X((\forall x) (x\notin X))$, et on l'appelle \emph{l'ensemble vide}. Soient $\bm{X}$ un terme et $\bm{x}$ une lettre ne figurant pas dans $\bm{X}$. On a $(\forall \bm{x}) (\bm{x}\notin \bm{X})$ est équivalente \'{a} $\bm{X}=\varnothing$. Si c'est le cas on dit que \emph{l'ensemble $\bm{X}$ est vide}. Puisque $\varnothing=\varnothing$ est $(\forall x) (x\notin \varnothing)$ est un théorème. On a aussi les théorèmes $\varnothing\subset\bm{X}$, $\complement_{\bm{X}}\bm{X}=\varnothing$, et d'après C52 et A1, $\complement_{\bm{X}}\varnothing=\bm{X}$.

\section*{\S2. COUPLES}

\subsection*{1. Définition des couples}

Soient $\bm{T}$ et $\bm{U}$ deux termes. On dit que le terme obtenu par remplacement simultané de $x$ par $\bm{T}$ et de $y$ par $\bm{U}$ dans $(x,y)$ (qui est par définition $\{\{x\},\{x,y\}\}$), est le couple (ordonné) formé de $\bm{T}$ et de $\bm{U}$, et on le note par $(\bm{T},\bm{U})$.

Soit $\bm{z}$ une lettre ne figurant pas ni dans $\bm{T}$ ni dans $\bm{U}$ (donc ni dans $\{\bm{T}\}$, ni dans $\{\bm{T},\bm{U}\}$, ni dans $\{\{\bm{T}\},\{\bm{T},\bm{U}\}\}$).On sait que $z\in\{\{x\},\{x,y\}\}\Leftrightarrow (z=\{x\} \text{ ou } z=\{x,y\})$ est un théorème. D'après C3 (I, p. 23), si on remplace simultanément $x$ par $\bm{T}$, $y$ par $\bm{U}$ et $z$ par $\bm{z}$ dans cette équivalence on obtient un théorème
$$\bm{z}\in (\bm{T},\bm{U})\Leftrightarrow (\bm{z}=\{\bm{T}\}\text{ ou }\bm{z}=\{\bm{T},\bm{U}\}) .$$
D'où
$$(\bm{T},\bm{U})=\{\{\bm{T}\},\{\bm{T},\bm{U}\}\}.$$
Soient $\bm{T}$, $\bm{T}'$, $\bm{U}$ et $\bm{U}'$ des termes. Encore par C3 (I, p. 23), la prop. 1 de II, p. 7, se généralise ainsi
$$(\bm{T},\bm{U})=(\bm{T}',\bm{U}')\Leftrightarrow (\bm{T}=\bm{T}'\text{ et }\bm{U}=\bm{U}')$$ est un théorème (par remplacement simultané de $x$ par $\bm{T}$, de $y$ par $\bm{U}$, de $x'$ par $\bm{T}'$ et de $y'$ par $\bm{U}'$).

Soient $\bm{T}$ un terme et $\bm{x}$, $\bm{y}$ deux lettres distinctes entre elles et ne figurant pas dans $\bm{T}$. La relation $(\exists\bm{x})(\exists\bm{y})(\bm{T}=(\bm{x},\bm{y}))$ se désigne par ``$\bm{T}$ \emph{est un couple}''.

Si $\bm{T}$ est un couple, la relation $(\exists\bm{y})(\bm{T}=(\bm{x},\bm{y}))$ est fonctionnelle en $\bm{x}$ (elle univoque en $\bm{x}$ d'après les généralisations du th. 3 de I, p. 40, et de la prop. 1 de II, p. 7). De même la relation $(\exists\bm{x})(\bm{T}=(\bm{x},\bm{y}))$ est fonctionnelle en $\bm{y}$. Les symboles fonctionels $\bm{\tau}_{\bm{x}}((\exists\bm{y})(\bm{T}=(\bm{x},\bm{y})))$ et $\bm{\tau}_{\bm{y}}((\exists\bm{x})(\bm{T}=(\bm{x},\bm{y})))$ sont appelés respectivement \emph{première coordonée} (ou \emph{première projection}) et \emph{seconde coordonée} (ou \emph{seconde projection}) et on les note respectivement par $\pr_1\bm{T}$ et $\pr_2\bm{T}$.

\begin{remark}
Soient $\bm{z}$ une lettre et $\bm{x}$ et $\bm{y}$ deux lettres distinctes entre elles et distinctes de $\bm{z}$. La relation $(\exists\bm{x})(\exists\bm{y})(\bm{z}=(\bm{x},\bm{y}))$ est identique à la relation obtenue en rempla\c{c}ant $\bm{x}$ par $\bm{y}$ et $\bm{y}$ par $\bm{x}$.

Soient $\bm{x}'$ et $\bm{y}'$ deux lettres distinctes entre elles et distinctes de $\bm{z}$. Les relations $(\exists\bm{x})(\exists\bm{y})(\bm{z}=(\bm{x},\bm{y}))$, $(\exists\bm{x}')(\exists\bm{y}')(\bm{z}=(\bm{x}',\bm{y}'))$ sont identiques. En effet,

1er cas. $\bm{x}$ et $\bm{x}'$ sont identiques et $\bm{y}$ et $\bm{y}'$ sont de même. C'est trivial.

2ème cas. $\bm{x}'$ et $\bm{y}'$ sont distinctes de $\bm{x}$ et $\bm{y}$. On le voit en utilisant CS8 et CS9.

3ème cas. L'une des $\bm{x}'$ et $\bm{y}'$ ent identique à $\bm{x}$ ou bien à $\bm{y}$, et l'autre lettre est distincte de la lettre restante. Par exemple, supposons que $\bm{y}'$ et $\bm{y}$ sont identiques et $\bm{x}'$ est distincte de $\bm{x}$. Donc $\bm{y}$ est distincte de $\bm{x}$ et de $\bm{x}'$ (qui sont distinctes entre elles), et tous distinctes de $\bm{z}$. On le voit de même en utilisant CS8 et CS9.

On a la même chose concernant les termes $\pr_1\bm{z}$ et $\pr_2\bm{z}$ (On utilise CS3, CS8 et CS9).
\end{remark}

On sait que  $z=(x,y)\Leftrightarrow (z \text{ est un couple et } x=\pr_1z\text{ et } y=\pr_2z)$ est un théorème. Soient $\bm{T}$, $\bm{U}$ et $\bm{V}$ des termes. Encore par C3 (I, p. 23),
\begin{equation}\label{couple1}
\bm{V}=(\bm{T},\bm{U})\Leftrightarrow (\bm{V} \text{ est un couple et } \bm{T}=\pr_1\bm{V}\text{ et } \bm{U}=\pr_2\bm{V})
\end{equation}
est aussi un théorème (par remplacement simultané de $x$ par $\bm{T}$, de $y$ par $\bm{U}$ et de $z$ par $\bm{V}$).

En particulier,
\begin{equation}\label{couple2}
(\bm{T},\bm{U})\text{ est un couple et } \bm{T}=\pr_1(\bm{T},\bm{U})\text{ et } \bm{U}=\pr_2(\bm{T},\bm{U}))
\end{equation}
est un théorème (d'après la généralisation du th. 1 de I, p. 39, aux termes).

\subsection*{2. Produit de deux ensembles}

Avec une démonstration analogue à celle du th. 1 de II, p. 8, on montre:

\begin{theorem}
Soient $\bm{x}$, $\bm{y}$, $\bm{z}$, $\bm{t}$, $\bm{u}$ et $\bm{v}$ des lettres deux à deux distintes entre elles. La relation
$$(\forall \bm{u})(\forall \bm{v})(\exists \bm{t})(\forall \bm{z})(( \bm{z}\in\bm{t})\Leftrightarrow \underbrace{(\exists \bm{x})(\exists \bm{y})(\bm{z}=(\bm{x},\bm{y}) \text{ et } \bm{x}\in \bm{u}\text{ et } \bm{y}\in \bm{v})})$$
est vraie. Autrement dit, quels que soit $\bm{u}$ et $\bm{v}$, la relation (équivalente à la relation entre l'acolade) ``$\bm{z} \text{ est un couple et } \pr_1\bm{z}\in\bm{u}\text{ et } \pr_2\bm{z}\in\bm{v}$'' est collectivisante en $\bm{z}$.
\end{theorem}

\begin{definition}
Soient $\bm{X}$, $\bm{Y}$ deux ensembles, et $\bm{x}$, $\bm{y}$, $\bm{z}$ des lettres deux à deux distinctes entre elles et ne figurant pas dans  $\bm{X}$ ni dans $\bm{Y}$. L'ensemble
$$\{\bm{z}\mid (\exists \bm{x})(\exists \bm{y})(\bm{z}=(\bm{x},\bm{y}) \text{ et } \bm{x}\in \bm{X}\text{ et } \bm{y}\in \bm{Y})\}=\{\bm{z}\mid\bm{z} \text{ est un couple et } \pr_1\bm{z}\in\bm{X}\text{ et } \pr_2\bm{z}\in\bm{Y}\}$$
s'appelle le \emph{produit cartésien} de $\bm{X}$ et $\bm{Y}$ et se désigne par $\bm{X}\times\bm{Y}$.
\end{definition}
On a donc
\begin{equation}\label{produit1}
\bm{z}\in \bm{X}\times\bm{Y}\Leftrightarrow(\bm{z} \text{ est un couple et } \pr_1\bm{z}\in\bm{X}\text{ et } \pr_2\bm{z}\in\bm{Y})
\end{equation}
est un théorème.
Les ensembles $\bm{X}$ et $\bm{Y}$ s'appellent le \emph{premier} et le \emph{second ensemble facteur} de $\bm{X}\times\bm{Y}$.

D'après \eqref{produit1} et C3 (I, p. 23), si $\bm{T}$ est un ensemble,
\begin{equation}\label{produit2}
\bm{T}\in \bm{X}\times\bm{Y}\Leftrightarrow(\bm{T} \text{ est un couple et } \pr_1\bm{T}\in\bm{X}\text{ et } \pr_2\bm{T}\in\bm{Y})
\end{equation}
est un théorème.

\begin{proposition}\label{II,8,prod2}
Soient $\bm{A}$ et $\bm{B}$ deux ensembles \underline{non vides}. La relation $\bm{A}'\times\bm{B}'\subset\bm{A}\times\bm{B}$ est équivalente \'{a} ``$\bm{A}'\subset\bm{A}\text{ et }\bm{B}'\subset\bm{B}$''.
\end{proposition}

On la montre par utilisation de C14.

$\Leftarrow$. Soit $\bm{z}$ une lettre ne figurant pas dans $\bm{A}$, ni dans $\bm{A}'$, ni dans $\bm{B}$, ni dans $\bm{B}'$. Les lettres figurant dans $\bm{A}\times\bm{B}$ sont celles qui figurent dans $\bm{A}$ ou dans $\bm{B}$. On a d'après \eqref{produit1}, $\bm{z}\in\bm{A}'\times\bm{B}'\Rightarrow\bm{z}\in\bm{A}\times\bm{B}$ est un théorème. D'où l'implication.

$\Rightarrow$. Soit $\bm{x}$ une lettre ne figurant pas dans $\bm{A}$ ni dans $\bm{A}'$. D'après C14, \eqref{couple2} et \eqref{produit2}, on a ``$\bm{x}\in\bm{A}'\Rightarrow\bm{x}\in\bm{A}$'' est un théorème. Donc $\bm{A}'\subset\bm{A}$ l'est encore. De même on montre $\bm{B}'\subset\bm{B}$.

\begin{proposition}\label{prod=vide}
Soient $\bm{A}$ et $\bm{B}$ deux ensembles. La relation $\bm{A}\times\bm{B}=\varnothing$ est équivalente \'{a} ``$\bm{A}=\varnothing\text{ ou }\bm{B}=\varnothing$''.
\end{proposition}

Il (faut et il) suffit de montrer que la relation $\bm{A}\times\bm{B}\neq\varnothing$ est équivalente \'{a} ``$\bm{A}\neq\varnothing\text{ et }\bm{B}\neq\varnothing$'' or cette équivalence est immédiate d'après C14, \eqref{produit2} et \eqref{couple2}.

\begin{remark}
la prop. \ref{II,8,prod2} n'est pas vraie en général si l'une des $\bm{A}$ et $\bm{B}$ est vide. Il suffit de prendre $\bm{A}=\bm{A}'=\varnothing$ et $\bm{B}$ et $\bm{B}'$ deux ensembles tels que $\bm{B}'\not\subset\bm{B}$.
\end{remark}

\section*{\S3. CORRESPONDANCES}

\subsection*{1. Graphes et correspondances}

\begin{definition}
Soit $\bm{G}$ un ensemble et $\bm{z}$ une lettre ne figurant pas dans $\bm{G}$. On dit que $\bm{G}$ est \emph{un graphe} si tout élément de $\bm{G}$ est un couple, autrement dit si la relation $$(\forall\bm{z})(\bm{z}\in\bm{G}\Rightarrow (\bm{z}\text{ est un couple}))$$ est vraie.
\end{definition}

Si $\bm{G}$ est un graphe, la relation $(\bm{x},\bm{y})\in\bm{G}$ s'exprime encore en disant que ``$\bm{y}$ correspond à $\bm{x}$ par $\bm{G}$''.

Soit $\bm{R}\{\bm{x},\bm{y}\}$ une relation, $\bm{x}$ et $\bm{y}$ étant des lettres distinctes. S'il existe un graphe $\bm{G}$ où ne figurent ni $\bm{x}$ ni $\bm{y}$, tel que $(\forall\bm{x})(\forall\bm{y})(\bm{R}\Leftrightarrow(\bm{x},\bm{y})\in\bm{G})$ soit vraie, on dit que $\bm{R}$ \emph{admet un graphe} (par rapport aux lettres $\bm{x}$ et $\bm{y}$), autrement dit (S5), si la relation $$(\bm{g}\text{ est un graphe) et }(\forall\bm{x})(\forall\bm{y})(\bm{R}\Leftrightarrow(\bm{x},\bm{y})\in\bm{g})$$ où $\bm{g}$ est une lettre distincte de $\bm{x}$, de $\bm{y}$ et ne figurant pas dans $\bm{R}$. D'après la remarque \ref{=2graphes} ci-dessous, $\bm{G}$ est unique, et s'appelle le \emph{graphe} de $\bm{R}$ (ou l'\emph{ensemble représentatif} de $\bm{R}$) par rapport à $\bm{x}$ et $\bm{y}$.

\begin{remark}\label{=2graphes}
Soient $\bm{G}$ et $\bm{G}'$ deux graphes (chacun peut être vide), $\bm{x}$ et $\bm{y}$ deux lettres ne figurant pas dans $\bm{G}$ ni dans $\bm{G}'$.

Si $(\forall\bm{x})(\forall\bm{y})((\bm{x},\bm{y})\in\bm{G}\Rightarrow(\bm{x},\bm{y})\in\bm{G}')$ est vraie, on a $\bm{G}\subset\bm{G}'$ (on distingue deux cas $\bm{G}=\varnothing$, $\bm{G}\neq\varnothing$).

Si $(\forall\bm{x})(\forall\bm{y})((\bm{x},\bm{y})\in\bm{G}\Leftrightarrow(\bm{x},\bm{y})\in\bm{G}')$ est vraie, par l'axiome d'extensionalité, on a $\bm{G}=\bm{G}'$.
\end{remark}

\begin{remark}
Soient $\bm{T}$ un terme et $\bm{x}'$ et $\bm{y}'$ deux lettres distinctes entre elles, distinctes de $\bm{x}$, de $\bm{y}$ et de $\bm{z}$ et ne figurant ni dans $\bm{R}$ ni dans $\bm{T}$. La relation $(\bm{T}\mid\bm{z})\bm{R}\{\pr_1\bm{z},\pr_2\bm{z}\}$

\noindent est identique à $(\bm{T}\mid\bm{z})(\pr_1\bm{z}\mid\bm{x}')(\pr_2\bm{z}\mid\bm{y}')(\bm{x}'\mid\bm{x})(\bm{y}'\mid\bm{y})\bm{R}$ (I, p. 16);

\noindent est identique à $(\pr_1\bm{T}\mid\bm{x}')(\pr_2\bm{T}\mid\bm{y}')(\bm{T}\mid\bm{z})(\bm{x}'\mid\bm{x})(\bm{y}'\mid\bm{y})\bm{R}$ (CS2);

\noindent est identique à
$(\pr_1\bm{T}\mid\bm{x}')(\pr_2\bm{T}\mid\bm{y}')(\bm{x}'\mid\bm{x})(\bm{y}'\mid\bm{y})\bm{R}$;

\noindent est identique à
$\bm{R}\{\pr_1\bm{T},\pr_2\bm{T}\}$.
\end{remark}

Maintenant, soit $\bm{T}$ un ensemble, où ne figurent ni $\bm{x}$ ni $\bm{y}$, tel que $\bm{R}\Rightarrow ((\bm{x},\bm{y})\in\bm{T})$ soit vraie. On a $\bm{R}$ admet un graphe. En effet, soit $\bm{z}$ une lettre distincte de $\bm{x}$, de $\bm{y}$ et ne figurant pas ni dans $\bm{R}$ ni dans $\bm{T}$. Soit $\bm{G}$ l'ensemble $\{\bm{z}\mid (\bm{z}\text{ est un couple) et } \bm{z}\in\bm{T}\text{ et } \bm{R}\{\pr_1\bm{z},\pr_2\bm{z}\}\}$. Évidemment $\bm{G}$ est un graphe et d'après la remarque 11, si $\bm{U}$, $\bm{V}$ deux ensembles,
$$(\bm{U},\bm{V})\in\bm{G}\Leftrightarrow ((\bm{U},\bm{V})\in\bm{T}\text{ et }\bm{R}\{\bm{U},\bm{V}\})$$ est vraie. Donc (I, p. 16) $\bm{G}$ est le graphe de $\bm{R}$.

Soit $\bm{G}$ un graphe. Soient $\bm{z}$ et $\bm{t}$ deux lettres distinctes et ne figurant pas dans $\bm{G}$. Les ensembles (II, p. 6) $$\{\bm{t}\mid (\exists\bm{z})(\bm{t}=\pr_1\bm{z}\text{ et }\bm{z}\in\bm{G})\},\{\bm{t}\mid (\exists\bm{z})(\bm{t}=\pr_2\bm{z}\text{ et }\bm{z}\in\bm{G})\},$$ s'appellent \emph{la première} et \emph{la seconde projection} du graphe $\bm{G}$ ou encore \emph{l'ensemble de définition} et \emph{l'ensemble des valeurs} de $\bm{G}$, et on le désigne par $\pr_1\langle\bm{G}\rangle$ et $\pr_2\langle\bm{G}\rangle$ (ou $\pr_1\bm{G}$ et $\pr_2\bm{G}$ lorsqu'aucune confusion n'en résulte).

\begin{remark}\label{C53}
On reprend l'énoncé de C53 (II, pp. 5, 6). Soit $\bm{B}$ l'ensemble $\{\bm{y}\mid (\exists\bm{x})(\bm{y}=\bm{T}\text{ et }\bm{x}\in\bm{A})\}$. Soient $\bm{Y}$ un ensemble et $\bm{x}'$ une lettre ne figurant ni dans $\bm{A}$, ni dans $\bm{T}$, ni dans $\bm{Y}$. On a $$\bm{Y}\in\bm{B}\Leftrightarrow (\exists\bm{x}')(\bm{Y}=(\bm{x}'\mid\bm{x})\bm{T}\text{ et }\bm{x}'\in\bm{A})$$ (on prend $\bm{x}'$ distincte de $\bm{x}$ et le cas contraire). Ce qui revient \'{a} trouver un ensemble $\bm{X}$ vérifiant $\bm{Y}=(\bm{X}\mid\bm{x})\bm{T}\text{ et }\bm{X}\in\bm{A}$ (d'après CS1).
\end{remark}

Donc si $\bm{T}$ un ensemble et $\bm{y}$ une lettre ne figurant ni dans $\bm{T}$ ni dans $\bm{G}$, $$(\bm{T}\in\pr_1\bm{G})\Leftrightarrow (\exists\bm{y})((\bm{T},\bm{y})\in\bm{G})$$ est vraie. Une similaire équivalence est vraie pour $\pr_2\bm{G}$.

Soient $\bm{x}$, $\bm{y}$, $\bm{x}'$ des lettres ne figurant pas dans $\bm{G}$, tels que $\bm{y}$ est distincte de $\bm{x}$ avec $\bm{x}'$ est distincte de $\bm{x}$ et de $\bm{z}$. On a $$(\bm{x}\in\pr_1\bm{G})\Leftrightarrow (\exists\bm{x}')(\bm{x}'\in\bm{G}\text{ et } \bm{x}=\underbrace{(\bm{x}'\mid\bm{z})\pr_1\bm{z}}_{\pr_1\bm{x}'})\Leftrightarrow^{S5} (\exists\bm{y})((\bm{x},\bm{y})\in\bm{G}).$$

Donc d'après C52 et l'axiome d'extensionalité, $$\pr_1\bm{G}=\{\bm{x}\mid(\exists\bm{y})((\bm{x},\bm{y})\in\bm{G})\}.$$

On a $\bm{G}\subset (\pr_1\bm{G})\times (\pr_2\bm{G})$ (cf Remarque \ref{=2graphes}): \textbf{tout ensemble de couples est donc une partie d'un produit, et réciproquement}. Si l'une des deux ensembles $\pr_1\bm{G}$, $\pr_2\bm{G}$ est vide, on a $\bm{G}=\varnothing$ (Prop. \ref{prod=vide}).

\begin{example}
Soit $\bm{G}$ l'ensemble $\bm{A}\times\bm{B}$ où $\bm{A}$ et $\bm{B}$ sont deux ensembles. Évidemment $\bm{G}$ est un graphe.
\begin{enumerate}
  \item Si $\bm{A}\neq\varnothing$ et $\bm{B}\neq\varnothing$: $\bm{G}\neq\varnothing$ (Prop. \ref{prod=vide}), $\pr_1\bm{G}=\bm{A}$ et $\pr_2\bm{G}=\bm{B}$.
  \item Si $\bm{A}$ ou $\bm{B}$ est vide: $\bm{G}=\varnothing$ (Prop. \ref{prod=vide}), $\pr_1\bm{G}=\varnothing$ et $\pr_2\bm{G}=\varnothing$.
\end{enumerate}
\end{example}

\begin{remark}
La relation $x=y$ n'admet pas de graphe; car la première projection de ce graphe, s'il existait, serait l'ensemble de tous les objets (cf. II, p. 6, Remarque). En effet, si $\bm{G}$ un tel graphe, on aurait $(x=y)\Leftrightarrow (x,y)\in \bm{G}$. Donc $(\forall x)(x\in \pr_1\bm{G})$ serait vraie.
\end{remark}

\begin{definition}
On appelle \emph{correspondance} entre un ensemble $\bm{A}$ et un ensemble $\bm{B}$ un triplet $\Gamma=(\bm{G},\bm{A},\bm{B})$ où $\bm{G}$ est un graphe tel que $\pr_1\bm{G}\subset\bm{A}$ et $\pr_2\bm{G}\subset\bm{B}$, \textbf{autrement dit $\bm{G}$ est un ensemble tel que $\bm{G}\subset\bm{A}\times\bm{B}$}. On dit que $\bm{G}$ est le \emph{graphe} de $\Gamma$, $\bm{A}$ \emph{l'ensemble de départ} et $\bm{B}$ \emph{l'ensemble d'arrivée} de $\Gamma$.
\end{definition}

Si $(\bm{x},\bm{y})\in\bm{G}$, on dit encore que ``$\bm{y}$ correspond à $\bm{x}$ par la correspondance $\Gamma$''. Pour tout $\bm{x}\in\pr_1\bm{G}$, on dit que la correspondance $\Gamma$ est \emph{définie por l'objet} $\bm{x}$, et $\pr_1\bm{G}$ est appelé \emph{l'ensemble de définition} (ou \emph{domaine}) de $\Gamma$; pour tout $\bm{y}\in\pr_2\bm{G}$, on dit que $\bm{y}$ est \emph{une valeur prise par} et $\pr_2\bm{G}$ est appelé \emph{l'ensemble des valeurs} (ou \emph{image}) de $\Gamma$.

Si $\bm{R}\{\bm{x},\bm{y}\}$ est une relation admettant un graphe $\bm{G}$ (par rapport aux lettres $\bm{x}$ et $\bm{y}$), et si $\bm{A}$ et $\bm{B}$ sont deux ensembles tels que $\pr_1\bm{G}\subset\bm{A}$ et $\pr_2\bm{G}\subset\bm{B}$, on dit que $\bm{R}$ est une \emph{relation entre un élément de $\bm{A}$ et un élément de} $\bm{B}$ (relativement aux lettres $\bm{x}$, $\bm{y}$). On dit que la correspondance $\Gamma=(\bm{G},\bm{A},\bm{B})$ est la correspondance entre $\bm{A}$ et $\bm{B}$ \emph{définie par la relation} $\bm{R}$ (par rapport à $\bm{x}$ et $\bm{y}$).

Maintenant, soient $\bm{R}$ une relation, $\bm{T}$ un ensemble, où ne figurent ni $\bm{x}$ ni $\bm{y}$, tel que $\bm{R}\Rightarrow ((\bm{x},\bm{y})\in\bm{T})$ soit vraie. Soit $\bm{G}$ le graphe de $\bm{R}$. Soient $\bm{x}'$ ni $\bm{y}'$ sont deux lettres distinctes et ne figurant pas dans $\bm{G}$ (i.e. ne figurant ni dans $\bm{R}\{\square,\square\}$ ni dans $\bm{T}$). On a $$\pr_2\bm{G}=\{\bm{y}'\mid (\exists\bm{x}')((\bm{x}',\bm{y}')\in\bm{G})\}=\{\bm{y}'\mid (\exists\bm{x}')\bm{R}\{\bm{x}',\bm{y}'\}\}.$$

\begin{definition}
Soient $\bm{G}$ un graphe et $\bm{X}$ un ensemble. Soient $\bm{x}$ et $\bm{y}$ deux lettres distinctes et ne figurant ni dans $\bm{G}$ ni dans $\bm{X}$. La relation ``$\bm{x}\in\bm{X}\text{ et } (\bm{x},\bm{y})\in\bm{G}$'' entraîne $ (\bm{x},\bm{y})\in\bm{G}$ et admet par suite un graphe $\bm{G}'$. La seconde projection de $\bm{G}'$ se compose de tous les objets qui correspondent par $\bm{G}$ \'{a} des objets de $\bm{X}$:  $\{\bm{y}\mid (\exists\bm{x})(\bm{x}\in\bm{X}\text{ et } (\bm{x},\bm{y})\in\bm{G})\}$, s'appelle \emph{l'image de $\bm{X}$ par} $\bm{G}$ et se désigne par $\bm{G}\langle\bm{X}\rangle$ ou $\bm{G}(\bm{X})$.

Soient $\Gamma=(\bm{G},\bm{A},\bm{B})$ une correspondance, et $\bm{X}$ une partie de $\bm{A}$. L'ensemble $\bm{G}\langle\bm{X}\rangle$ se note encore $\Gamma\langle\bm{X}\rangle$ ou $\Gamma(\bm{X})$ et s'appelle \emph{l'image de $\bm{X}$ par} $\Gamma$.
\end{definition}

Soit $\bm{G}$ un graphe. Soient $\bm{x}$ et $\bm{y}$ des objets. Comme la relation $(\bm{x},\bm{y})\in\bm{G}$ entraîne $\bm{y}\in\pr_2\bm{G}$, on a $\bm{G}\langle\bm{X}\rangle\subset\pr_2\bm{G}$ pour tout ensemble $\bm{X}$; comme $(\bm{x},\bm{y})\in\bm{G}$ entraîne $\bm{x}\in\pr_1\bm{G}$, on a $\bm{G}\langle\pr_1\bm{G}\rangle=\pr_2\bm{G}$. On a  $\bm{G}\langle\varnothing\rangle=\varnothing$, puisque puisque $x\not\in\varnothing$ est un théorème. Si $\bm{X}\subset\pr_1\bm{G}$ et $\bm{X}\neq\varnothing$, on a $\bm{G}\langle\bm{X}\rangle\neq\varnothing$.

La proposition suivante est évidente:

\begin{proposition}
Soient $\bm{G}$ un graphe, $\bm{X}$ et $\bm{Y}$ deux ensembles; la relation $\bm{X}\subset\bm{Y}$ entraîne $\bm{G}\langle\bm{X}\rangle\subset\bm{G}\langle\bm{Y}\rangle$.
\end{proposition}

\begin{corollary}
Si $\bm{A}\supset\pr_1\bm{G}$, on a $\bm{G}\langle\bm{A}\rangle=\pr_2\bm{G}$.
\end{corollary}

\begin{definition}
Soient $\bm{G}$ un graphe et $\bm{x}$ un objet. On appelle \emph{coupe} de $\bm{G}$ suivant $\bm{x}$ l'ensemble $\bm{G}\langle\{\bm{x}\}\rangle$ (qu'on désigne aussi parfois $\bm{G}(\bm{x})$, par abus de notation).
\end{definition}

Si $\bm{y}$ une lettre ne figurant ni dans $\bm{G}$ ni dans $\bm{x}$, la relation $\bm{y}\in\bm{G}\langle\{\bm{x}\}\rangle$ est équivalente à $(\bm{x},\bm{y})\in\bm{G}$. Si $\bm{G}$ et $\bm{G}'$ sont deux graphes, la relation $\bm{G}\subset\bm{G}'$ est donc équivalente à $(\forall\bm{x})(\bm{G}\langle\{\bm{x}\}\rangle\subset\bm{G}'\langle\{\bm{x}\}\rangle)$ où $\bm{x}$ est une lettre ne figurant ni dans $\bm{G}$ ni dans $\bm{G}'$. En effet, cette dernière est identique à $(\forall\bm{x})(\forall\bm{y})(\bm{y}\in\bm{G}\langle\{\bm{x}\}\rangle\Rightarrow\bm{y}\in\bm{G}'\langle\{\bm{x}\}\rangle)$ où $\bm{y}$ est une lettre distincte de $\bm{x}$ et ne figurant ni dans $\bm{G}$ ni dans $\bm{G}'$, donc équivalente à $(\forall\bm{x})(\forall\bm{y})((\bm{x},\bm{y})\in\bm{G}\Rightarrow(\bm{x},\bm{y})\in\bm{G}')$ (cf.
 la remarque \ref{=2graphes}).

Si $\Gamma=(\bm{G},\bm{A},\bm{B})$ est une correspondance entre $\bm{A}$ et $\bm{B}$, pour tout $\bm{x}\in\bm{A}$ la coupe de $\bm{G}$ suivant $\bm{x}$ s'appelle encore la \emph{coupe} de $\Gamma$ suivant $\bm{x}$ et se note également $\Gamma\langle\{\bm{x}\}\rangle$ (ou $\Gamma(\bm{x})$).

\subsection*{4. Fonctions}

D'abord on donne deux remarques:

\begin{remarks}\label{remarques fonct}
\begin{enumerate}[1)]
  \item Soient $\bm{R}$ une relation et $\bm{x}$, $\bm{x}'$ des lettres tels que $\bm{x}'$ ne figure pas dans $\bm{R}$. La relation ``$\bm{R}$ est univoque (resp. fonctionnelle) en $\bm{x}$'' est identique à ``$(\bm{x}'\mid\bm{x})\bm{R}$ est univoque (resp. fonctionnelle) en $\bm{x}'$''. La première partie d'après CS1 (en prenant $\bm{y}$ et $\bm{z}$ deux lettres distinctes, distinctes de $\bm{x}$, de $\bm{x}'$ et ne figurant pas dans $\bm{R}$) et la seconde d'après CS8.
  \item Soient $\bm{R}$ une relation, $\bm{B}$ un terme, et $\bm{x}$, $\bm{t}$ deux lettres distinctes tels que $\bm{x}$ ne figurant pas dans $\bm{B}$.

      D'après CS9, CS5, CS6, CS2, la relation (cf. CS8) $$(\bm{B}\mid\bm{t})(\text{il existe au plus un } \bm{x} \text{ tel que }\bm{R})$$ est identique à $$\text{il existe au plus un } \bm{x} \text{ tel que } (\bm{B}\mid\bm{t})\bm{R}$$
      (en prenant $\bm{y}$ et $\bm{z}$ deux lettres distinctes, distinctes de $\bm{x}$, de $\bm{t}$ et ne figurant ni $\bm{R}$ ni dans $\bm{B}$);

la relation $$(\bm{B}\mid\bm{t})(\text{il existe un } \bm{x} \text{ et un seul tel que }\bm{R})$$ est identique à $$\text{il existe un } \bm{x} \text{ et un seul tel que }(\bm{B}\mid\bm{t})\bm{R}.$$
\end{enumerate}
\end{remarks}

\begin{definition}
On dit qu'un graphe $\bm{F}$ est un graphe \emph{fonctionnel} si, pour tout\footnote{On prend $\bm{x}$ et $\bm{y}$ deux lettres distinctes et ne figurant pas dans $\bm{F}$} $\bm{x}$, il existe au plus un objet correspondant à $\bm{x}$ par $\bm{F}$, plus présisement, si la relation
\begin{equation}\label{graphe fonct}
(\forall\bm{x})(\text{il existe au plus un } \bm{y} \text{ tel que } (\bm{x},\bm{y})\in\bm{F})
\end{equation}
est vraie (cf. I, p. 40). La relation \eqref{graphe fonct} est indépendante du choix de $\bm{x}$ et de $\bm{y}$. On le voit évidemment lorsque on considère l'assemblage de cette relation.

Soit $\bm{f}=(\bm{F},\bm{A},\bm{B})$ une correspondance. Dans ce cas la relation \eqref{graphe fonct} est équivalente d'après I, p. 48, exerc. 7, à la relation obtenue en rempla\c{c}ant $\forall$ par $\forall_{(\bm{x}\in\bm{A})}$\footnote{On prend $\bm{x}$ et $\bm{y}$ deux lettres distinctes et ne figurant ni dans $\bm{F}$ ni dans $\bm{A}$} dans \eqref{graphe fonct}.

On dit que $\bm{f}$ est une \emph{fonction} ou une \emph{application} si son graphe $\bm{F}$ est un graphe fonctionnel et si $\bm{A}=\pr_1\bm{F}$; autrement dit, si pour tout $\bm{x}\in\bm{A}$, la relation $(\bm{x},\bm{y})\in\bm{F}$ est fonctionnelle en $\bm{y}$ (d'après C40, C35).

Soit $\bm{X}\in\bm{A}$. On prend $\bm{y}$ ne figurant pas aussi dans $\bm{X}$. D'après la remarque \ref{remarques fonct}, (2), la relation $(\bm{X},\bm{y})\in\bm{F}$ est fonctionnelle en $\bm{y}$. Son symbole fonctionnel $\bm{\tau}_{\bm{y}}((\bm{X},\bm{y})\in\bm{F})$ s'appelle la \emph{valeur} de $\bm{f}$ pour l'élément $\bm{X}$ de $\bm{A}$, et on le désigne par $\bm{f}(\bm{X})$ ou $\bm{f}_{\bm{X}}$ (ou $\bm{F}(\bm{X})$, ou $\bm{F}_{\bm{X}}$). La relation $\bm{y}=\bm{f}(\bm{X})$ est donc équivalente à $(\bm{X},\bm{y})\in\bm{F}$ (I, p. 41, C46). Donc si $\bm{Y}$ est un ensemble, la relation $\bm{Y}=\bm{f}(\bm{X})$ est équivalente à $(\bm{X},\bm{Y})\in\bm{F}$ (I, C27, C30).

Soient $\bm{x}$ une lettre ne figurant pas dans $\bm{F}$ et $\bm{X}\in\bm{A}$. On désigne $\bm{\tau}_{\bm{y}}((\bm{x},\bm{y})\in\bm{F})$ par $\bm{f}(\bm{x})$ où $\bm{y}$ une lettre distincte de $\bm{x}$ et ne figurant pas dans $\bm{F}$. La relation $(\bm{X}\mid\bm{x})\bm{f}(\bm{x})$ est donc identique à $(\bm{X}\mid\bm{x})\bm{\tau}_{\bm{y}'}((\bm{X},\bm{y}')\in\bm{F})$ où $\bm{y}'$ est une lettre distincte de $\bm{x}$ et $\bm{y}$ et ne figurant pas dans $\bm{X}$ et $\bm{F}$ (CS3). Donc identique à $\bm{\tau}_{\bm{y}''}((\bm{X},\bm{y}'')\in\bm{F})$ où $\bm{y}''$ est une lettre ne figurant pas dans $\bm{X}$ et $\bm{F}$ (CS4); c'est à dire à $\bm{f}(\bm{X})$.
\end{definition}

\textbf{Dans la suite et par commodité on cessera d'utiliser des lettres italiques grasses pour désigner des assemblages indéterminés. Dans un assemblage les lettres considérées sont sous-entendus deux à deux distinctes et ne figurant pas dans les ensembles qui y figurent.}

\begin{remark}\label{appl}
Soit $F$ un graphe fonctionnel. Soient $x,y\in \pr_1F$ et $z$, $t$ deux ensembles, tels que $(x,z)\in F$, $(y,t)\in F$. On a $x=y$ entraîne $z=t$. En effet, $(x,t)=(y,t)$ (II, p. 7, Prop. 1). D'après S6 (en considérant la relation ``$u\in F$'' où $u$ est une lettre ne figurant pas dans $F$), $(x,t)\in F$. Donc $z=t$.

Soit $f=(F,A,B)$ une application. Pour tout $x,y\in A$, $x=y$ entraîne $f(x)=f(y)$. En effet, $(x,f(y))=(y,f(y))\in F$. Donc $f(x)=f(y)$.
\end{remark}

\subsection*{7. Composée de deux fonctions. Fonction réciproque}

Soient $f=(F,A,B)$ et $g=(G,B,C)$ deux applications. Soient $x\in A$ et $z\in C$ tels que $(x,z)\in G\circ F$. Donc il existe un ensemble $y$ tel que $(x,y)\in F$ et $(y,z)\in G)$. Ce qui entraîne $y=f(x)\in B$ et $z=g(y)$.  D'après la remarque \ref{appl}, $z=g(f(x))$. Finalement $(g\circ f)(x)=g(f(x))$. D'où l'application $g\circ f$ est égale à l'application $x\mapsto g(f(x))$.

\vspace{0.5cm}
Le résultat suivant est une généralisation de la prop. 7 de II, p. 17.
\begin{proposition}\label{fonct recip}
Soient $f=(F,A,B)$ une application et $f^{-1}=(F^{-1},B,A)$ sa correspondance réciproque.
\begin{enumerate}[(1)]
  \item $B=\pr_1F^{-1} \Leftrightarrow f \text{ est surjective}$.
  \item $ F^{-1} \text{ est fonctionnel} \Leftrightarrow f \text{ est injective}$.
  \item Pour que $f^{-1}$ soit une fonction, il faut et il suffit que $f$ soit bijective.
\end{enumerate}
\end{proposition}
(1) On a $\pr_1F^{-1}=\pr_2F=F\langle\pr_1F\rangle=F\langle A\rangle=f\langle A\rangle$. L'avant dernière égalité est immédiate du fait que deux ensembles sont égaux si et seulement si chacun d'entre eux est contenu dans l'autre (d'après S6 et l'axiome d'extensionalité) et en appliquant (II, p. 10, Prop. 2). D'où (1).

(2) ``$\Rightarrow$''. Soient $x$ et $y$ deux éléments de $A$ tels que $f(x)=f(y)$. On a $(f(x),x)\in F^{-1}$ et $(f(y),y)\in F^{-1}$. D'après la remarque \ref{appl}, $x=y$.

``$\Leftarrow$''. Soient $x$ et $y$ deux ensembles tels que $(x,y)\in F^{-1}$ et $(x,z)\in F^{-1}$. Donc $x=f(y)=f(z)$. Donc $y=z$.

(3) est immédiate en appliquant (1) et (2).

\vspace{0.5cm}
La remarque suivante est la remarque de II, p. 18, mais avec des démostrations plus élégantes.

\begin{remark}
Soit $f$ une application de $A$ dans $B$.
\begin{enumerate}
  \item Pour toute partie $X$ de $A$, on a vu (II, p. 12) que l'on a $X\subset f^{-1}\langle f\langle X\rangle\rangle$.
  \item Pour toute partie $Y$ de $B$, on a, d'après II, p. 50, exerc. 7, $f\langle f^{-1}\langle Y\rangle\rangle\subset Y$.
  \item Si $f$ est une surjection, on a $f\langle f^{-1}\langle Y\rangle\rangle=Y$ pour toute partie $Y$ de $B$. Cela est immédiat du fait (II, p. 12): si $G$ est un graphe et $X\subset\pr_1G$, $G\subset G^{-1}\langle G\langle G\rangle\rangle$. ($\pr_1F^{-1}=B$.)
  \item Si $f$ est une injection, d'après II, p. 50, exerc. 7, pour toute partie $X$ de $A$, on a $f^{-1}\langle f\langle X\rangle\rangle=X$.
\end{enumerate}
\end{remark}

\subsection*{8. Rétractions et sections}

On donne les détails de la démonstration de la prop. 8.

Supposons que $f$ est surjective. Soient $x$ une lettre ne figurant ni dans $A$, ni dans $F$, ni dans $B$, et $y$ une lettre distincte de $x$ et ne figurant ni dans $A$ ni dans $F$. Désignons par $T$ le terme $\bm{\tau}_y(y\in A \text{ et } f(y)=x)$. Soit $X\in B$, on a $(X\mid x)T$ est identique \'{a} $\bm{\tau}_z(z\in A \text{ et } f(z)=X)$ pour $z$ une lettre (distincte de x, de y et) ne figurant ni dans $A$ ni dans $F$ ni dans $X$ (CS3, CS4). On le désigne par $s(X)$. La relation $(\exists z)(z\in A \text{ et } f(z)=X)$ est vraie et n'est autre que $(s(X)\mid z)(z\in A \text{ et } f(z)=X)$  qui est identique \'{a} ``$s(X)\in A \text{ et }f(s(X))=X$''. Considérons donc l'application $s:x\mapsto T\,(x\in B, T\in A)$. Finalement $f\circ s=\text{Id}_B$.

Supposons que $f$ est injective et $A\neq\varnothing$. Soient $a\in A$ et $x$, $y$ deux lettres distinctes et ne figurant pas dans chacun des ensembles $A$, $B$, $F$, $a$.

La relation $(y\in A \text{ et } x=f(y))\text{ ou } (y=a \text{ et }x\in B-f(A))$ entraîne $(x,y)\in B\times A$ d'après (11), p. 3. Donc elle admet un graphe disons $R$.

D'après (16), (18), pp. 3, 4, si $R$, $S$, $R'$ et $S'$ des relations, les relations $$(R\text{ ou }S)\text{ et }(R'\text{ ou }S'),$$ $$(R\text{ et }R')\text{ ou }(S\text{ et }R')\text{ ou }(R\text{ et }S')\text{ ou }(S\text{ et }S')$$ sont équivalentes.

Soient $x$ et $y$ deux ensembles tels que $(x,y)\in R$ et $(x,y')\in R$. On a
\begin{enumerate}
  \item $\Big((y\in A \text{ et } x=f(y))\text{ et }(y'\in A \text{ et } x=f(y'))\Big)\Rightarrow (y=y')$ (car $f$ est injective).
  \item $\Big(\underbrace{(y=a \text{ et } x\in B-f(A))\text{ et }(y'\in A \text{ et } x=f(y')}_{(*)})\Big)\Rightarrow (\underbrace{x\in f(A)\text{ et } x\in B-f(A)}_{(**)})$.
  \item $\Big((y\in A \text{ et } x=f(y))\text{ et }(y'=a \text{ et } x\in B-f(A))\Big)\Rightarrow (x\in f(A)\text{ et } x\in B-f(A))$.
  \item $\Big((y=a \text{ et } x\in B-f(A))\text{ et }(y'=a \text{ et } x\in B-f(A))\Big)\Rightarrow (y=y')$.
  \end{enumerate}
La négation de la relation $(**)$ est vraie, donc celle de $(*)$ est aussi vraie. Même chose pour 3. D'après (5), p. 3, $y=y'$ est vraie. Donc $R$ est fonctionnel. On a évidemment $\pr_1R\subset B$. Inversement, soit $x\in B$. Si $x\in f(A)$, il existe $y\in A$ tel que $x=f(y)$. Donc $(x,y)\in R$. Si $x\in B-f(A)$, $(x,a)\in R$. D'où $\pr_1R=B$. Soit donc l'application $r=(R,B,A)$. On a $r(x)=a$ si $x\in B-f(A)$ et $f(r(x))=x$ si $x\in f(A)$. D'où $r\circ f=\text{Id}_A$.

\vspace{0.5cm}
Finalement, Dans la démonstration du th. 1 de II, p. 19, on se sert de la remarque suivante.

\begin{remark}
Soit $B$ un ensemble. La seule application $\varnothing\to B$ est $(\emptyset,\emptyset,B)$ (car un graphe est vide si son ensemble de définition est vide, cf. p. 6). Elle injective. En effet, si $R$ une relation et $x$ une lettre, on a ``$x\not\in\emptyset$'' est un théorème et d'après C35, S2, C27, la relation $(\forall_{(x\in\emptyset)} x)R$ est vraie.
\end{remark}

\subsection*{9. Fonctions de deux arguments}

Soient $D$ un graphe et $f:D\to C$ une application. Soient $y\in \pr_1D$ et $A_y$ l'ensemble des $x$ tels que $(x,y)\in D$. L'application $f(.,y)$ est la composition de la bijection canonique $A_y\to A_y\times \{y\}$ (cf. II, p. 17, Exemple 5) et la restriction de $f$ à $A_y\times \{y\}$.

\vspace{0.5cm}
Dans la définition du produit de deux applications nous avons besoin de:

Soient $f:A\to B$ et $g:A\to C$ deux applications. L'application $h:A\to B\times C$, $x\mapsto (f(x),g(x))$ est l'unique application satifaisant $\pr_1\circ h=f$ et $\pr_2\circ h=g$. En effet, si on prend $x$ et $y$ deux lettres distinctes et ne figurant ni dans $A$ ni dans les graphes de $f$ et de $g$. Soit $X\in A$, $(X\mid x)(f(x),g(x))$ est identique à $((X\mid x)f(x),(X\mid x)g(x))$ qui est $(f(X),g(X))$. Ce dernier est un élément de $B\times C$.

\section*{\S4. RÉUNION ET INTERSECTION D'UNE FAMILLE D'ENSEMBLES}

\subsection*{2. Propriétés de la réunion et de l'intersection}

On verra que l'introduction de la sous-section: 1. L'axiome de l'ensemble des parties (II, pp. 30, 31) est obligatoire avant les énoncés de la prop. 2 de II, p. 24, et ceux des résultas de la sous-section 3.

On reprend l'énoncé de la Proposition 2. On considère la composée des applications $L\to \mathfrak{P}(I),\,\lambda\mapsto J_\lambda$ et $\mathfrak{P}(I)\to \mathfrak{P}(\bigcup\limits_{i\in I}X_i),\,J\mapsto\bigcup\limits_{i\in J}X_i$ (cf. II, p. 15, C54). Donc $(\bigcup\limits_{i\in J_\lambda}X_i)_{\lambda\in L}$ est une famille de parties de $\bigcup\limits_{i\in I}X_i$. De même en considérant l'application $L\to\mathfrak{P}(\bigcup\limits_{i\in I}X_i),\,\lambda\mapsto\bigcap\limits_{i\in J_\lambda}X_i$, $(\bigcap\limits_{i\in J_\lambda}X_i)_{\lambda\in L}$ est donc une famille de parties de $\bigcup\limits_{i\in I}X_i$. Ce qui donne un sens à la prop. 2 (II, p. 24).

\subsection*{3. Images d'une réunion et d'une intersection}

Soient $\Gamma$ une correspondance entre $A$ et $B$, $(X_i)_{i\in I}$ une famille de parties de $A$, et $(Y_i)_{i\in I}$ une famille de parties de $B$. Soient $X:I\to\mathfrak{B}$ et $Y:I\to\mathfrak{B}'$ les applications définissant ces deux familles respectivement. L'application composée $I\xrightarrow X \mathfrak{B}\hookrightarrow \mathfrak{P}(A)\xrightarrow {\hat\Gamma}\mathfrak{P}(B)$ définit la famille $(\Gamma(X_i))_{i\in I}$ de parties de $B$ (l'application $\hat\Gamma$ est définit dans II, p. 30). Noter que les lettres figurant dans le graphe de l'injection canonique $\mathfrak{B}\hookrightarrow \mathfrak{P}(A)$ sont celles figurant dans $\mathfrak{B}$; et les lettres figurant dans le graphe de $\hat\Gamma$ sont celles figurant dans le graphe de $\Gamma$ ou dans $A$.

De même l'application composée $I\xrightarrow Y \mathfrak{B}'\hookrightarrow \mathfrak{P}(B)\xrightarrow {\hat{\overset{-1}\Gamma}}\mathfrak{P}(A)$ définit la famille $(\Gamma^{-1}(Y_i))_{i\in I}$ de parties de $A$.

\begin{proposition}
\begin{enumerate}[(1)]
Soient $\Gamma$ une correspondance entre $A$ et $B$, $(X_i)_{i\in I}$ une famille de parties de $A$, et $(Y_i)_{i\in I}$ une famille de parties de $B$. On a donc
  \item $\Gamma\langle\bigcup\limits_{i\in I}X_i\rangle=\bigcup\limits_{i\in I}\Gamma\langle X_i\rangle$.
  \item $\Gamma\langle\bigcap\limits_{i\in I}X_i\rangle\subset\bigcap\limits_{i\in I}\Gamma\langle X_i\rangle$.
  \item Si le graphe de $\Gamma$ est fonctionnel et $I\not=\varnothing$, $\Gamma^{-1}\langle\bigcap\limits_{i\in I}Y_i\rangle=\bigcap\limits_{i\in I}\Gamma^{-1}\langle Y_i\rangle$. Si de plus $\Gamma^{-1}\langle B\rangle=A$, cette égalité est encore vraie dans le cas où $I=\varnothing$.
\end{enumerate}
\end{proposition}
Les assertions (1) et (2) sont déjà démontrées (cf. la prop. 3 de II, p. 25). La première inclusion de (3) s'obtient en appliquant (2) et l'autre s'obtient facilement.

\begin{corollary}
\begin{enumerate}[(a)]
  \item Soient $f$ une application de $A$ dans $B$, et $(Y_i)_{i\in I}$ une famille de parties de $B$. On a $f^{-1}\langle\bigcap\limits_{i\in I}Y_i\rangle=\bigcap\limits_{i\in I}f^{-1}\langle Y_i\rangle$.
  \item Si $f$ une injection de $A$ dans $B$ et si $(X_i)_{i\in I}$ une famille de parties de $A$ avec $I\not=\varnothing$, on a $f\langle\bigcap\limits_{i\in I}X_i\rangle=\bigcap\limits_{i\in I}f\langle X_i\rangle$.
\end{enumerate}
\end{corollary}
S'obtient immédiatement en appliquant (3) et la prop. \ref{fonct recip}, (2).

\subsection*{4. Complémentaire d'une réunion ou d'une intersection}

Soit $E$ un ensemble. On définit l'application $\mathfrak{P}(E)\to\mathfrak{P}(E)$, $X\mapsto \complement_E X$ (cf. II, p. 15, C54). On utilise pour cela le fait que si $R$ est une relation collectivisante en une certaine lettre $x$, et si $y$ est une lettre et $U$ est un ensemble tels que $x$ est distincte de $y$ et ne figurant pas dans $U$, on a par l'axiome d'extensionalité $(U\mid y)\{x\mid R\}=\{x\mid (U\mid y)R\}$.

\begin{proposition}
Pour toute famille $(X_i)_{i\in I}$ de parties de $E$, on a
$$\complement_E(\bigcup\limits_{i\in I}X_i)=\bigcap\limits_{i\in I}(\complement_E X_i)\quad\text{et}\quad\complement_E(\bigcap\limits_{i\in I}X_i)=\bigcup\limits_{i\in I}(\complement_E X_i).$$
\end{proposition}
D'après la seconde partie de C38 et C35, pour que $x\notin \bigcup\limits_{i\in I}X_i$ il faut et il suffit que pour tout $i\in I$, $x\notin X_i$. Donc la première égalité est vraie pour $I\neq\varnothing$. Cette égalité est encore vraie dans le cas où $I=\varnothing$: $\complement_E(\bigcup\limits_{i\in \varnothing}X_i)=\complement_E\varnothing=E=\bigcap\limits_{i\in \varnothing}(\complement_E X_i)$. La seconde résulte immédiatement, en vertu de la relation $\complement_E (\complement_E X)=X$ pour toute partie $X$ de $E$.

\subsection*{5. Réunion et intersection de deux ensembles}

Soit $(X_i)_{i\in I}$ une famille de parties d'un ensemble $E$ avec $I=\bigcup\limits_{\lambda\in L}J_\lambda$. On va vérifier que dans ce cas la seconde partie de la prop. 2 de II, p. 24, est encore vraie sans les restrictions sur $L$ et $J_\lambda$. On considère la composée des applications $L\to \mathfrak{P}(I),\,\lambda\mapsto J_\lambda$ et $\mathfrak{P}(I)\to \mathfrak{P}(E),\,J\mapsto\bigcap\limits_{i\in J}X_i$ (cf. II, p. 15, C54). Donc $(\bigcap\limits_{i\in J_\lambda}X_i)_{\lambda\in L}$ est une famille de parties de $E$. Soit $L'=\{\lambda\mid\lambda\in L \text{ et }J_\lambda=\varnothing\}\subset L$. On distingue les cas suivants:
\begin{itemize}
  \item $L=\varnothing$. Dans ce cas $I=\varnothing$. L'égalité est trivial: $E=E$.
  \item $L'=L$. Dans ce cas aussi l'égalité est trivial: $E=E$.
  \item $\varnothing\neq L'\subsetneqq L$. On a $I=(\underbrace{\bigcup\limits_{\lambda\in L'}J_\lambda}_{=\varnothing})\bigcup(\bigcup\limits_{\lambda\in L-L'}J_\lambda)=\bigcup\limits_{\lambda\in L-L'}J_\lambda.$
      Donc $\bigcap\limits_{\lambda\in L}(\bigcap\limits_{i\in J_\lambda}X_i)=(\underbrace{\bigcap\limits_{\lambda\in L'}(\underbrace{\bigcap\limits_{i\in J_\lambda}X_i}_{=E})}_{=E})\bigcap(\bigcap\limits_{\lambda\in L-L'}(\bigcap\limits_{i\in J_\lambda}X_i))=\bigcap\limits_{\lambda\in L-L'}(\bigcap\limits_{i\in J_\lambda}X_i)=\bigcap\limits_{i\in I}X_i.$
\end{itemize}

\begin{proposition}
Soit $\Gamma$ une correspondance entre $A$ et $B$ dont le graphe est fonctionnel; pour toute partie $Y$ de $B$, on a $\Gamma^{-1}\langle\complement_B Y\rangle=\complement_{\Gamma^{-1}\langle B\rangle} \Gamma^{-1}\langle Y\rangle$.
\end{proposition}

En effet, soit $x\in\complement_{\Gamma^{-1}\langle B\rangle} \Gamma^{-1}\langle Y\rangle$. Donc il existe $x'\in B$ tel que $(x',x)\in\Gamma^{-1}$. On a $x'\notin Y$ car $x\notin\Gamma^{-1}\langle Y\rangle$. Ce qui donne $\Gamma^{-1}\langle\complement_B Y\rangle$. Inversement, soit $x\in\Gamma^{-1}\langle\complement_B Y\rangle$. Donc il existe $x'\in \complement_B Y$ tel que $(x',x)\in\Gamma^{-1}$. On a $x\in\Gamma^{-1}\langle B\rangle$. On a aussi $x\notin\Gamma^{-1}\langle Y\rangle$ car sinon il existerait $x''\in Y$ tel que $(x'',x)\in \Gamma^{-1}$. Ce qui entraînerait $x'=x''\in Y$; ce qui donne une contradiction. Donc $x\in\complement_{\Gamma^{-1}\langle B\rangle} \Gamma^{-1}\langle Y\rangle$.

\begin{corollary}
\begin{enumerate}
  \item Soit $f$ une application de $A$ dans $B$; pour toute partie $Y$ de $B$, on a $f^{-1}\langle\complement_B Y\rangle=\complement_{f^{-1}\langle B\rangle}f^{-1}\langle Y\rangle$.
  \item Soit $f$ une injection de $A$ dans $B$; pour toute partie $X$ de $A$, on a $f\langle\complement_A X\rangle=\complement_{f\langle A\rangle}f\langle X\rangle$.
\end{enumerate}
\end{corollary}
S'obtient immédiatement en appliquant la proposition précédente et la prop. \ref{fonct recip}, (2).

\subsection*{6. Recouvrements}

Soient $(X_i)_{i\in I}$ et $(Y_k)_{k\in K}$ deux familles d'ensembles. Alors il existe un ensemble $E$ (tout ensemble contenant tous les $X_i$ et tous les $Y_k$, c'est-\'{a}-dire, contenant $(\bigcup\limits_{i\in I}X_i)\cup(\bigcup\limits_{k\in K}Y_k)$), tel que $X$  est aussi une application de $I$ dans $\mathfrak{P}(E)$ et $Y$ une application de $K$ dans $\mathfrak{P}(E)$. D'après C51 on définit une application $I\times K\to \mathfrak{P}(E\times E)$, $(i,k)\mapsto X_i\times Y_k$, et une autre $I\times K\to \mathfrak{P}(E)$, $(i,k)\mapsto X_i\cap Y_k$. Cette dernière application est encore la composée de $I\times K\to \mathfrak{P}(E)\times\mathfrak{P}(E)$, $(i,k)\mapsto (X_i,Y_k)$ (produit de $X$ et $Y$, cf. II, p. 21), et $\mathfrak{P}(E)\times\mathfrak{P}(E)\to\mathfrak{P}(E)$, $(X,Y)\mapsto X\cap Y$ (définie d'après C51). Donc on a deux familles d'ensembles $(X_i\times Y_k)_{(i,k)\in I\times K}$ et $(X_i\cap Y_k)_{(i,k)\in I\times K}$.

\subsection*{8. Somme d'une famille d'ensembles}

On commence par un lemme dont la vérification est facile.
\begin{lemma}\label{II, p. 29, prop. 8'}
Soient $(X_i)_{i\in I}$ et $(Y_i)_{i\in I}$ deux familles d'ensembles mutuellement disjoints indexées par un même ensemble $I$, et $f_i:X_i\to Y_i$ une application pour tout $i\in I$. Si les applications $f_i$ sont injectives (resp. surjective, resp. bijective) alors l'application canonique $f: \bigcup\limits_{i\in I}X_i\to \bigcup\limits_{i\in I}Y_i$ (cf. II, p. 29, Prop. 8)  est injective (resp. surjective, resp. bijective).
\end{lemma}

Le résultat suivant est une amélioration de la prop. 9 de II, p. 29.

\begin{proposition}\label{som-fam-ens}
Soit $(X_i)_{i\in I}$ une famille d'ensembles. Il existe un ensemble $X$ unique à une bijection près possédant la propriété: $X$ est réunion d'une famille $(X'_i)_{i\in I}$ d'ensembles mutuellement disjoints, telle que, pour tout $i\in I$, il existe une application bijective de $X_i$ sur $X'_i$.
\end{proposition}

Évidemment la famille d'ensembles $(X_i\times\{i\})_{i\in I}$ montre l'existence. L'unicité découle du lemme précédent.

\begin{definition}
Soit $(X_i)_{i\in I}$ une famille d'ensembles. La somme de cette famille d'ensembles est $\bigcup\limits_{i\in I}X_i\times\{i\}$, notée $\bigsqcup\limits_{i\in I}X_i$.
\end{definition}

\begin{corollary} [II, p. 29, Prop. 10]\label{som-fam-ens-prop10}
Soit $(X_i)_{i\in I}$ une famille d'ensembles mutuellement disjoints. Il existe une bijection de $\bigcup\limits_{i\in I}X_i$ sur $\bigsqcup\limits_{i\in I}X_i$.
\end{corollary}

Comme conséquence directe du lemme précédent on a:

\begin{proposition}\label{som-fam-ens-bij}
Soit $(X_i)_{i\in I}$ et $(Y_i)_{i\in I}$ deux familles d'ensembles indexées par un même ensemble $I$ telle qu'il existe une injection (resp. une surjection, resp. une bijection) de $X_i$ dans (resp. sur) $Y_i$, pour tout $i\in I$. Il existe une injection (resp. une surjection, resp. une bijection) de $\bigsqcup\limits_{i\in I}X_i$ dans (resp. sur) $\bigsqcup\limits_{i\in I}Y_i$.
\end{proposition}

\section*{\S5. PRODUIT D'UNE FAMILLE D'ENSEMBLES}

\subsection*{3. Définition du produit d'une famille d'ensembles}

\begin{definition}
Soit $(X_i)_{i\in I}$ une famille d'ensembles. L'ensemble (C51) des $F\in \mathfrak{P}(I\times\bigcup\limits_{i\in I}X_i)$ tels que $F$ est (un graphe) fonctionnel et pour tout $i\in I$, on ait $F(i)\in X_i$, s'appelle le \emph{produit de la famille d'ensembles} $(X_i)_{i\in I}$ et se désigne par $\prod\limits_{i\in I}X_i$. Pour tout $i\in I$, $X_i$ s'appelle le \emph{facteur d'indice} $i$ du produit $\prod\limits_{i\in I}X_i$; l'application
$$F\mapsto F(i):\bm{\tau}_y((i,y)\in F)\quad (F\in\prod\limits_{i\in I}X_i, F(i)\in X_i)$$
s'appelle la \emph{fonction coordonée} (ou \emph{projection}) d'indice $i$, et se note $\pr_i$.
\end{definition}

On dit que $F(i)$ est la \emph{coordonée d'indice} $i$ (ou \emph{projection d'indice $i$}) de $F$; l'image $\pr_i\langle A\rangle$ d'une partie $A$ de $\prod\limits_{i\in I}X_i$ par la fonction coordonée d'indice $i$ s'appelle la \emph{projection d'indice} $i$ de $A$. En considérant l'application $i\mapsto \pr_i\langle A\rangle$ de $I$ dans $\bigcup\limits_{i\in I}X_i$, on a $(\pr_i\langle A\rangle)_{i\in I}$ est une famille d'ensembles. On a en plus $A\subset \prod\limits_{i\in I}\pr_i\langle A\rangle$.

Si, pour tout $i\in I$, $X_i\subset E$ où $E$ est un ensemble, on a $\prod\limits_{i\in I}X_i\subset E^I$. Donc $\prod\limits_{i\in I}X_i$ est en bijection avec son image par l'application canonique $E^I\to\mathcal{F}(I;E)$ (cf. \S5, Sous-section 2).

\bigskip
Le résultat suivant est une amélioration de la prop. 4 de II, p. 33.

\begin{proposition}\label{prod_u}
Soient $(X_i)_{i\in I}$ une famille d'ensembles et $u:K\to I$ une application dont le graphe est $U$. On a une application $F\mapsto F\circ U$ de $\prod\limits_{i\in I}X_i$ dans $\prod\limits_{k\in K}X_{u(k)}$. En plus, cette application est injective (resp. bijective) si $u$ est surjective (resp. bijective).
\end{proposition}

En effet, soit $F\in\prod\limits_{i\in I}X_i$. On considère l'application composée $K\xrightarrow {u} I\xrightarrow {f=(F,I,\pr_2F)}\pr_2F$. Son graphe est $F\circ U$ et pour tout $k\in K$, $F\circ U(k)=f\circ u(k)=F(u(k))\in X_{u(k)}$.

Si $u$ est surjective. Soit $s:I\to K$ une section de $u$ dont le graphe est $S$. La composée de notre application et de l'application $G\mapsto G\circ S$ de $\prod\limits_{k\in K}X_{u(k)}$ dans $\prod\limits_{k\in K}X_{u(s(k))}=\prod\limits_{i\in I}X_i$ est l'application identique de $\prod\limits_{i\in I}X_i$. Donc elle est injective.

Si $u$ est bijective. Soit $s:I\to K$ l'inverse de $u$. L'application définie dans le cas précédent est bien l'inverse de notre application.

\subsection*{4. Produit partiels}

\begin{proposition}\label{prolong-prod-ens}
Soit $(X_i)_{i\in I}$ une famille d'ensembles telle que $X_i\neq\varnothing$ pour tout $i\in I$. étant donnée une application $g$ de $J\subset I$ dans $A=\bigcup\limits_{i\in I}X_i$, telle que $g(i)\in X_i$ pour tout $i\in J$, il existe un prolongement $f$ de $g$ à $I$, tel que $f(i)\in X_i$ pour tout $i\in I$.
\end{proposition}
C'est la prop. 6 de II, p. 34. Dans sa démonstration, on vérifie que le graphe $G\cup (\bigcup\limits_{i\in I-J}\{(i,T_i)\})$ est fonctionnel en s'aidant de C27 (I, p. 32) et S7 (I, p. 38).

\begin{corollary} [Axiome de choix]
Soit $(X_i)_{i\in I}$ une famille d'ensembles. Pour que $\prod\limits_{i\in I}X_i\neq\varnothing$, il faut et il suffit que pour tout $i\in I$, $X_i\neq\varnothing$.
\end{corollary}

En effet, si pour tout $i\in I$, on a $X_i\neq\varnothing$, on a aussi $\prod\limits_{i\in I}X_i\neq\varnothing$ en appliquant la prop. \ref{prolong-prod-ens} au cas où $J=\varnothing$. Inversement, si $\prod\limits_{i\in I}X_i\neq\varnothing$, par définition de $\prod\limits_{i\in I}X_i$, on a, pour tout $i\in I$, $X_i\neq\varnothing$.

\begin{corollary}
Soient $(X_i)_{i\in I}$ et $(Y_i)_{i\in I}$ deux familles d'ensembles ayant le même ensemble d'indices $I$. Si, pour tout $i\in I$, on a $X_i\subset Y_i$, on a aussi $\prod\limits_{i\in I}X_i\subset\prod\limits_{i\in I}Y_i$. Réciproquement, si $\prod\limits_{i\in I}X_i\subset\prod\limits_{i\in I}Y_i$ et si, pour tout $i\in I$, on a $X_i\neq\varnothing$, on a $X_i\subset Y_i$ pour tout $i\in I$.
\end{corollary}

La première partie est évidente. Pour la seconde, soient $A=\bigcup\limits_{i\in I}X_i$, $\alpha\in I$ et $x\in X_\alpha$. Soit $g$ l'application de $\{\alpha\}\subset I$ dans $A$ tel que $g(\alpha)=x$. D'après la prop. \ref{prolong-prod-ens}, on a $\{(\alpha,x)\}\cup(\bigcup\limits_{i\in I-\{\alpha\}}\{(i,T_i)\})\in\prod\limits_{i\in I}X_i$. Donc $x\in Y_\alpha$.

Soient $(X_i)_{i\in I}$ une famille d'ensembles et $J$ une partie de $I$. L'application de $\prod\limits_{i\in I}X_i$ dans $\prod\limits_{i\in J}X_i$ définie dans la prop. \ref{prod_u} qui est associée à l'injection canonique $J\to I$, s'appelle la \emph{projection d'indice} $J$ et se note  $\pr_J$. D'après la prop. \ref{prolong-prod-ens} elle surjective.

\begin{proposition}
Soient $(X_i)_{i\in I}$ une famille d'ensembles et $J$ une partie de $I$. Si, pour tout $i\in I$, on a $X_i\neq\varnothing$, l'application $\operatorname{pr}_J$ est surjective.
\end{proposition}

\begin{corollary}
Soient $(X_i)_{i\in I}$ une famille d'ensembles telle que, pour tout $i\in I$, on ait $X_i\neq\varnothing$. Alors, pour tout $\alpha\in I$, la projection $\operatorname{pr}_\alpha$ est une application surjective de $\prod\limits_{i\in I}X_i$ sur $X_\alpha$.
\end{corollary}

\section*{\S6. RELATIONS D'ÉQUIVALENCE}

\subsection*{1. Définition d'une relation d'équivalence}

Soit $R\{x,y\}$ une relation, $x$ et $y$ étant \underline{deux lettres distinctes figurant dans $R$} (le cas contraire est sans importance et conduit \'{a} plus de vérifications). On dit que $R$ est \emph{symétrique} (par rapport aux lettres $x$ et $y$) si l'on a $R\{x,y\}\Rightarrow R\{y,x\}$; c'est-à-dire, pour tout $X$ et $Y$ deux ensembles, $R\{X,Y\}\Rightarrow R\{Y,X\}$ (cf. I, p. 16). Donc $R\{x,y\}$ et $R\{y,z\}$ sont équivalentes.

Soit $z$ une lettre ne figurant pas dans $R$. On dit que $R$ est \emph{transitive} (par rapport aux lettres $x$ et $y$) si l'on a $R\{x,y\}\text{ et }R\{y,z\} \Rightarrow R\{x,z\}$; c'est-à-dire, pour tout $X$, $Y$ et $Z$ des ensembles, $(R\{X,Y\}\text{ et }R\{Y,Z\}) \Rightarrow R\{X,Z\}$.

\begin{examples}
La relation $x=y$ est symétrique et transitive. La relation $X\subset Y: (\forall x)((x\in X)\Rightarrow (x\in Y))$ est transitive, mais non symétrique. La relation $X\cap Y=\varnothing$ est symétrique, mais non transitive (la relation $X\cap Y$ est identique \'{a} la relation $\bm{\tau}_y(\forall x)((x\in y)\Leftrightarrow (x\in X\text{ ou }x\in Y))$).
\end{examples}

Si $R\{x,y\}$ est \'{a} la fois symétrique et transitive, on dit qu'elle est une \emph{relation d'équivalence} (par rapport aux lettres $x$ et $y$). Dans ce cas

\begin{equation}\label{R{x,x}}
R\{x,y\}\Rightarrow (R\{x,x\}\text{ et }R\{y,y\}).
\end{equation}

Donc pour tout $X$ et $Y$ deux ensembles, $R\{X,Y\}\Rightarrow (R\{X,X\}\text{ et }R\{Y,Y\})$.

Soient $R\{x,y\}$ une relation et $E$ un ensemble tels que $x$ et $y$ deux lettres distinctes figurant dans $R$ mais pas dans $E$. On dit que $R$ est \emph{réflexive dans} $E$ (par rapport aux lettres $x$ et $y$) si $R\{x,x\}\Leftrightarrow (x\in E)$; c'est-à-dire, pour tout $X$ ensemble, $R\{X,X\}\Leftrightarrow (X\in E)$.

Si en plus $R$ est une relation d'équivalence, on dit que $R$ est une \emph{relation d'équivalence dans} $E$.

Soient en plus $x'$ et $y'$ deux lettres distinctes et ne figurant pas dans $R\{\square,\square\}$. D'après C3 (I, p. 23), $R\{x,y\}\Leftrightarrow R\{x',y'\}$. Donc $R\{x,y\}$ est symétrique (resp. transitive, resp. relation d'équivalence) par rapport aux lettres $x$ et $y$ est équivalente \'{a} $R\{x',y'\}$ est symétrique (resp. transitive, resp. relation d'équivalence) par rapport aux lettres $x'$ et $y'$. Si $R$ est une relation d'équivalence dans $E$ par rapport aux lettres $x$ et $y$, et si $x'$ et $y'$ ne figurant pas dans $E$, alors $R\{x',y'\}$ est une relation d'équivalence dans $E$ par rapport aux lettres $x'$ et $y'$.

Soient $\Gamma=(F,E,E)$ une correspondance et $x$ et $y$ deux lettres distinctes ne figurant pas dans $F$. On dit que $\Gamma$ est une \emph{équivalence} si la relation $(x,y)\in F$ est une relation d'équivalence dans $E$.

Soit $R\{x,y\}$ une relation vérifiant \eqref{R{x,x}} (par exemple si $R\{x,y\}$ est une relation d'équivalence ou relation de préordre ou d'ordre, cf. III). Si $R$ est réflexive dans un ensemble $E$, $R\{x,y\}\Rightarrow (x,y)\in E\times E$, donc $R$ admet un graphe (par rapport aux lettres $x$ et $y$). Supposons que $R$ admet un graphe $G$ (par rapport aux lettres $x$ et $y$). On a $R\{x,x\}\Leftrightarrow (\exists y)R\{x,y\}\Leftrightarrow (\exists y)((x,y)\in G)\Leftrightarrow (x\in \pr_1G)$ (pour la première équivalence on utilise S5 (I, p. 33) et \eqref{R{x,x}}, pour la deuxième C31 (I, p. 34), et pour la dernière cf. notre \S3, 1), de fa\c{c}on que $R$ est réflexive dans $\pr_1G$. Donc, pour que $R$ admet un graphe (par rapport aux lettres $x$ et $y$) il faut et il suffit que $R$ est réflexive dans un ensemble $E$. D'où, si $R\{x,y\}$ est une relation d'équivalence (resp. de préordre, resp. d'ordre) dans un ensemble $E$, alors $(G,E,E)$ est une équivalence (resp. un préordre, resp. un ordre) dans $E$ où $G$ est le graphe de $R$. (D'après la réflexivité de $R$ dans $E$ et \eqref{R{x,x}}, $E=\pr_1G=\pr_2G$.)

\subsection*{2. Classes d'équivalence; ensemble quotient}

Soit $f=(F,E,E')$ une application telle que $E\neq\varnothing$. La relation ``$x\in E \text{ et } y\in E \text{ et } f(x)=f(y)$'' est une relation d'équivalence dans $E$; nous la définirons comme la relation d'équivalence \emph{associée} à $f$. Elle est équivalente à la relation $(\exists z)((x,y)\in F\text{ et }(y,z)\in F)$, c'est-à-dire à $(x,y)\in F^{-1}\circ F$. Donc admet $F^{-1}\circ F$ comme graphe.

Nous allons maintenant voir que toute relation d'équivalence $R\{x',y'\}$ dans un ensemble $E$ est de ce type. En effet, soient $G$ le graphe de $R$ et $x\in E$. L'ensemble $G(x)$ (la coupe de $G$ suivant $x$) (ou tout ensemble lui est égal) s'appelle la \emph{classe d'équivalence de $x$ suivant $R$}. Soit $y$ une lettre ne figurant ni dans $x$ ni dans $G$. On a $$G(x)=\{y\mid (x,y)\in G\}=\{y\mid R\{x,y\}\}=\{y\mid (y\in E) \text{ et }R\{x,y\}\}\subset E;$$ et $G(x)\neq\varnothing$ car $x\in G(x)$. Chaque élément de $G(x)$ s'appelle un \emph{représentant} de cette classe. L'ensemble $$E/R=\{y\mid (\exists x)((x\in E)\text{ et }y=G(x))\} \text{ (C53)}$$ s'appelle l'\emph{ensemble quotient} de $E$ par $R$ (ses éléments sont $G(x)$ où $x\in E$). L'application $p:E\to E/R$ définit par $x\mapsto G(x)$ s'appelle  l'\emph{application canonique} de $E$ sur $E/R$ (on prend $x$ une lettre ne figurant ni dans $E$ ni dans $G$, donc, pour tout $X\in E$, $(X\mid x)G(x)$ est identique à $G(X)$).

\vspace{0.5cm}
C55. $R\{x',y'\}\Leftrightarrow (x'\in E \text{ et } y'\in E \text{ et } p(x')=p(y'))\Leftrightarrow (\exists X)(X\in E/R\text{ et }x'\in X\text{ et }y'\in X).$

En effet, $R\{x',y'\}\Rightarrow (R\{x',x'\}\text{ et }R\{y',y'\})\Rightarrow (x'\in E\text{ et }y'\in E)$.
D'autre part,
\begin{equation*}
\begin{split}
R\{x',y'\} &\Leftrightarrow (x',y')\in G
\\ &\Leftrightarrow y'\in G(x')
\\ & \Leftrightarrow G(y')\subset G(x')
\end{split}
\end{equation*}
$(y'\in G(x'))\Rightarrow (G(y')\subset G(x'))$: Soit $z\in G(y')$, donc $(y',z)\in G$. On a $R\{x',y'\}$ et $R\{y',z\}$. Puisque $R$ est transitive, $R\{x',z\}$, i.e. $(x',z)\in G$, c'est-à-dire, $z\in G(x')$.

En vertu de $R\{x',y'\}\Leftrightarrow R\{y',x'\}$, on montre $R\{x',y'\}\Leftrightarrow (G(x')=G(y'))$.

Les lettres $x'$ et $y'$ ne figurent ni dans $E$ ni dans $G$, donc non plus dans le graphe de $p$ (les lettres figurant dans ce dernier sont celles figurant dans $E$ ou $G$).

La seconde équivalence de C55 est évidente, et montre que $R$ est déterminée à une équivalence près par la donnée de la partition $E/R$ de $E$.

\subsection*{3. Relations compatibles avec une relation d'équivalence}

Soient $R\{x,x'\}$ une relation d'équivalence, et $P\{x\}$ une relation. On dit que $P\{x\}$ est \emph{compatible avec la relation d'équivalence} $R\{x,x'\}$ (par rapport à $x$) si, $y$ désignant une lettre qui ne figure ni dans $P$ ni dans $R$, on a
$$(P\{x\}\text{ et }R\{x,y\})\Rightarrow P\{y\};$$
c'est-à-dire, pour tout $a$ et $b$ deux ensembles,  $(P\{a\}\text{ et }R\{a,b\})\Rightarrow P\{b\}.$

Nous donnons notre version de C56 où nous avons supprimés une condition qui n'est pas nécessaire et nous avons ajouté une autre afin d'aboutir le résultat voulu.

\vspace{0.5cm}
C56. Soient $R\{x,x'\}$ une relation d'équivalence dans un ensemble $E$, $P\{x\}$ une relation compatible (par rapport à $x$) avec la relation d'équivalence $R\{x,x'\}$; alors si $t$ ne figure pas dans $P\{x\}$, \underline{ni dans $R$, ni dans $E$ (donc non plus dans $E/R$),} la relation ``$t\in E/R\text{ et }(\exists x)(x\in t \text{ et }P\{x\})$'' est équivalente à la relation ``$t\in E/R\text{ et }(\forall x)((x\in t)\Rightarrow P\{x\})$''.

Nous complétons son démonstration en vérifiant la dernière équivalence. Soit $y\in E$, donc $f(y)\in E/R$. On a
\begin{equation*}
\begin{split}
P'\{f(y)\} &\Leftrightarrow (f(t)\mid t)(\forall x)((x\in t)\Rightarrow P\{x\})
\\ &\Leftrightarrow (f(t)\mid t)(\forall x')((x'\in t)\Rightarrow P\{x'\}) \\ & \qquad\text{où $x'$ est une lettre distincte de $t$ et ne figurant ni dans $P\{x\}$ ni dans $f(y)$}
\\ & \Rightarrow P\{y\}.
\end{split}
\end{equation*}
L'implication $P\{y\}\Rightarrow P'\{f(y)\}$ est évidente.

\subsection*{4. Parties saturées}

Soient $R\{x,y\}$ une relation d'équivalence dans un ensemble $E$, et $A$ une partie de $E$ ne contenant pas $x$. Soit $f$ l'application canonique de $E$ dans $E/R$. On a

\begin{equation}\label{part.sat}
A\subset f^{-1}\langle f\langle A\rangle\rangle=\bigcup\limits_{x\in A}f^{-1}\langle \{f(x)\}\rangle=\bigcup\limits_{x\in A}f(x) \text{ (car $f\langle A\rangle=\bigcup\limits_{x\in A}\{f(x)\}$)}.
\end{equation}
On dit que $A$ est \emph{saturée pour} $R$ si la relation $x\in A$ est compatible (par rapport à $x$) avec $R\{x,y\}$; c'est-à-dire, \emph{pour tout} $x\in A$, \emph{la classe d'équivalence de $x$ est contenu dans $A$}. D'après \eqref{part.sat}, on a
\begin{proposition}
\begin{enumerate}[1)]
Les relations suivantes sont équivalentes
  \item $A$ est saturée pour $R$;
  \item $A=\bigcup\limits_{x\in A}f(x);$
  \item $A$ est réunion d'un ensemble de classes d'équivalence suivant $R$;
  \item $A= f^{-1}\langle f\langle A\rangle\rangle$;
  \item $A= f^{-1}\langle B\rangle$ où $B\subset E/R$.
\end{enumerate}
\end{proposition}

\subsection*{7. Quotients de relations d'équivalence}

Soient $R$ et $S$ deux relations d'équivalence, par rapport à deux lettres $x$ et $y$. Nous dirons que $S$ est \emph{plus fine} que $R$ (ou que $R$ est \emph{moins fine} que $S$) si la relation $S\Rightarrow R$ est vraie.

\begin{proposition}
Soient $R$ et $S$ deux relations d'équivalence dans un même ensemble $E$ dont les graphes respectives sont $G$ et $G'$. Les relations suivantes sont équivalentes:
\begin{enumerate}[1)]
  \item $S\Rightarrow R$;
  \item $G'\subset G$;
  \item $(\forall x\in E)(G'(x)\subset G(x))$;
  \item $(\forall x\in E)(\exists y\in E)(G'(x)\subset G(y))$;
  \item $(\forall x\in E)(\forall y\in G(x))(G'(y)\subset G(x))$.
\end{enumerate}
\end{proposition}

En effet, d'après la remarque \ref{=2graphes}, 1) et 2) sont équivalentes. On a vu (cf. \S3, 3) que 2) et 3) sont équivalentes. L'implication $3)\Rightarrow 4)$ est évidente. L'implication $4)\Rightarrow 3)$ est vraie car $x\in G'(x)$. L'implication $3)\Rightarrow 5)$: Soient $x\in E$ et $y\in G(x)$. On a $G'(y)\subset G(y)=G(x)$. L'implication $5)\Rightarrow 3)$ est vraie car $x\in G(x)$.

\subsection*{8. Produit de deux relations d'équivalence}

Soit $R$ une relation où figurent des lettres deux à deux distinctes $x$, $y$, $x'$ et $y'$. Soient $z$ et $z'$ deux lettres distinctes et ne figurant pas dans $R$. D'après II, p. 49, exerc. 1 de \S2, on a
\begin{equation*}
\begin{split}
(\exists x)(\exists y)(\exists x')(\exists y')R\{x,y,x',y'\} &\Leftrightarrow (\exists x)(\exists y)(\exists z')(z'\text{ est un couple et }R\{x,y,\pr_1z',\pr_2z'\}
\\ &\Leftrightarrow (\exists z)(z\text{ est un couple et }(\exists z')(z'\text{ est un couple et }\\ & \qquad R\{x,y,\pr_1z',\pr_2z'\}\{\pr_1z,\pr_2z\})) \text{ (cf. I, p. 16)}
\\ &\Leftrightarrow (\exists z)(\exists z')(z\text{ est un couple et }z'\text{est un couple et } \\ & \qquad R\{\pr_1z,\pr_2z,\pr_1z',\pr_2z'\}) \text{ (d'après C33 et C31)}
\\ &\Leftrightarrow (\exists t)(t\text{ est un quadruplet et }R\{\pr_1t,\pr_2t,\pr_3t,\pr_4t\})),
\end{split}
\end{equation*}
où $t$ est une lettre (distincte de $z$ et $z'$ et) ne figurant pas dans $R\{\square,\square,\square,\square\}$.

Maintenant, soient $R\{x,y\}$ et $R'\{x',y'\}$ deux relations d'équivalence où $x$, $y$, $x'$ et $y'$ sont des lettres deux à deux distinctes telles que $x'$ et $y'$ ne figurent pas dans $R\{x,y\}$ et $x$ et $y$ ne figurent pas dans $R'\{x',y'\}$. Soient $u$ et $v$ deux lettres distinctes et ne figurant ni dans $R$ ni dans $R'$. Désignons par $S\{u,v\}$ la relation
$$(\exists x)(\exists y)(\exists x')(\exists y')(u=(x,x')\text{ et }v=(y,y')\text{ et }R\{x,y\}\text{ et }R'\{x',y'\});$$ qui est équivalent à
$$(\exists t)(t\text{ est un quadruplet et } u=(\pr_1t,\pr_3t)\text{ et }v=(\pr_2t,\pr_4t)\text{ et }R\{\pr_1t,\pr_2t\}\text{ et }R'\{\pr_3t,\pr_4t\});$$
où $t$ est une lettre distincte de $u$ et $v$ et ne figurant pas dans $R\{\square,\square\}$ ni dans $R'\{\square,\square\}$. Donc si $U$ et $V$ sont deux ensembles, $S\{U,V\}$ est équivalent à $$(\exists t)(t\text{ est un quadruplet et } U=(\pr_1t,\pr_3t)\text{ et }V=(\pr_2t,\pr_4t)\text{ et }R\{\pr_1t,\pr_2t\}\text{ et }R'\{\pr_3t,\pr_4t\});$$
où $t$ est une lettre ne figurant pas dans aucun des assemblages $U$, $V$, $R\{\square,\square\}$ et $R'\{\square,\square\}$.

La relation $S\{u,v\}$ est symétrique car $R$ et $R'$ le sont aussi. En se servant du fait que $R$ et $R'$ sont transitive et de la remarque suivante, elle est encore transitive. Donc la relation $S\{u,v\}$ est une relation d'équivalence, que l'on appelle \emph{produit} de $R$ et $R'$ et qu'on désigne par $R\times R'$.

\begin{remark}\label{S6'}
Soient $A\{x,y\}$ une relation où $x$ et $y$ sont deux lettres distinctes, et $B$, $B'$ et $C$ sont des termes. D'après S6, si $B=B'$, on a $A\{B,C\}\Leftrightarrow A\{B',C\}$ (cf. I, p. 16).
\end{remark}

Si $R$ est une relation réflexive dans un ensemble $E$ et $R'$ est une relation réflexive dans un ensemble $E'$, alors $R\times R'$ est réflexive par rapport à $u$ et $v$ qui sont pris ne figurant ni dans $E$ ni dans $E'$ (donc non plus dans $E\times E'$) (on se sert de la remarque \ref{S6'}).

Soient $R$ une relation d'équivalence dans un ensemble $E$ dont le graphe est $G$ et $R'$ une relation d'équivalence dans un ensemble $E'$ dont le graphe est $G'$. Donc $R\times R'$ est aussi une relation d'équivalence. Si $b=(a,a')\in E\times E'$, la classe d'équivalence de $b$ suivant $R\times R'$ est $G(a)\times G'(a')$. Donc \emph{toute classe d'équivalence suivant $R\times R'$ est le produit d'une classe d'équivalence suivant $R$ et d'une classe d'équivalence suivant $R'$, et réciproquement}.

Soient $f$ et $f'$ les applications canoniques de $E$ sur $E/R$ et de $E'$ sur $E'/R'$, et soit $f\times f'$ l'extension canonique de $f$ et $f'$ aux ensembles produits (II, p. 21). On a $(f\times f')^{-1}(f(a),f(a'))=f(a)\times f(a')$ pour $(a,a')\in E\times E'$.  Si on prend de plus $u$ et $v$ ne figurant pas dans les graphes de $R$ et de $R'$ (donc non plus dans les graphes de $f$, $f'$ et $f\times f'$), alors la relation d'équivalence associée à $f\times f'$: `` $u\in E\times E'\text{ et }v\in E\times E'\text{ et } (f\times f')(u)=(f\times f')(v)$ '' est équivalente à $R\times R'$. L'application $f\times f'$ peut donc se mettre sous la forme $h\circ g$, où $g$ est l'application canonique de $E\times E'$ sur $(E\times E')/(R\times R')$ et où $h$ est une application bijective de $(E\times E')/(R\times R')$ sur $(E/R)\times (E'/R')$; cette application et son application réciproque sont dites \emph{canoniques}.

Dans la remarque de II, p. 47, dire que $P$ est compatible avec les relations d'équivalence $R\{x,y\}$ et $R'\{x',y'\}$ revient au même de dire que pour tout $a$, $b$, $a'$ et $b'$ des ensembles, $(P\{a,a'\}\text{ et }R\{a,b\}\text{ et }R'\{a',b'\})\Rightarrow P\{b,b'\}$.

La lettre $u$ est supposée ne figurant pas dans $P\{x,x'\}$; et on a
\begin{equation*}
\begin{split}
Q\{u\} &\Leftrightarrow (\exists z)(z \text{ est un couple et }(u=(\pr_1z,\pr_2z))\text{ et }P\{\pr_1z,\pr_2z\})
\\ &\Leftrightarrow (u \text{ est un couple et }P\{\pr_1u,\pr_2u\}).
\end{split}
\end{equation*}

\subsection*{9. Classes d'objets équivalents}

Soient $x$, $x'$ et $y$ trois lettres distinctes, et $R\{x,y\}$ une relations d'équivalence. D'après C27, $R\{x,y\}$ est équivalente à $(\forall y)(R\{x,y\}\Leftrightarrow R\{x',y\})$, donc, d'après S7, \'{a} $\bm{\tau}_y(R\{x,y\})=\bm{\tau}_y(R\{x',y\})$ c'est-à-dire à $\theta\{x\}=\theta\{x'\}$ où $\theta\{x\}$ est le terme $\bm{\tau}_y(R\{x,y\})$. Soit $X$ un ensemble, $\theta\{X\}$ est identique à $\bm{\tau}_{y'}(R\{X,y'\})$ où $y'$ est une lettre ne figurant ni dans $R$ ni dans $X$. D'après C31 et S5, $R\{x,\theta\{x\}\}$, qui n'est autre que $(\exists y)R\{x,y\}$, est équivalente à $R\{x,x\}$ (si $T$ est un ensemble, $R\{x,T\}$ est identique à $(T\mid y)R\{x,y\}$).

Donc, d'après la remarque \ref{S6'}, $$(R\{x,x\}\text{ et }R\{x',x'\}\text{ et }\theta\{x\}=\theta\{x'\})\Leftrightarrow R\{x,x'\}.$$
On dit que l'ensemble $\theta\{x\}$ est \emph{la classe d'objets équivalents à} $x$ (suivant la relation $R$).

Soient $X$ et $Y$ deux ensembles. Si on pend en plus $x'$ ne figurant ni dans $R$, on obtient
\begin{equation}\label{theta1}
R\{X,Y\}\Rightarrow (\theta\{X\}=\theta\{Y\})
\end{equation} et
\begin{equation}\label{theta2}
(R\{X,X\}\text{ et }R\{Y,Y\}\text{ et }\theta\{X\}=\theta\{Y\})\Leftrightarrow R\{X,Y\}.
\end{equation}
En plus on a
\begin{equation}\label{theta3}
R\{Y,Y\}\Leftrightarrow R\{Y,\theta\{Y\}\}
\end{equation} (d'après C27, C30, CS7).

En combinant \eqref{theta1} et \eqref{theta3}, on a

\begin{equation}\label{theta4}
R\{Y,Y\}\Rightarrow (\theta\{Y\}=\theta\{\theta\{Y\}\}).
\end{equation}

Soit maintenant $T$ un ensemble ne contenant ni $x$ ni $y$, vérifiant la relation
\begin{equation}\label{T}
(\forall y)(R\{y,y\}\Rightarrow (\exists x)(x\in T\text{ et }R\{x,y\})).
\end{equation}
Soit $\Theta$ l'ensemble $\{z\mid (\exists x)(z=\theta\{x\}\text{ et }x\in T)\}$ (C53, pp. 5, 6). D'après la remarque \ref{C53}, les éléments de $\Theta$ sont de la forme $\theta\{X\}$ où $X\in T$. On a $(\exists x)(R\{x,x\}\text{ et }z=\theta\{x\})\Rightarrow (z\in \Theta)$. Donc, d'après C52 (II, p. 5), la relation $(\exists x)(R\{x,x\}\text{ et }z=\theta\{x\})$ est collectivisante en $z$.

On peut supposer que $(x\in T)\Rightarrow R\{x,x\}$ est vraie; il suffit de remplacer $T$ par l'ensemble $T'=\{x\in T\mid R\{x,x\}\}$. Supposons que $R\{Y,Y\}$ soit vraie où $Y$ est un ensemble. Donc il existe $X\in T$ tel que l'on ait $R\{X,Y\}$. Alors $\theta\{Y\}=\theta\{X\}\in\Theta$. On dit que $\Theta$ est \emph{l'ensemble des classes d'objets équivalents} suivant la relation $R$. D'où, d'après \eqref{theta3}; la remarque \ref{S6'}, \eqref{theta1} et \eqref{theta4}, $\theta\{Y\}$ est l'unique élément $z\in\Theta$ tel que l'on ait $R\{Y,z\}$. Pour l'unicité, soient $z,z'\in \Theta$ tels que $R\{Y,z\}$ et $R\{Y,z'\}$. On a $z=\theta\{X\}$ et $z'=\theta\{X'\}$ où $X,X'\in T$. Donc $z=\theta\{X\}=\theta\{\theta\{X\}\}=\theta\{Y\}=\theta\{\theta\{X'\}\}=\theta\{X'\}=z'$.

Sous les mêmes hypothèses, soit $A\{x\}$ un terme ne contenant pas $y$ tel que $R\{x,y\}$ entraîne $A\{x\}=A\{y\}$. Alors $(\exists x)(R\{x,x\}\text{ et } z=A\{x\})$ est aussi collectivisante en $z$, car, d'après \eqref{theta3}, $R\{X,X\}\Leftrightarrow (A\{X\}=A\{\theta\{X\}\})\Rightarrow (A\{X\}\in E)$, où $E$ est l'ensemble des objets de la forme $A\{t\}$ pour $t\in \Theta$. Si $f$ est l'application $t\mapsto A\{t\}\quad (t\in \Theta, A\{t\}\in E)$, on a $R\{x,x\}$ entraîne $A\{x\}=f(A\{x\})$.

En particulier, si $R$ est une relation d'équivalence dans un ensemble $F$, on peut prendre $T=F$ et pour $A\{x\}$ la classe d'équivalence de $x$ suivant $R$. Dans ce cas $E=E/R$ et la fonction $f$ est une bijection (l'injectivité: si $A\{\theta\{X\}\}=A\{\theta\{Y\}\}$ pour $X,Y\in F$, alors $A\{X\}=A\{\theta\{X\}\}=A\{\theta\{Y\}\}$ et $\theta\{X\},\theta\{Y\}\in \Theta$, et par suite $\theta\{X\}=\theta\{Y\}$). Ce qui justifie la terminologie introduite.

\begin{example}
Soit $R\{x,y\}$ la relation d'équivalence ``$x$ et $y$ sont deux espaces vectoriels de même dimension finie sur $\mathbb{C}$'', qui n'admet pas de graphe (car sinon elle serait une relation d'équivalence dans ``l'ensemble'' des espaces vectoriels de dimension finie sur $\mathbb{C}$). Elle vérifie \eqref{T} en prenant pour $T$ l'ensemble des objets de la forme $\mathbb{C}^n$ pour $n\in \N$ ($\mathbb{C}^0=0$).
\end{example}

\chapter{Notes concernant le chapitre III}\label{app3}

\section*{\S1. RELATIONS D'ORDRES. ENSEMBLES ORDONNÉS}
III, p. 2:
\begin{itemize}
\item Exemple 2: On a supposé que $x$ et $y$ ne figurent pas dans $E$. Maintenant, soient $R\{x,y\}$ une relation transitive et antisymétrique (2ème relation de III, p. 1) entre $x$ et $y$, et $E$ un ensemble tel que \begin{equation}\label{rel-ind}
\text{pour tout } X\in E,\; (X\in E)\Rightarrow R\{X,X\}.
\end{equation} Soient $x'$ et $y'$ deux lettres distinctes et ne figurant pas ni dans $R\{\square,\square\}$ ni dans $E$. La relation $R\{x',y'\}$ est équivalente à $R\{x,y\}$ donc est une relation d'ordre. Donc, en appliquant l'exemple 2 à $R\{x',y'\}$, la relation ``$R\{x',y'\}\text{ et } x'\in E\text{ et } y'\in E$'' est une relation d'ordre dans $E$, appelée relation d'ordre \emph{induite} par $R\{x,y\}$ dans $E$. Si en plus $R\{x,y\}$ est une relation d'ordre dans un ensemble $E'$ de graphe $G$. La condition \eqref{rel-ind} est équivalente à $E\subset E'$. La relation d'ordre induite par $R\{x,y\}$ dans $E$ est équivalente à $(x',y')\in G\cap (E\times E)$ (les relations $R\{x,y\}$ et $R\{x',y'\}$ sont équivalentes donc elles sont des relations d'ordres sur $E'$ et ont même graphe $G$). Donc elle n'est autre que celle induite par la relation sur $E$ définie dans III, p. 5.
\item Prop. 1: Je crois qu'il faut ajouter la condition ``$G\subset \Delta\circ G\text{ et }G\subset G\circ \Delta$'' à b) pour avoir l'équivalence.
\end{itemize}

\vspace{0.5cm}
III, p. 3: Plus généralement, soit $R\{x,y\}$ une relation de préordre. Soit $S\{x,y\}$ une relation d'équivalence telle que $R\{x,y\}$ soit compatible avec les relations d'équivalence $S\{x,x'\}$ et $S\{y,y'\}$  (par rapport à $x$ et $y$) où $x'$ et $y'$ sont des lettres ne figurant pas dans $S$, c'est-à-dire, pour tout ensembles $X$, $Y$ et $Z$, $(R\{X,Y\}\text{ et } S\{X,Z\})\Rightarrow R\{Z,Y\}$  et $(R\{X,Y\}\text{ et } S\{Y,Z\})\Rightarrow R\{X,Z\}$.

Soit $R'\{X,Y\}$ la relation $$X\in E/S\text{ et }Y\in E/S\text{ et }(\exists x)(\exists y)(x\in X\text{ et }y\in Y\text{ et }R\{x,y\}).$$
On propose ici une démonstration directe et simple du fait que $R'\{X,Y\}$ est une relation d'ordre entre éléments de E/S sans se servir de C56, II, p. 43. Soient $X$, $Y$ et $Z$ des éléments de $E/S$ tels que $R'\{X,Y\}$ et $R'\{Y,Z\}$. Donc, d'après la compatibilité et II, p. 49, \S2, exerc. 1, il existe $x\in X$, $y\in Y$ et $z\in Z$ tels que $R\{x,y\}$ et $R\{y,z\}$. D'où $E/S$ transitive. Soient $X$, $Y$ deux éléments de $E/S$ tels que $R'\{X,Y\}$ et $R'\{Y,X\}$. Donc, d'après ce qui précède, il existe $x\in X$, $y\in Y$ et $z\in X$ tels que $R\{x,y\}$ et $R\{y,z\}$. D'après la seconde implication de la compatibilité (5ème ligne), on a $R\{y,x\}$. D'où $X=Y$. Le reste est facile à voir.

\vspace{0.5cm}
III, p. 4: Je crois que, pour que $\Gamma$ soit un préordre sur $E$ il faut et il suffit que $\Delta\subset G$, $G\subset \Delta\circ G$, $G\subset G\circ \Delta$ et $G\circ G\subset G$.

\vspace{0.5cm}
III, p. 5: Je crois que C58 n'est pas clair. À mon avis C58 est plus clair avec cette version: Soient $(E,\leq)$ un ensemble ordonné et $x$, $y$ et $z$ des éléments de $E$. La relation $x\leq y$ est équivalente à ``$x< y\text{ ou }x=y$''. Chacune des relations ``$x\leq y\text{ et }y< z$'', ``$x< y\text{ et }y\leq z$'' entraîne $x<z$.

\vspace{0.5cm}
III, p. 14: Je crois que, pour qu'un ordre sur $E$ soit total, il faut et il suffit que son graphe $G$, en plus des relations $G\circ G= G$, $G\subset \Delta\circ G$, $G\subset G\circ \Delta$ et $G\cap G^{-1}=\Delta$, satisfasse à la relation $G\cup G^{-1}=E\times E$.

\section*{\S2. ENSEMBLES BIEN ORDONNÉS}

\vspace{0.5cm}
III, p. 16: Si $E$ un ensemble ordonné et $x\in E$, on dit aussi que $S_x$ est le \emph{segment initial} d'extémité (ou déterminé par) $x$.

\vspace{0.5cm}
III, p. 16, Prop. 3: Il faut ajouter une hypothèse: \dots pour tout couple d'indices $(\iota,\kappa)$ l'un des ensembles $X_\iota$ et $X_\kappa$ soit un segment de l'autre et l'ordre induit est identique à l'ordre donné.

\vspace{0.5cm}
III, p. 17: La démonstration de la prop. 3:
\begin{itemize}
\item La famille d'ensembles $(X_\iota)_{\iota\in I}$ vérifie les conditions du lemme 1. La première est évidente. Pour la deuxième, soit $(\iota,\kappa)$ un couple d'indices tel que $X_\iota\subset X_\kappa$. Si $X_\iota$ est un segment de $X_\kappa$, c'est évident. Si $X_\kappa$ est un segment de $X_\iota$, dans ce cas $X_\iota=X_\kappa$ et ces deux ensembles ont le même ordre. On conclut à l'aide de C18 (I, p. 28).

\item Montrons que si $x\in X_\iota$, $y\in E$ et $y\leq x$, alors $y\in X_\iota$. On a $y\in X_\kappa$ pour un certain $\kappa\in I$. Si $X_\kappa$ est un segment de $X_\iota$, on a évidemment $y\in X_\iota$. Si $X_\iota$ est un segment de $X_\kappa$, on a aussi $y\in X_\iota$. On conclut en utilisant C18 (I, p. 28). Donc chaque $X_\iota$ est un segment de $E$. Le fait que, pour tout $x\in X_\iota$, le segment d'extrémité $x$ dans $X_\iota$ est identique à l'intervalle $]\leftarrow,x[$ dans $E$ est une simple conséquence de la définition d'un segment (si $S$ est un segment d'un ensemble ordonné $E$ et $x\in S$, on a $]\leftarrow,x[_S=]\leftarrow,x[_E$).

\item Dans la démonstration du fait que $E$ est bien ordonné et pour monter que $a$ est aussi le plus petit élément de $H$ dans $E$, soit $x\in H$. On a $x\in X_\kappa$ pour un certain $\kappa\in I$. Si $X_\kappa$ est un segment de $X_\iota$, on a évidemment $x\in H\cap X_\iota$ et donc $a\leq x$. Si $X_\iota$ est un segment de $X_\kappa$. On ne peut avoir $x<a$ car si c'est le cas on aurait $x\in H\cap X_\iota$ et donc $a \leq x$ et $a<a$. Finalement $a\leq x$ car $X_\kappa$ est totalement ordonné (raisonnement par l'absurde, cf. C15 (I, p. 27)). On conclut à l'aide de C18 (I, p. 28).
\end{itemize}

\vspace{0.5cm}
III, p. 17,  Lemme 2: En vertu de la condition 1, l'ensemble $\mathfrak{S}$ doit vérifier $\varnothing=\bigcup\limits_{i\in \varnothing}S_i\in\mathfrak{S}$ (donc $\mathfrak{S}$ est non vide).

\vspace{0.5cm}
III, p. 18: C59 est équivalent à

C59$'$ Soient $E$ un ensemble bien ordonné et $A$ une partie de $E$. Si $$(\forall x)((x\in E\text{ et }S_x\subset A)\Rightarrow x\in A),$$ on a $A=E$.

En effet, C59 $\Rightarrow$ C59$'$. Il suffit de considérer le cas où $R\{x\}$ est la relation ``$x\in A$''.

C59$'$ $\Rightarrow$ C59. Il suffit de considérer le cas où $A=\{x\in E\mid R\{x\}\}$.

\vspace{0.5cm}
III, p. 18, Démonstration de C60: Toute réunion de segments appartenant à $\mathfrak{S}_1$ appartient à $\mathfrak{S}_1$: D'abord, si $S',S''$ deux segments de $E$ tels que $S''\in\mathfrak{S}_1$ et $S'\subset S''$, alors $S'\in\mathfrak{S}_1$ via $f_{S''}:S'\to f_{S''}(S')$. Clairement $\varnothing\in\mathfrak{S}_1$. D'autre part, soit $(S_i)_{i\in I}$ une famille d'éléments de $\mathfrak{S}_1$ telle que $I\neq\varnothing$. Soit l'application composée $S_i\xrightarrow {f_{S_i}}f_{S_i}(S_i)\hookrightarrow \bigcup\limits_{i\in I}f_{S_i}(S_i)$, pour tout $i$. D'après la remarque précédente, la seconde partie de C60 (l'unicité) et la prop.~7 de II, p. 28, ces applications s'étendent de façon unique à une application de $\bigcup\limits_{i\in I}S_i$ dans $\bigcup\limits_{i\in I}f_{S_i}(S_i)$. Donc $\bigcup\limits_{i\in I}S_i\in\mathfrak{S}_1$ via cette application.

\vspace{0.5cm}
III, p. 18: Une application de C60:

C60$'$ Soient $E$ un ensemble bien ordonné, $F$ un ensemble, $\Gamma$ l'ensemble des applications de segments initiaux de $E$ dans $F$ et $\varphi$ une application de $\Gamma$ dans $F$. Donc il existe une application unique $f:E\to F$ telle que pour tout $x\in E$ on ait $f(x)=\varphi(f|_{S_x})$.

On applique C60 et la considération qui se trouve juste après ce critère, en identifiant toute application $g:S_a\to F$ où $a\in E$, avec l'application surjective $g:S_a\to g(S_a)$. On prend $T\{u\}$ le terme $\varphi(u)=\bm{\tau}_y((u,y)\in G)$ où $G$ est le graphe de $\varphi$. L'unicité de $f$ se démontre d'une façon similaire à de C60.

Il existe aussi une démontration directe de C60$'$ analogue à celle de C60.

\vspace{0.5cm}
III, p. 19, Lemme 3:
\begin{itemize}
\item En outre on peut montrer facilement que la partie $M$ de $E$ et le bon ordre sur $M$ sont uniques à l'aide de la construction de $M$ (cf. les archives de N. Bourbaki, l'état 6 du chapitre III, p. 29).
\item Une autre application intéressante de ce lemme (autre que les théorèmes de Zermelo et de Zorn), le théorème de Bourbaki \cite[Exercice 4, p. 206]{Birkhoff: 1984}:

Soit $E$ un ensemble ordonné tel que toute partie bien ordonné admet une borne supérieure dans $E$. soit $f$ une application de $E$ dans $E$ telle que $f(x)\geq x$ pour tout $x\in E$. Sans assumer aucune forme de l'axiome de choix, $f$ admet un point fixe.

Soient $a\in E$ et $\mathfrak{S}$ l'ensemble formé de l'ensemble vide et des parties $X$ de $E$ contenant $a$, admettant une borne supérieure $m$ dans $E$, et telles que l'on ait $m\not\in X\text{ ou }f(m)>m$. Soit $p$ l'application de $\mathfrak{S}$ dans $E$ telle $p(\varnothing)=a$ et $$p(X)=
\begin{cases}m & \text{si}\quad m=\sup_EX\not\in X\\
f(m) & \text{si}\quad m=\sup_EX\in X.
\end{cases}$$
D'après C25 (I, p. 31), si $m=\sup_EX\in X$, alors $f(m)>m$, et par suite $f(m)\not\in X$. D'après le lemme 3 (III, p. 19), il existe une partie $M$ de $E$ et un bon ordre sur $M$ tels que 1 et 2. On montre, comme dans la démonstration de la prop. 4 de III, p. 20, que le bon ordre sur $M$ est l'ordre induit sur $M$ par l'ordre de $E$. On a $M\neq\varnothing$ car $M\notin\mathfrak{S}$. On a aussi $a\in M$, car dans le cas contraire il existerait $x\in M-\{a\}$. On aurait $S_x^M=\{y\in M\mid y<x\}\neq\varnothing$ car $p(S_x^M)=x$. Donc $a\in S_x^M$. Enfin, d'après les hypothèses, $M$ admet une borne supérieure $b$ vérifiant $b\in M$ et $f(b)=b$.
\end{itemize}

\vspace{0.5cm}
III, p. 20, Définition 3: Bien sûr un ensemble inductif $E$ est non vide car la partie vide de $E$ est totalement ordonnée (cf. III, p. 14, Exemple 2).

\vspace{0.5cm}
III, p. 21,  Th. 3: À l'aide de ce résultat et du th. de Zermelo, on obtient le théorème suivant: Soient $E$ et $F$ deux ensembles; l'une au moins des deux propositions suivantes est vraie:
\begin{enumerate}[1)]
  \item il existe une injection de $E$ dans $F$;
  \item il existe une injection de $F$ dans $E$.
\end{enumerate}

\vspace{0.5cm}
III, p. 22,  Produits lexicographiques: Soit $(E_i)_{i\in I}$ une famille d'ensembles ordonnés, dont l'ensemble d'indices $I$ soit bien ordonné. On appelle \emph{produit lexicographique} des $E_i$ l'ensemble $E=\prod\limits_{i\in I}E_i$ muni de la relation d'ordre: ``$x=(x_i)_{i\in I}\leq y=(y_i)_{i\in I}$'' si $$x=y \text{ ou } (x\neq y \text{ et } x_{i_0}<y_{i_0} \text{ où } i_0=\inf\{i\in I\mid x_i\neq y_i\}).$$ L'antisymétrie (la 2ème condition dans III, p. 1) est immédiate car si $x \leq y$ et $y \leq x$ et $x\neq y$ on aurait $x_{i_0}<y_{i_0}$ et $y_{i_0}<x_{i_0}$ où $i_0=\inf\{i\in I\mid x_i\neq y_i\}$, et donc $x_{i_0}=y_{i_0}$. Pour la transitivité, soient $x,y,z\in E$ tels que $x\leq y\leq z$ vérifiant $x\neq y$ et $y\neq z$. Soient $i_0=\inf\{i\in I\mid x_i\neq y_i\}$, $i_1=\inf\{i\in I\mid y_i\neq z_i\}$ et $i_2=\inf(i_0,i_1)$. Si $i<i_2$, on a $x_i=y_i=z_i$. En plus $x_{i_2}<z_{i_2}$. En effet, si $i_0=i_1$, on a $x_{i_0}<y_{i_0}<z_{i_0}$. Si $i_0<i_1$, on a $x_{i_0}<y_{i_0}=z_{i_0}$. Si $i_1<i_0$, on a $x_{i_1}=y_{i_1}<z_{i_1}$. On conclut à l'aide de C58 (III, p. 5) et de l'exerc. 3, b), \S 3, p. 48 (qui est une généralisation de C18 (I, p. 28)). Donc c'est bien une relation d'ordre sur $E$. Si en plus les $E_i$ sont totalement ordonnés, $E$ est évidemment totalement ordonné.

\section*{\S3. ENSEMBLES ÉQUIPOTENTS. CARDINAUX}

Soit $X$ un ensemble. On dit que $X$ est \emph{transitif} si la relation $x\in X$ implique $x\subset X$. On dit que $X$ est \emph{décent} si la relation $x\in X$ implique $x\not\in x$. L'ensemble vide est évidemment transitif et décent. Le \emph{successeur} de $X$ est l'ensemble $X\cup\{X\}$ qu'on note $X^+$. Noter que $X\subset X^+$ et $(X^+=X \Leftrightarrow X\in X)$. Si $X$ est décent, on a $X\not\in X$ et donc $X^+\neq X$. Le lemme suivant est évident.

\begin{lemma}\label{transitif}
Si $Y$ est un ensemble transitif (resp. décent), il en est de même de $Y^+$. Si $(Y_i)_{i\in I}$ est une famille d'ensembles transitifs (resp. décents), les ensembles $\bigcup\limits_{i\in I}Y_i$ et $\bigcap\limits_{i\in I}Y_i$ sont transitifs (resp. décents). La condition $I\neq\varnothing$ est implicite pour l'intersection. On rappelle que si $I=\varnothing$, $\bigcup\limits_{i\in I}Y_i=\varnothing$.
\end{lemma}

On dit que $X$ est un \emph{nombre ordinal} ou simplement un \emph{ordinal} si tout ensemble transitif $Y$ tel que $Y\subset X$ et $Y\neq X$ est un élément de $X$. Par exemple l'ensemble vide est un ordinal. Noter que l'ensemble vide appartient à tout ordinal non vide.

\begin{theorem}\label{ordinal}
\begin{enumerate}[(1)]
\item Tout ordinal est transitif et décent.
\item Si $X$ et $Y$ sont des ordinaux distincts, on a
\begin{enumerate}[a)]
\item $X\in Y\Leftrightarrow X\subset Y$.
\item $X\in Y\text{ ou }Y\in X$.
\end{enumerate}
\item Si $X$ est un ensemble transitif tel que tout $x\in X$ est un ordinal, alors $X$ est un ordinal.
\item Tout élément d'un ordinal $X$ est un ordinal.
\item \begin{enumerate}[a)]
\item Pour toute famille $(X_i)_{i\in I}$ d'ordinaux, l'intersection $\bigcap\limits_{i\in I}X_i$ est le plus petit élément de cette famille pour la relation d'inclusion. Soit $X$ est un ordinal. La relation $x\subset y$ entre éléments de $X$ est une relation de bon ordre.
\item \begin{enumerate}[i)]
\item Soit $X$ est un ordinal. Pour tout $x\in X$, $S_x=x$. En plus, une partie $Y$ de $X$ est un segment de $X$ si et seulement si $Y$ est transitif.
\item Si on a deux ordinaux, l'un est un segment de l'autre.
\item Deux ordinaux isomorphes sont identiques.
\end{enumerate}
\end{enumerate}
\item Si $X$ est un ordinal, il en est de même de $X^+$. Pour toute famille $(X_i)_{i\in I}$ d'ordinaux, la réunion $\bigcup\limits_{i\in I}X_i$ est un ordinal.
\end{enumerate}
\end{theorem}

(1) Soit $X$ un ordinal. Soient $(X_i)_{i\in I}$ la famille des sous-ensembles de $X$ transitifs et décents et $Y=\bigcup\limits_{i\in I}X_i$. D'après le lemme \ref{transitif}, $Y$ et $Y^+$ sont transitifs et décents. Supposons que $Y\neq X$. Donc $Y\in X$. Ce qui entraîne $Y^+\subset X$. Puisque $Y$ est décent, $Y\subsetneqq Y^+$, ce qui est absurde.

%Soit $Y\subsetneqq X^+$ et transitif. On distingue trois cas:
%\begin{itemize}
%\item $Y\subsetneqq X$. Donc $Y\in X\subset X^+$.
%\item $Y=X$. Or $X\in X^+$, donc $Y\in X^+$.
%\item $Y=A\cup\{X\}$ où $A\subsetneqq X$. On va montrer que $A$ est transitif. Soit $a\in A$. Donc $a\in Y$. Ce qui entraîne que $a\subset Y$.
%On a $X\not\in a$ car sinon, on aurait $X\in a\in X$, et donc $X\in X$ car $X$ est transitif, ce qui serait absurde ($X$ est décent). Finalement $a\subset A$. D'où $A\in X\subset X^+$ et $X\in X^+$.
%\end{itemize}

(2) Soient $X$ et $Y$ sont des ordinaux distincts. L'assertion i) est évidente: l'implication ``$\Rightarrow$'' puisque $Y$ est transitif, et l'implication ``$\Leftarrow$'' puisque $Y$ est un ordinal et $X$ est transitif.
Montrons ii). Supposons que $X\not\subset Y$ et $Y\not\subset X$. On a $X\cap Y\subsetneqq X$ et $X\cap Y\subsetneqq Y$. D'après (1) et le lemme \ref{transitif}, $X\cap Y$ est transitif et décent. Donc $X\cap Y\in X\cap Y$. Or $X\cap Y\not\in X\cap Y$ puisque $X\cap Y$ est décent. Ce qui constitue une contradiction. D'où $X\subset Y\text{ ou }Y\subset X$.

(3) Soit $X$ un ensemble transitif tel que tout $x\in X$ est un ordinal. Soient $Y\subsetneqq X$ transitif et $x\in X-Y$. On a $Y\subset x$. En effet, soit $y\in Y$. Supposons que $x\in y$. Donc $x\in Y$ car $Y$ est transitif, ce qui est absurde. D'après (2), $y\in x$. Si $Y=x$, on a évidemment $Y\in X$. Si $Y\neq x$, alors $Y\in x$ car $x$ est un ordinal et $Y$ est transitif. Donc $Y\in X$ car $X$ est transitif.

(4) La démonstration est similaire à celle de (1). Soit $X$ un ordinal. Soient $(X_i)_{i\in I}$ la famille des sous-ensembles de $X$ transitifs dont les éléments sont des ordinaux, et $Y=\bigcup\limits_{i\in I}X_i$. D'après le lemme \ref{transitif}, $Y$ et $Y^+$ sont transitifs et $Y$ est décent. D'après (3), $Y$ est un ordinal et les éléments de $Y^+$ sont des ordinaux. Supposons que $Y\neq X$. Donc $Y\in X$. Ce qui entraîne $Y^+\subset X$. Puisque $Y$ est décent, $Y\subsetneqq Y^+$, ce qui est absurde.

(5) a) Soit $(X_i)_{i\in I}$ une famille d'ordinaux. D'après (1) et le lemme \ref{transitif}, les $X_i$ et $\bigcap\limits_{i\in I}X_i$ sont transitifs et décents. Il suffit de montrer qu'il existe $i\in I$ tel que $\bigcap\limits_{i\in I}X_i=X_i$. En effet, supposons que pour tout $i\in I$, $\bigcap\limits_{i\in I}X_i\subsetneqq X_i$. Donc, pour tout $i\in I$, $\bigcap\limits_{i\in I}X_i\in X_i$. D'où $\bigcap\limits_{i\in I}X_i\in \bigcap\limits_{i\in I}X_i$. Or ceci est impossible car $\bigcap\limits_{i\in I}X_i$ est décent. La seconde partie découle de la première et de (4).

b) i) On a, pour tout $x\in X$, $S_x=\{y\in X\mid y<x\}=\{y\in X\mid y\in x\}=X\cap x=x$, car $X$ est transitif. La deuxième assertion découle de la première car, $Y$ est un segment $X$ veut dire que pour tout $x\in Y$, $S_x\subset Y$.

L'assertion ii) est évidente d'après i) et (2); et l'assertion iii) est évidente d'après ii) et le cor. 1 du th. 3, III, p. 21.

(6) est immédiate de (1), du lemme \ref{transitif}, de (4) et (3).

\begin{theorem}\label{caract. ordinal}
Soit $X$ un ensemble. Les conditions suivantes sont équivalentes:
\begin{enumerate}[1)]
\item $X$ est un ordinal.
\item La relation ``$(x\in y\text{ ou }x=y)\text{ et }x\in X\text{ et }y\in X$'' est une relation de bon ordre sur $X$, et $X$ est transitif.
\item $X$ est bien ordonné, et pour tout $x\in X$, $x$ est égal au segment d'extrémité $x$.
\end{enumerate}
\end{theorem}
$1)\Rightarrow 2)$ est contenue dans le th. \ref{ordinal}.

$2)\Rightarrow 3)$ est évidente.

$3)\Rightarrow 1)$. Soient $x,y\in X$. On a $(x<y)\Leftrightarrow (x\in S_y=y)$. Soit $x\in X$. On a $S_x=X\cap x=x$, c'est-à-dire, $x\subset X$. Donc $X$ est transitif.

Soit $Y\subsetneqq X$ transitif. On a $Y$ est un segment de $X$. D'après la prop. 1 de III, p. 16, il existe $z\in X$ tel que $Y=S_z=z$.

\begin{remark}
 N. Bourbaki a adopté au début la condition 3) du th. \ref{caract. ordinal} comme définition d'un ordinal (cf. les archives de N. Bourbaki, l'état 5 du chapitre III), et dans III, p. 79, exerc. 20, a appellé un ensemble vérifiant la condition 1) du th. \ref{caract. ordinal} un pseudo-ordinal.
\end{remark}

\begin{theorem}
\begin{enumerate}[1)]
\item Soit $E$ un ensemble d'ordinaux. $E$ muni de l'inclusion est bien ordonné.
\item Soit $\alpha$ un ordinal. L'ordinal $\alpha$ (resp. $\alpha^+$) est l'esemble des nombres ordinaux stritement plus petit (resp. plus petit) que $\alpha$.
\item Soit $E$ un ensemble d'ordinaux. La réunion $\bigcup\limits_{\alpha\in E}\alpha$ est la borne supérieure de $E$.
\item (Burali-Forti) Les ordinaux ne forment pas un ensemble.
\end{enumerate}
\end{theorem}
1) est évidente d'après la première partie de 5, a) du th. \ref{ordinal}.

2) est évidente d'après 4 et 6 du th. \ref{ordinal} et l'axiome d'extensionalité.

3) et 4) sont évidentes d'après 6 du th. \ref{ordinal}.

\begin{remarks}\label{prop. ordinaux}
Soient $\alpha$ et $\beta$ deux ordinaux.
\begin{enumerate}[i)]
\item Si $\beta<\alpha$, alors $\beta<\beta^+=\beta\cup\{\beta\}\leq\alpha<\alpha^+$. Donc $\alpha^+$ est le plus petit ordinal $>\alpha$.
\item Si $\beta<\alpha^+$, alors $\beta\leq\alpha<\alpha^+$.
\item Si $\alpha^+=\beta^+$ alors $\alpha=\beta$. En effet, supposons que $\alpha^+=\beta^+$ et $\alpha\neq\beta$. On a $\alpha,\beta\in\alpha^+=\beta^+$. Donc $\alpha\in\beta$ et $\beta\in\alpha$. D'où $\alpha\subset\beta$ et $\beta\subset\alpha$, c'est-à-dire, $\alpha=\beta$, ce qui est absurde.
\end{enumerate}
\end{remarks}

\begin{theorem}[Principe de récurrence sur les ordinaux]\label{réc. ordinaux}
Soit $R\{x\}$ une relation où figure la lettre $x$. Si pour tout ordinal $\alpha$, on a, $R\{\beta\}$ est vraie pour tout ordinal $\beta<\alpha$ entraîne $R\{\alpha\}$ est vraie. Alors $R\{\alpha\}$ est vraie pour tout ordinal $\alpha$.
\end{theorem}

Il suffit d'appliquer le principe de récurrence transfinie sur l'ensemble des ordinaux $\leq\xi$ pour un ordinal $\xi$ quelconque.

\begin{theorem}\label{type d'ordre}
Soit $X$ un ensemble bien ordonné. Il existe un ordinal unique isomorphe à $X$, qu'on appelle le \emph{type d'ordre} de $X$.
\end{theorem}

L'unicité découle de th. \ref{ordinal}, 5, b), iii).

L'existence. Soit $Y$ l'ensemble des $x\in X$ tels que $S_x$ est isomorphe à un ordinal. Cet ordinal est unique, notons le par $\alpha(x)$. Montrons que les $\alpha(x)$ forment un ensemble. Soient $x\in Y$ et $R$ la relation ``$y$ est un ordinal et $S_x$ est isomorphe à $y$''. On a $(\exists y)R$. Donc $\alpha(x)=\bm{\tau}_y(R)$. L'ensemble $\beta$ des objets de la forme $\alpha(x)$ pour $x\in Y$ montre notre assertion (cf. II, p. 6).

Maintenant, soit $f:E\to F$ un isomorphisme d'ensembles ordonnés. Si $A\subset E$ et $B=f(A)$, alors $f:A\to B$ est un isomorphisme. En particulier, pour tout $x\in E$, $f:S_x\to S_{f(x)}=f(S_x)$ est un isomorphisme.

D'après 3 du th. \ref{ordinal}, pour montrer que $\beta$ est un ordinal il suffit de vérifier qu'il est transitif. Soient $x\in Y$ et $\xi\in \alpha(x)$. Soit $f:S_x\to \alpha(x)$ l'isomorphisme canonique. Les segments initiaux d'extrémité $f^{-1}(\xi)$ dans $S_X$ et dans $X$ coincident. Donc $S_{f^{-1}(\xi)}$ est isomorphe à $S_\xi=\xi$. Par suite $f^{-1}(\xi)\in Y$ et $\alpha(f^{-1}(\xi))=\xi\in\beta$. D'où $\alpha(x)\subset\beta$.

Soient $x\in X$ et $y\in Y$ tels que $x<y$. Soit $f:S_y\to \alpha(y)$ l'isomorphisme canonique. Les segments initiaux d'extrémité $x$ dans $S_y$ et dans $X$ coincident, et $S_x\subsetneqq S_y$ car $x\in S_y-S_x$. Donc $f:S_x\to f(S_x)=S_{f(x)}=f(x)$ est un isomorphisme et $f(S_x)\subsetneqq\alpha(y)$. Donc $x\in Y$ et $\alpha(x)<\alpha(y)$. D'où $Y$ est un segment de $X$ et l'application $\alpha:Y\to \beta$ définie par $x\mapsto \alpha(x)$ est strictement croissante. D'après la prop. 11 de III, p. 14, l'application $\alpha$ est un isomorphisme.

Finalement, On a $Y=X$ car sinon on aurait $Y=S_x$ où $x\in X$ (d'après III, p. 16, Prop. 1), donc $x\in Y$, ce qui contredit le fait que $x\not\in S_x$.

\vspace{0.5cm}
\textbf{Les entiers naturels.} Maintenant on va définir les entiers naturels. On définit $0=\varnothing$, $1=\varnothing^+=\{0\}$, $2=1^+=\{0,1\}$, $3=2^+=\{0,1,2\}$,\dots, 95, 96=95$^+$,\dots.

On sait que l'ensemble vide est un ordinal. D'après 6 du th. \ref{ordinal}, les entiers naturels sont des ordinaux.

Soit $\alpha$ un ordinal $\neq0$. Soit $\beta$ la borne supérieure de l'ensemble des ordinaux $\xi<\alpha$. On a $\beta\leq\alpha$. Si $\beta=\alpha$, on dit que $\alpha$ est un \emph{ordinal limite}. Supposons que $\beta<\alpha$. Donc $\beta<\beta^+\leq\alpha$. On ne peut avoir $\beta^+<\alpha$; ce qui entraîne que $\beta^+=\alpha$. Dans ce cas, $\beta$ est appellé \emph{le prédécesseur} de $\alpha$. Réciproquement, supposons qu'il existe un ordinal $\gamma$ tel que $\gamma^+=\alpha$. D'après la remarque \ref{prop. ordinaux} ii), $\gamma$ la borne supérieure de l'ensemble des ordinaux $\xi<\alpha\;(\Leftrightarrow\xi\leq\gamma$) et il est le prédécesseur de $\alpha$. D'où, on a l'équivalence de
\begin{itemize}
  \item $\alpha$ est un ordinal limite;
  \item pour tout ordinal $\xi$, $\xi<\alpha\Rightarrow\xi^+<\alpha$;
  \item pour tout ordinal $\xi$, $\xi<\alpha\Rightarrow\xi^+\neq\alpha$
\end{itemize}
(car d'après la remarque \ref{prop. ordinaux} i), $\xi<\alpha\Rightarrow\xi^+\leq\alpha$).

Soit $n$ un entier $\neq0$. Le prédécesseur de $n$ existe. On le note par $n-1$. Si $n\geq2$, on note le prédécesseur de $n-1$ par $n-2$.

\vspace{0.5cm}
Par la suite on propose de calculer les nombres de signes et de liens figurant dans 1 et 2. On sait que 0 contient $12$ signes et $3$ liens (cf. II, p. 6). Soit $\mathcal{T}$ une théorie mathématique. Soient $\bm{A}$, $\bm{B}$ et $\bm{R}$ des assemblages et $\bm{x}$ une lettre. On a
\begin{itemize}
  \item $\bm{A}\Rightarrow \bm{B}$ est l'assemblage $\neg\,\bm{A}\,\bm{B}$.
  \item $\bm{A}\text{ et }\bm{B}$ est l'assemblage $\neg\,\vee\,\neg\,\bm{A}\,\neg\,\bm{B}$.
  \item $\bm{A}\Leftrightarrow\bm{B}$ est l'assemblage $\neg\,\vee\,\neg\,\vee\,\neg\,\bm{A}\,\bm{B}\,\neg\,\vee\,\neg\,\bm{B}\,\bm{A}$.
  \item $(\forall\bm{x})\bm{R}$ est l'assemblage $\neg\,\neg\,(\bm{\tau}_{\bm{x}}(\neg\bm{R})\mid\bm{x})\bm{R}$.
\end{itemize}

Revenons nous à la théorie des ensembles. Soient $\bm{R}$ la relation ``$z=x\text{ ou }z=y$'' et $\bm{S}$ la relation $(z\in t)\Leftrightarrow \bm{R}$. Donc $\bm{R}$ est l'assemblage $\vee=z\,x=z\,y$, et $\bm{S}$ est l'assemblage $$\neg\,\vee\,\neg\,\vee\,\neg\in z\,t\, \vee=z\,x=z\,y \,\neg\,\vee\,\neg \,\vee=z\,x=z\,y \in z\,t.$$
Donc $\bm{\tau}_z(\neg\bm{S})$ est l'assemblage $$\bm{\tau}\,\neg\,\vee\,\neg\,\vee\,\neg\in \square \,t \,\vee=\square \,x=\square \,y \,\neg\,\vee\,\neg \,\vee=\square \,x=\square \,y \in \square \,t$$
avec des liens entre $\bm{\tau}$ et les carrés, et $(\forall z)\bm{S}$ est

$\neg\,\neg\,\neg\,\vee\,\neg\,\vee\,\neg\in \bm{\tau}_z(\neg\bm{S}) \,t \,\vee=\bm{\tau}_z(\neg\bm{S}) \,x=\bm{\tau}_z(\neg\bm{S}) \,y \,\neg\,\vee\,\neg \,\vee=\bm{\tau}_z(\neg\bm{S}) \,x=\bm{\tau}_z(\neg\bm{S}) \,y \in \bm{\tau}_z(\neg\bm{S}) \,t$

dont figurent $24+(6×30)=204$ signes et $6×6=36$ liens.

Finalement $\{z\mid \bm{R}\}$, qui est par définition $\bm{\tau}_t(\forall z)\bm{S}$, contient 204+1=205 signes et $36+(2×6)+2=50$ liens avec $(6×4)+4=28$ occurences de $x$ et de $y$. D'où 1 contient $205-28+(28×12)=513$ signes et $50+(28×3)=134$ liens, car 0 contient $12$ signes et $3$ liens. Ce qui est acceptable pour 1. En outre 2 contient $205-28+(14×12)+(14×513)=7527$ signes et $50+(14×3)+(14×134)=1968$ liens.

\vspace{0.5cm}
Soient $X$ et $Y$ deux ensembles. On dit qu'ils sont \emph{équipotents} s'il existe une bijection de l'un sur l'autre. On dit que $X$ est un ensemble \emph{fini} s'il est équipotent à un entier naturel. Sinon on dit qu'il est \emph{infini}.

\vspace{0.5cm}
\textbf{L'axiome de l'infini.} Il existe un ensemble infini. Ce qui est équivalent, à l'aide des théorèmes de Zermelo et \ref{type d'ordre}, à ``il existe un ordinal infini''.

Soit $\alpha$ un ordinal infini. Les entiers naturels appartiennent à $\alpha$. L'intersection des ordinaux infinis $\xi\leq\alpha^+$ est un ordinal infini qu'on note $\N$ ou $\omega$. L'ordinal $\N$ est le plus petit ordinal infini. En effet, soit $\alpha'$ un ordinal infini. Si $\alpha\leq\alpha'$. Donc $\N\leq\alpha\leq\alpha'$. Si $\alpha'\leq\alpha<\alpha^+$. Donc $\N\leq\alpha'$. L'ordinal $\N$ est aussi l'ensemble des entiers naturels. En effet, les entiers naturels appartiennent à $\N$. Soit $\beta\in\N$. Donc $\beta$ est un ordinal et $\beta<\N\leq\alpha^+$. En plus $\beta$ est fini car dans le cas contraire on aurait $\N\leq\beta$. Alors $\beta$ est équipotent à un certain entier $n$. Si $n=0$, $\beta=0$. Supposons que $n\geq1$. On applique C60$'$ (cf. cette annexe). On prend $E=n=\{0,\dots,n-1\}$, $F=\beta$, et si $g:S_i\to\beta$ une application, $\varphi(g)$ est le plus petit élément de $\beta-g(S_i)$. Donc il existe une application $f:n\to\beta$ telle que $f(i)(=\varphi(f|S_i))$ est le plus petit élément de $\beta-f(S_i)$, pour tout $i\in n$. Si $n\geq2$ et d'après C59, on a $f(i)<f(i+1)$ pour tout $0\leq i\leq n-2$. Donc, à l'aide de la prop. 11 de III, p. 14, $f$ est un isomorphisme. D'où $\beta=n$.

Noter qu'on vient de démontrer aussi que tout ordinal fini est un entier. Les nombres ordinaux finis sont exactement les entiers. Les autres sont appellés \emph{transfinis}. L'ensemble $\N$ est le plus petit nombre ordinal transfini et aussi un ordinal limite. Puisque les entiers naturels ne sont pas des ordinaux limites alors un ordinal limite est nécessairement infini. Donc l'axiome de l'infini est équivalent à ``il existe un ordinal limite''.

\begin{theorem}[Le principe de récurrence]
Soit $R\{x\}$ une relation où figure la lettre $x$. Si $R\{0\}$ est vraie, et pour tout entier $n$, $R\{n\}\Rightarrow R\{n^+\}$ est vraie; alors $R\{n\}$ est vraie pour tout entier $n$.
\end{theorem}

En effet, supposons qu'il existe un entier $m$ tel que $\text{non }R\{m\}$ est vraie; et soit $m_0$ le plus petit des ces entiers. Donc $m_0\neq0$ car $R\{0\}$ est vraie. D'où $R\{m_0-1\}$ est vraie; ce qui entraîne que $R\{m_0\}$ est vraie, ce qui est absurde.

\begin{proposition}
Soit $n$ un entier naturel. Toute partie propre de $n$ est équipotente à un entier strictement plus petit que $n$.
\end{proposition}

On montre notre assertion par récurrence sur $n$. Le cas où $n=0$ est trivial car $0$ n'a pas de partie propre. Supposons que l'assertion est vrai pour $n$. Soit $A\subsetneqq n^+$. Si $A\subsetneqq n$, par l'hypothèse de récurrence, $A$ est équipotent à $m\in n\subset n^+$, donc $m\in n^+$. Si $A=\{n\}$, donc $A$ est équipotent à $1\in n^+$. Si $A=A'\cup\{n\}$ où $\varnothing\neq A'\subsetneqq n$. Par l'hypothèse de récurrence, $A'$ est en bijection avec $m\in n$. D'où $A$ est en bijection avec $m^+$ et $m^+\in n^+$ (cf. la remarque \ref{prop. ordinaux} i)).

\begin{corollary}
\begin{enumerate}[1)]
  \item Toute partie d'un ensemble fini est finie.
  \item Tout entier $n$ ne peut être équipotent à une partie propre.
\end{enumerate}
\end{corollary}

\begin{theorem}[Définition d'une application définie dans $\N$ par récurrence] Soient $F$ un ensemble non vide, $f:F\to F$ une application et $a\in F$. Il existe une suite unique $(u_n)_{n\in\N}$ d'éléments de $F$ vérifiant:
$$\begin{cases} u_0=a\in F\\
u_{n^+}=f(u_n) \quad n\geq 0.
\end{cases}$$
\end{theorem}
On applique C60$'$ (cf. cette annexe). On prend $E=\N$ et si $g:S_0=0\to F$ une application, $\varphi(g)= a\in F$. Si $g:S_{n^+}=n^+\to F$ une application où $n\geq 0$, $\varphi(g)= f(g(n))$.

En utilisant ce théorème, il est très facile de définir l'addition, la multiplication et l'exponentiation des entiers naturels, cf. \cite{Halmos:1960} ou bien \cite{Schwartz:1991}.

\begin{theorem}[Cantor]\label{cantor}
Soit $E$ un ensemble. Il existe une injection de $E$ dans $\mathfrak{P}(E)$ mais il n'y a pas de surjection (donc de bijection) de $E$ sur $\mathfrak{P}(E)$.
\end{theorem}

L'application $E$ dans $\mathfrak{P}(E)$ définie par $x\mapsto\{x\}$ est bien une injection. Soit $f$ une application de $E$ dans $\mathfrak{P}(E)$. On a $\{x\in E\mid x\not\in f(x)\}\not\in f(E)$; ce qui achève la démonstration du théorème.

\vspace{0.5cm}
Soit $X$ un ensemble. D'après le th. \ref{type d'ordre}, il existe un ordinal $\alpha$ équipotent à $X$. Montrons que les ordinaux équipotents à $X$ forment un ensemble. Soit $\gamma$ un ordinal équipotent à $\mathfrak{P}(X)$. On a $\gamma\not\subset\alpha$ car sinon il existerait une injection de $\mathfrak{P}(X)$ dans $X$ et donc une surjection de $X$ dans $\mathfrak{P}(X)$ (II, p. 18, Prop. 8), ce qui contredit le th. \ref{cantor}. D'où $\alpha\in\gamma$. On conclut à l'aide de C51 (II, p. 5).

Le plus petit des ordinaux équipotents à $X$ s'appelle le \emph{cardinal} de $X$ (aussi la \emph{puissance} de $X$) et on le note par $\card(X)$.

Un \emph{nombre cardinal} ou simplement un \emph{cardinal} est un nombre ordinal $\alpha$ tel que pour tout nombre ordinal $\beta$ équipotent à $\alpha$ on ait $\alpha\leq\beta$. Le cardinal d'un ensemble est bien un nombre ordinal. Inversement, si $\alpha$ un nombre cardinal, $\card(\alpha)=\alpha$.

Tout nombre cardinal infini $\alpha$ est un ordinal limite. En effet, supposons que le contraire. Donc il existe un ordinal infini $\beta$ tel que $\alpha=\beta^+=\beta\cup\{\beta\}$. On a $\N\subset\beta$ et $\beta\cup\{\beta\}=(\beta-\N^*)\cup (\N^*\cup\{\beta\})$ est équipotent à $(\beta-\N^*)\cup\N^*=\beta$. Ce qui entraîne que $\alpha$ et $\beta$ sont équipotents et $\alpha\leq\beta$; ce qui contredit $\beta<\beta^+$.

Soient $X$ et $Y$ deux ensembles. On a évidemment, $X$ et $Y$ sont équipotents si et seulement si $\card(X)=\card(Y)$.

Si $X$ est un ensemble fini, donc il est équipotent à un entier $n$. Donc $\card(X)=n$ car $n$ est le seul ordinal équipotent à $n$.

Si $X=\N$. Donc $\card(\N)=\N$ car $\N$ est le plus petit nombre ordinal transfini.

Soient $X$ et $Y$ deux ensembles. On a $\card(X)\subset\card(Y)$ si et seulement s'il existe une injection de $X$ dans $Y$. En effet, l'implication réciproque est évidente. Supposons que $\card(X)\not\subset\card(Y)$. Donc $\card(Y)\subsetneqq\card(X)$. Ce qui entraîne qu'il existe une injection de $Y$ dans $X$; et $X$ et $Y$ ne sont pas équipotents. D'après le th. de Cantor-Bernstein, il n'existe pas d'injection de $X$ dans $Y$.

Donc la relation $R\{\mathfrak{r},\mathfrak{n}\}$:

``$\mathfrak{r}\text{ et }\mathfrak{n}\text{ sont des cardinaux et }\mathfrak{r}\text{ est équipotent à une partie de } \mathfrak{n}$''
(cf. III, p. 24)

est équivalente à

``$\mathfrak{r}\text{ et }\mathfrak{n}\text{ sont des cardinaux et }\mathfrak{r}\subset\mathfrak{n}$''; et est une relation de bon ordre (III, p. 15). Notons $\mathfrak{r}\leq\mathfrak{n}$ la relation $R\{\mathfrak{r},\mathfrak{n}\}$.

Soit $\mathfrak{a}$ un cardinal. La relation ``$\mathfrak{r}\text{ un cardinal et }\mathfrak{r}\leq\mathfrak{a}$'' implique $\mathfrak{r}\in\mathfrak{a}^+$, donc elle collectivisante en $\mathfrak{r}$ (C52, II, p. 5). L'ensemble des $\mathfrak{r}$ vérifiant cette relation s'appelle l'\emph{ensemble des cardinaux} $\leq\mathfrak{a}$.

\vspace{0.5cm}
III, p. 24,  Démonstration du th. 1: on propose une méthode simple et directe pour montrer que ``$\mathfrak{r}\in E\text{ et }\mathfrak{n}\in E\text{ et }\mathfrak{r}\text{ est équipotent à une partie de } \mathfrak{n}$'' entraîne ``$\mathfrak{r}\in E\text{ et }\mathfrak{n}\in E\text{ et }\varphi(\mathfrak{r})\subset\varphi(\mathfrak{n})$''. En effet, soient $\mathfrak{r}\in E$ et $\mathfrak{n}\in E$. Supposons que $\mathfrak{r}$ est équipotent à une partie de $\mathfrak{n}$. Cette dernière relation est équivalente à ``$\mathfrak{r}$ est équipotent à une partie de $\varphi(\mathfrak{n})$''. Soit $Y$ une partie de $\varphi(\mathfrak{n})$ équipotent à $\mathfrak{r}$. Soit $S$ un segment de $\varphi(\mathfrak{n})$ équipotent à $Y$ (III, p. 22, Cor. 3). Or $\varphi(\mathfrak{n})$ est un segment de $A$, donc $S$ un segment de $A$ équipotent à $\mathfrak{r}$. D'où $\varphi(\mathfrak{r})\subset S\subset\varphi(\mathfrak{n})$.

Pour montrer que $R\{\mathfrak{r},\mathfrak{n}\}$ est bien une relation d'ordre entre $\mathfrak{r}$ et $\mathfrak{n}$, on considère les cas où $E$ est un ensemble à deux et à trois éléments (cf. l'axiome de l'ensemble à deux éléments (II, p. 4) et II, p. 26).

\vspace{0.5cm}
III, p. 25,  Prop. 3: Si $Y\neq\varnothing$, on a l'implication réciproque.

\vspace{0.5cm}
III, p. 25,  Définition 3: Si $I=\varnothing$, alors $\sum\limits_{i\in I}\mathfrak{a}_{i}=0$ et $\underset{i\in I}{\p}\mathfrak{a}_{i}=\card(\{\varnothing\})=\card(1)=1$.

\vspace{0.5cm}
III, p. 26,  Prop. 4: La seconde partie découle de la prop. \ref{som-fam-ens-bij} de l'annexe B. Le corollaire de la prop. 4 découle de la prop. 8 de II, p. 29, et de la prop. 3 de III, p. 25.

\vspace{0.5cm}
III, p. 26,  Prop. 5: On va donner les détails de la démonstration:
\begin{enumerate}[a)]
  \item $\sum\limits_{k\in K}\mathfrak{a}_{f(k)}=\card(\bigsqcup\limits_{k\in K}\mathfrak{a}_{f(k)})=\card(\bigcup\limits_{k\in K}\mathfrak{a}_{f(k)}\times\{f(k)\})=\card(\bigcup\limits_{i\in I}\mathfrak{a}_{i}\times\{i\})=\sum\limits_{i\in I}\mathfrak{a}_{i}$. (La seconde égalité d'après la prop. \ref{som-fam-ens} de l'annexe B; la troisième d'après la prop. 1 de II, p. 23.)

$\underset{k\in K}{\p}\mathfrak{a}_{f(k)}=\card(\prod\limits_{k\in K}\mathfrak{a}_{f(k)})=\card(\prod\limits_{i\in I}\mathfrak{a}_{i})=\underset{i\in I}{\p}\mathfrak{a}_{i}$. (D'après la prop. 4 de II, p. 33.)
  \item $\sum\limits_{i\in I}\mathfrak{a}_{i}=\card(\bigcup\limits_{i\in I}\mathfrak{a}_{i}\times\{i\})=\card(\bigcup\limits_{\lambda\in L}(\bigcup\limits_{i\in J_\lambda}\mathfrak{a}_{i}\times\{i\}))=\card(\bigsqcup\limits_{\lambda\in L}(\bigcup\limits_{i\in J_\lambda}\mathfrak{a}_{i}\times\{i\}))=\card(\bigsqcup\limits_{\lambda\in L}(\bigsqcup\limits_{i\in J_\lambda}\mathfrak{a}_{i}))=\card(\bigsqcup\limits_{\lambda\in L}(\sum\limits_{i\in J_\lambda}\mathfrak{a}_{i}))=\sum\limits_{\lambda\in L}(\sum\limits_{i\in J_\lambda}\mathfrak{a}_{i})$. (La seconde égalité d'après la prop. 2 de II, p. 24; la troisième d'après le cor. \ref{som-fam-ens-prop10} de l'annexe B; la cinquième d'après la prop. \ref{som-fam-ens-bij} de l'annexe B.)

$\underset{i\in I}{\p}\mathfrak{a}_{i}=\card(\prod\limits_{i\in I}\mathfrak{a}_{i})=\card(\prod\limits_{\lambda\in L}(\prod\limits_{i\in J_\lambda}\mathfrak{a}_{i}))=\card(\prod\limits_{\lambda\in L}(\underset{i\in J_\lambda}{\p}\mathfrak{a}_{i}))=\underset{\lambda\in L}{\p}(\underset{i\in J_\lambda}{\p}\mathfrak{a}_{i})$. (D'après la prop. 7 de II, p. 35, et le cor. de la prop. 11 de II, p. 38.)
\item $\underset{\lambda\in L}{\p}(\sum\limits_{i\in J_\lambda}\mathfrak{a}_{\lambda i})=\card(\prod\limits_{\lambda\in L}(\sum\limits_{i\in J_\lambda}\mathfrak{a}_{\lambda i}))=\card(\prod\limits_{\lambda\in L}(\bigcup\limits_{i\in J_\lambda}\mathfrak{a}_{\lambda i}\times\{i\}))=\card(\bigcup\limits_{f\in I}(\prod\limits_{\lambda\in L}\mathfrak{a}_{\lambda f(\lambda)}\times\{f(\lambda)\}))=\card(\bigsqcup\limits_{f\in I}(\prod\limits_{\lambda\in L}\mathfrak{a}_{\lambda f(\lambda)}\times\{f(\lambda)\}))=\card(\bigsqcup\limits_{f\in I}(\underset{\lambda\in L}{\p}\mathfrak{a}_{\lambda f(\lambda)}))=\sum\limits_{f\in I}(\underset{\lambda\in L}{\p}\mathfrak{a}_{\lambda f(\lambda)})$. (La seconde égalité d'après le cor. de la prop. 11 de II, p. 38; la troisième d'après la prop. 9 de II, p. 36; la quatrième d'après le cor. 1 de la prop. 9 de II, p. 36, et le cor. \ref{som-fam-ens-prop10} de l'annexe B; la cinquième d'après le cor. de la prop. 11 de II, p. 38, et la prop. \ref{som-fam-ens-bij} de l'annexe B.)
\end{enumerate}

\vspace{0.5cm}
III, p. 27,  Cor. de la prop. 5: Soient $\mathfrak{a}$, $\mathfrak{b}$, $\mathfrak{c}$ des cardinaux et $(\mathfrak{a}_i)_{i\in I}$ et $(\mathfrak{b}_j)_{j\in J}$ deux familles de cardinaux; on a

\begin{enumerate}[(1)]
\item $\mathfrak{a}+\mathfrak{b}=\mathfrak{b}+\mathfrak{a},\qquad\mathfrak{a}\mathfrak{b}=\mathfrak{b}\mathfrak{a},$
\item $\mathfrak{a}+\mathfrak{b}+\mathfrak{c}=(\mathfrak{a}+\mathfrak{b})+\mathfrak{c}=\mathfrak{a}+(\mathfrak{b}+\mathfrak{c}),\qquad \mathfrak{a}\mathfrak{b}\mathfrak{c}=(\mathfrak{a}\mathfrak{b})\mathfrak{c}=\mathfrak{a}(\mathfrak{b}\mathfrak{c}),$
\item $(\sum\limits_{i\in I}\mathfrak{a}_i)(\sum\limits_{j\in J}\mathfrak{b}_j)=\sum\limits_{(i,j)\in I\times J}\mathfrak{a}_i\mathfrak{b}_j$. En particulier, $\mathfrak{a}(\mathfrak{b}+\mathfrak{c})=\mathfrak{a}\mathfrak{b}+\mathfrak{a}\mathfrak{c}$.
\end{enumerate}

(1) est évidente d'après la prop. 5, a).

Pour la première égalité de la première partie de (2) il suffit d'appliquer la prop. 5, b) au cas où $I=3=\{0,1,2\}$, $\mathfrak{a}_0=\mathfrak{a}$, $\mathfrak{a}_1=\mathfrak{b}$, $\mathfrak{a}_2=\mathfrak{c}$, $L=\{0,1\}$, $J_0=\{0,1\}$ et $J_1=\{2\}$. La seconde égalité est évidente de (1) et de la première égalité. La seconde partie de (2) se démontre de la même façon.

Pour la première égalité de (3) il suffit d'appliquer la prop. 5, c) au cas où $L=\{0,1\}$, $J_0=I$ et $J_1=J$, $\mathfrak{a}_{0,i}=\mathfrak{a}_i$, $\mathfrak{a}_{1,j}=\mathfrak{b}_j$, pour tout $i\in I$ et $j\in J$, en remarquant qu'il existe une bijection canonique de $I=\prod\limits_{\lambda\in L}J_\lambda$ sur $J_0\times J_1$ (cf. II, p. 33), et en utilisant la prop.~5,~a). Le cas particulier est évident.

\vspace{0.5cm}
III, p. 29,  Cor. 3: On peut démontrer ce corollaire directement à l'aide de la prop. 3 de II, p. 31.

\vspace{0.5cm}
III, p. 30: Une autre démonstration de la prop. 14 est évidente à l'aide de la prop. \ref{som-fam-ens-bij} de l'annexe B et le cor. de la prop. 11 de II, p. 38.

\vspace{0.5cm}
III, p. 30, Cor. 2: On peut montrer autrement que ${\mathfrak{a}'}^\mathfrak{b}\leq{\mathfrak{a}'}^{\mathfrak{b}'}$: On écrit d'abord à l'aide de la prop. 13 (III, p. 29), $\mathfrak{b}'=\mathfrak{b}+\mathfrak{c}$ pour un cardinal $\mathfrak{c}$. Donc, d'après la prop. 10 (III, p. 28), son cor. 1, la prop. 14 (III, p. 30), la prop. 11 (III, p. 29) et le cor. 1 de la prop. 6 (III, p. 27), ${\mathfrak{a}'}^{\mathfrak{b}'}={\mathfrak{a}'}^\mathfrak{b}{\mathfrak{a}'}^\mathfrak{c}\geq {\mathfrak{a}'}^\mathfrak{b}1^\mathfrak{c}={\mathfrak{a}'}^\mathfrak{b}1={\mathfrak{a}'}^\mathfrak{b}$.


\begin{thebibliography}{99}

\bibitem{Birkhoff: 1984} G. Birkhoff, \emph{Lattice theory}, Amer. Math. Soc. Colloq. Publ. vol. 25; 3rd Revised edition, 1984.
\bibitem{Bourbaki:1970} N. Bourbaki, \emph{Théorie des ensembles}, Hermann, 1970, Springer-Verlag, 2006 (ses archives sont disponibles à http://mathdoc.emath.fr/archives-bourbaki/).
\bibitem{Cantor:1915} G. Cantor, \emph{Contributions to the founding of the theory of transfinite numbers}, Dover publications, 1915.
\bibitem{Froidevaux:1983} C. Froidevaux, \emph{La fonction $\varepsilon$ de Hilbert à travers les ``Grundlagen der Mathematik"}, Mathématiques et sciences humaines, tome \textbf{84} (1983), pp. 65--82 (disponible à www.numdam.org).
\bibitem{Godement:1997} R. Godement, \emph{Cours d'algèbre}, Hermann, 2$^\text{e}$ édition, 1997.
\bibitem{Grothendieck:1957} A. Grothendieck, \emph{Quelques points d'algèbre homologique}, Tôhoku Math. J. \textbf{9} (1957), 119--221.
\bibitem{Grothendieck/Verdier:1963/64} A. Grothendieck, J. -L. Verdier, \emph{Théorie des Topos et Cohomologie Étale des Schémas}, SGA 4, Lecture Notes in Math. \textbf{269}, Springer-Verlag (1972).
\bibitem{Grothendieck:1965} A. Grothendieck, \emph{Introduction au langage fonctoriel}, Séminaires 1965--1966, Faculté des Sciences d'Alger.
\bibitem{Halmos:1960} P. R. Halmos, \emph{Naive set theory}, D. Van Nostrand Company, Princeton 1960.
\bibitem{Hilbert/Bernays:1939}  D. Hilbert, P. Bernays, \emph{Grundlagen der mathematik II}, Second edition, Springer, 1970. Traduction française: \emph{Fondements des Mathématiques 2}, L'Harmattan, 2001.
\bibitem{Krivine:2007} J. -L. Krivine, \emph{Théorie des ensembles}, Cassini, 2$^\text{e}$ édition, 2007.
\bibitem{Manin:1977} Yu. I. Manin, \emph{A course in mathematical logic}, Springer-Verlag, 1977.
\bibitem{Mathias:2002} A. R. D. Mathias, \emph{A term of length 4 523 659 424 929}, Synthese \textbf{133} (2002), 75--86.
\bibitem{Schwartz:1991} L. Schwartz, \emph{Analyse I: Théorie des ensembles et topologie}, Hermann, 1991.
\bibitem{Tarski:1955} A. Tarski, \emph{A lattice-theoretical fixpoint theorem and its applications}. Pacific Journal of Mathematics, vol. \textbf{5} (1955), 285--309.
\end{thebibliography}
\end{document}